\documentstyle[12pt]{article}
\input amssym.def
\input amssym.tex
\parskip=\medskipamount
\newtheorem{guia}{}[section] 
\newtheorem{lemma}[guia]{Lemma} 
\newtheorem{example}[guia]{Example}
 
\newtheorem{proposition}[guia]{Proposition}
\newtheorem{corollary}[guia]{Corollary}
\newtheorem{theorem}[guia]{Theorem}
\newtheorem{remark}[guia]{Remark}

 1   
\font\ddpp=msbm10  scaled \magstep 1  

\def\QED{\hskip0.1em\hfill\null\ \null\nobreak\hfill
\kern3pt\lower1.8pt\vbox{\hrule\hbox   {\vrule\kern1pt\vbox{\kern1.7pt
\hbox{$\scriptstyle   QED$}\kern0.2pt}\kern1pt\vrule}\hrule}}
\def\R{\hbox{\ddpp R}}               
\def\Z{\hbox{\ddpp Z}}      

\def\arco#1{\kern3pt \mathop{\vbox{\ialign{##\crcr\noalign{\kern1pt}
        $\braceld\leaders\vrule\hfill\leaders\vrule\hfill\bracerd$
    \crcr\noalign{\kern1pt\nointerlineskip}
        $\hfil\displaystyle{\kern-1pt#1\kern2pt}\hfil$\crcr}}}\limits}

\def\lcf{\lbrack\! \lbrack}
\def\rcf{\rbrack\! \rbrack}

\begin{document}

\baselineskip=.55cm
\title{LICHNEROWICZ-JACOBI COHOMOLOGY AND HOMOLOGY OF JACOBI MANIFOLDS: MODULAR CLASS AND
DUALITY}
\author{\large Manuel  de Le\'on$^1$, Bel\'en  L\'opez$^2$, Juan C.
Marrero$^{3}$ and Edith  Padr\'on$^{3} $  
\\[7pt]
{\small \it $^1$Instituto de Matem\'aticas y F{\'\i}sica Fundamental,}\\[-8pt]
{\small\it Consejo Superior de In\-ves\-ti\-ga\-cio\-nes
Ci\-en\-t{\'\i}\-fi\-cas,} \\[-8pt]
{\small\it Serrano 123, 28006 Madrid, SPAIN,}\\[-8pt]
{\small\it E-mail: mdeleon@fresno.csic.es} \\
{\small\it $^2$Departamento de Matem\'aticas, Edificio de
Matem\'aticas e Inform\'atica,}\\[-8pt]
{\small \it Campus Universitario de Tafira, Universidad de Las Palmas
de Gran Canaria}\\[-8pt]
{\small \it 35017 Las Palmas, Canary Islands, SPAIN,}\\[-8pt]
{\small \it  E-mail: blopez@dma.ulpgc.es}\\
{\small\it $^3$Departamento de Matem\'atica Fundamental,
Facultad de Matem\'aticas,}\\[-8pt]
{\small\it Universidad de la Laguna,
La Laguna,} \\[-8pt]
{\small\it Tenerife, Canary Islands, SPAIN,}\\[-8pt]
{\small\it E-mail: jcmarrer@ull.es, mepadron@ull.es }
}
\date{\empty}

\maketitle

\begin{abstract}
Lichnerowicz-Jacobi cohomology and homology of Jacobi manifolds are
reviewed. We present both in a unified approach using the
representation of the Lie algebra of functions on itself by means of
the hamiltonian vector fields. The use of the associated Lie
algebroid allows to prove that the Lichnerowicz-Jacobi cohomology and
homology are invariant under conformal changes of the Jacobi 
structure and 
to stablish the duality between Lichnerowicz-Jacobi cohomology and
homology when the modular class vanishes. We also compute the
Lichnerowicz-Jacobi cohomology and homology for a large variety of examples.
\end{abstract}

\begin{quote}

\end{quote}

\newpage
\baselineskip=.5cm

\tableofcontents

\baselineskip=.55cm

\section{Introduction}

Since their introduction by Lichnerowicz in \cite{L1,L2}, Poisson and Jacobi
manifolds have deserved a lot of interest in the mathematical physics
literature. Indeed, the need to use more general phase spaces for
hamiltonian systems lead to the consideration of Poisson brackets of
non-constant rank, and, more than this, brackets which do not satisfy
Leibniz rule (Jacobi brackets). 

From the viewpoint of Differential Geometry, both structures are of
great interest. The local and global structures of Poisson and Jacobi
manifolds were ellucidated by Dazord, Guedira, Lichnerowicz, Marle, 
Weinstein and many others (\cite{DLM,GL,W}; see also
\cite{BV,liber,Vai2}). A Poisson manifold is basically made of
symplectic 
pieces, but the structure of a Jacobi manifold is more complicated, and
it is made of pieces which are contact or locally conformal symplectic
manifolds. 

The Poisson structure of a Poisson manifold $M$ allows  to define some
homology and cohomology operators. Indeed, the Poisson bivector of $M$
determines the so-called Lichnerowicz-Poisson cohomology
(LP-cohomology) and the $1$-differentiable Chevalley-Eilenberg
co\-ho\-mo\-lo\-gy, which can be alternatively described as the cohomologies
of two subcomplexes of the Chevalley-Eilenberg complex associated with
the Lie algebra of differentiable functions endowed with its Poisson
bracket (see \cite{L1}). The first of these subcomplexes consists of
the linear skew-symmetric multidifferential operators which are
derivations in each argument with respect to the usual product of 
functions, that is, the multivectors on $M.$ The second one consists
of the $1$-differentiable cochains, that is, the linear
skew-symmetric multidifferential operators of order $1$ in each
argument. Note that the space of $k$-cochains of this subcomplex is
isomorphic to ${\cal V}^k(M)\oplus {\cal V}^{k-1}(M),$ where ${\cal
V}^r(M)$ is the space of $r$-vectors on $M.$ 
Computation of Poisson cohomology is generally quite difficult. For
regular Poisson manifolds and for the Lie-Poisson structure on the
dual space of the Lie algebra of a compact Lie group, some results
were obtained in \cite{GLu,GW,V,Xu2}. 
On the other hand, we remark that  
the $k$-th LP-cohomology group has interesting
interpretations for the first few values of $k$. Moreover, these
cohomology groups allow to describe important results about the
geometric quantization and the deformation quantization of Poisson
manifolds (for more information, we refer to \cite{Vai2} and to the
recent survey \cite{W3}; see also the references therein).
The Poisson tensor of $M$ also allows to define the canonical
homology operator on forms (see \cite{Br,Ko}). The duality between
the canonical homology and the LP-cohomology is directly related to the
vanishing of the modular class introduced by Weinstein \cite{W2} (see
also \cite{BZ,ELW,Xu}). 

The situation for a Jacobi manifold $M$ is more involved. Note that
the Jacobi bracket of functions on $M$ is not a derivation in each
argument with respect to the usual product of functions (this is the
difference with the Poisson bracket). It
is only a linear skew-symmetric $2$-differential operator of order
$1$ or, in other words, a $1$-differentiable $2$-cochain in the
Chevalley-Eilenberg complex of the Lie algebra of functions. Thus, for
the manifold $M$, we have two possibilities. The first one is to
consider the representation of the Lie algebra of functions on itself
given by the Jacobi bracket. The resultant cohomology, the
Chevalley-Eilenberg cohomology of the Lie algebra of functions, was
studied by Guedira and Lichnerowicz \cite{GL,L2}. Particularly,
Guedira and Lichnerowicz studied the $1$-differentiable
Chevalley-Eilenberg cohomology, that is, the cohomology of the
subcomplex of the Chevalley-Eilenberg complex which consists of the
$1$-differentiable cochains. The second possibility is to consider
the representation of the Lie algebra of functions on itself given by
the action of the hamiltonian vector fields. The resultant cohomology
was termed by the authors, in \cite{LME3,LME4}, the
H-Chevalley-Eilenberg cohomology of $M.$ 
As in the case of the Chevalley-Eilenberg complex of $M$, one can
consider also the subcomplex of the $1$-differentiable cochains. The
cohomology of this subcomplex was termed the Lichnerowicz-Jacobi
cohomology (LJ-cohomology, for brevity) of $M$ (see \cite{LME3,LME4}).
If $M$ is a Poisson manifold then the Chevalley-Eilenberg cohomology and
the H-Chevalley-Eilenberg cohomology coincide and the $1$-differentiable
Chevalley-Eilenberg cohomology is just the LJ-cohomology. On the
other hand, the Lichnerowicz-Poisson complex of 
$M$ is isomorphic to a subcomplex of the Lichnerowicz-Jacobi complex.
The H-Chevalley-Eilenberg
cohomology and the LJ-cohomology  of a Jacobi manifold $M$ play an
important role in the geometric quantization of $M$ and in the study
of the existence of prequantization representations for complex line
bundles over $M$ (for more details, see \cite{LME3,LME4}).

The LJ-cohomology can be also described using the Lie algebroid
associated with the Jacobi manifold. Indeed, it is just the Lie algebroid
cohomology with trivial coefficients (see \cite{LME3,LME4,Pva}).
Using this fact we prove, in this paper, that the LJ-cohomology is
invariant under conformal changes of the Jacobi structure. Moreover,
we show that in many cases, the LJ-cohomology can be related with the
de Rham cohomology of the manifold, and in some cases it results a
topological invariant.

On the other hand, using again the representation of the Lie algebra
of functions on itself given by the hamiltonian vector fields, we
introduce the H-Chevalley-Eilenberg homology of a Jacobi manifold $M.$
The H-Chevalley-Eilenberg homology operator permits to define an
homology operator $\delta$ on the complex $\Omega^*(M)\oplus
\Omega^{*-1}(M),$ where $\Omega^k(M)$ is the space of $k$-forms
on $M,$ $\Omega^*(M)=\displaystyle\bigoplus_{k=1}^n\Omega^k(M)$ and $n$ is the
dimension of $M.$ The resultant homology is termed the
Lichnerowicz-Jacobi homology (LJ-homology).
If $M$ is a Poisson manifold the canonical homology of $M$ is
isomorphic to the homology of a subcomplex of the Lichnerowicz-Jacobi
complex.

It is easy to prove that the LJ-homology of a Jacobi
manifold $M$ is isomorphic to the Jacobi homology introduced by
Vaisman \cite{Pva}, which is described using the Lie algebroid associated
with $M.$ In fact, the Jacobi homology is the Lie algebroid homology
with respect to a flat connection on the top exterior power of the
jet bundle $J^1(M,\R).$ Using this result, we show, in this paper,
that the LJ-homology is invariant under conformal changes of the
structure. In \cite{Pva}, Vaisman also introduces the modular class
of $M$ as an element of the first LJ-cohomology group. Moreover, he
proves that if such a class is null (that is, $M$ is unimodular) then
the Lie algebroid cohomology with trivial coefficients and the Jacobi
homology are dual one each other. However, as we will show in this paper, there
exist important examples of Jacobi manifolds which are not
unimodular. In these cases we will describe the LJ-homology and we will
conclude that in most of them, the LJ-cohomology and the
LJ-homology are not, in general, dual one each other.

The paper is organized as follows.

Section 2 is introductory, and it contains some generalities about
Jacobi manifolds: definitions, examples and the construction of the
Lie algebroid canonically associated to any Jacobi manifold.
In Section 3, we first recall the notions of H-Chevalley-Eilenberg
cohomology and LJ-cohomology of a Jacobi manifold $(M,\Lambda,E)$. 
In particular, we recall that the LJ-cohomology can be described as
the cohomology of the Lie algebroid associated with $M$ (see
\cite{LME3,LME4,Pva}) or alternatively introducing the 
operator $\sigma:{\cal V}^k(M)\oplus {\cal V}^{k-1}(M)\longrightarrow
{\cal V}^{k+1}(M)\oplus {\cal V}^{k}(M)$, given by 
$\sigma(P,Q)=(-[\Lambda,P]+kE\wedge P+\Lambda\wedge
Q,[\Lambda, Q]-(k-1)E\wedge Q+[E,P])$, where $[\;,\;]$ is the
Schouten-Nijenhuis bracket (see \cite{LME3,LME4}).
The first description allows us to prove that 
the LJ-cohomology is invariant under conformal changes of the Jacobi structure. 
The rest of the section is devoted to study the LJ-cohomology for
different examples of Poisson structures (symplectic and Lie-Poisson
structures and a quadratic Poisson structure on $\R^2$) and of
Jacobi structures. In particular, the 
LJ-cohomology of contact and locally conformal symplectic manifolds
is extensively studied. 
We also consider another interesting example: the Jacobi
structure of the unit sphere of a real Lie algebra of finite
dimension. In Table I, we summarize the main results obtained about
the LJ-cohomology of the different examples of Jacobi manifolds.

In Section 4, we introduce and study the LJ-homology of a Jacobi
manifold $(M,\Lambda,E).$ First of
all, the H-Chevalley-Eilenberg homology of $M$ is
defined as the homology of the Lie algebra of functions
on $M$ with respect to the representation given by the
hamiltonian vector fields. Since every $k$-chain of the
H-Chevalley-Eilenberg complex defines a pair $(\alpha,\beta)$, with
$\alpha$ a $k$-form and $\beta$ a $(k-1)$-form, the
H-Chevalley-Eilenberg boundary  operator induces a
boundary operator 
$\delta: \Omega^k(M)\oplus \Omega^{k-1}(M)\longrightarrow
\Omega^{k-1}(M)\oplus \Omega^{k-2}(M)$ given by
$\delta(\alpha,\beta) =
(i(\Lambda)d\alpha-di(\Lambda)\alpha+ki_E\alpha + (-1)^k{\cal L}_E\beta,
i(\Lambda)d\beta-di(\Lambda)\beta+(k-1)i_E\beta+(-1)^ki(\Lambda)\alpha)$,
where $i(\Lambda)$ denotes the contraction by $\Lambda$ and ${\cal L}$ is
the Lie derivative operator.
The resultant homology is called Lichnerowicz-Jacobi homology (LJ-homology).
As we have indicated above, there is an alternative description of
the LJ-homology which
is due to Vaisman \cite{Pva}. 
This description allows us to prove an  important first result: the
invariance of the 
LJ-homology under conformal changes of the Jacobi structure.
The rest of the section is devoted to study the LJ-homology of
the different examples of Jacobi manifolds considered in Section 3.
We show that the symplectic manifolds, the dual space of a unimodular
real Lie algebra ${\frak g}$ (endowed with the Lie-Poisson structure)
and the unit sphere on ${\frak g}$ are unimodular Jacobi manifolds.
Thus, using the results of Section 3 and the results of Vaisman
\cite{Pva} on the duality between the LJ-cohomology and the
LJ-homology, we describe the LJ-homology of the above examples of
Jacobi manifolds (at least for the case when ${\frak g}$ is the Lie
algebra of a compact Lie group). On the other hand, the contact
manifolds and the locally (non-globally) conformal symplectic
manifolds are not unimodular Jacobi manifolds. In fact, we deduce
that in these cases the LJ-homology and the LJ-cohomology are not, in
general, dual one each other. In Table II, we summarize the main
results obtained on the LJ-homology (and its relation with the
LJ-cohomology) of the different examples of Jacobi manifolds.

\section{Jacobi manifolds}

All the manifolds considered in this paper are assumed to be
connected. Moreover, if $M$ is a differentiable manifold, we will use
the following notation:
\begin{itemize}
\item
$C^{\infty}(M, \R)$ is the algebra of $C^{\infty}$ real-valued
functions on $M$.
\item
$\frak X(M)$ is the Lie algebra of the vector fields on $M$.
\item
$\Omega^{k}(M)$ is the space of $k$-forms on $M$.
\item
${\cal V}^{k}(M)$ is the space of $k$-vectors on $M$.
\end{itemize}

\subsection{Local Lie algebras and Jacobi manifolds}

A {\it Jacobi structure } on a $n$-dimensional manifold $M$ is a pair
$(\Lambda, E)$ where $\Lambda$ is a $2$-vector  and $E$ a vector
field on $M$ satisfying the following properties:
\begin{equation}
\label{e1}
[\Lambda, \Lambda] = 2 E \wedge \Lambda \; , \makebox[1cm]{}
{\cal L}_E\Lambda =[E,\Lambda] = 0 \; .
\end{equation}
Here $[\;,\;]$ denotes the Schouten-Nijenhuis bracket
(\cite{BV,Vai2})
and ${\cal L}$ is the Lie derivative operator. The
manifold $M$ endowed with a Jacobi structure is called a {\it Jacobi
manifold}. A bracket of functions (the {\it Jacobi bracket}) is
defined by
\begin{equation}\label{cf}
\{f,g\} = \Lambda(df,dg) + fE(g) -gE(f) \; ,\makebox[2cm]{for all}
f,g\in C^\infty(M,\R).
\end{equation}

The Jacobi bracket $\{\;,\;\}$ is skew-symmetric, satisfies the Jacobi
identity and, in addition, we have 
\[
{\rm support }\{f,g\}\subseteq ({\rm support }\; f)\cap ({\rm support
 }\; g ) \; ,\makebox[2cm]{for all} f,g\in C^\infty(M,\R).
\]

Thus, the space $C^{\infty}(M,\R)$ endowed with the Jacobi bracket is
{\it a local Lie algebra} in the sense of Kirillov (see \cite{Ki}).
Conversely, a structure of local Lie algebra on $C^{\infty}(M,\R)$
defines a Jacobi structure on $M$ (see \cite{GL,Ki}).
If the vector field $E$ identically vanishes then $\{\; , \;\}$ is a
derivation in each argument and, therefore, $\{ \;,\;\}$ defines
a {\it Poisson bracket } on $M$. In this case, (\ref{e1}) reduces to
$[\Lambda,\Lambda] = 0$
and $(M,\Lambda)$ is a {\it Poisson
manifold.} Jacobi and Poisson manifolds were introduced by
Lichnerowicz (\cite{L1,L2}).

\begin{remark}\label{2.0}{\rm Let $(\Lambda,E)$ be a Jacobi structure
on a manifold $M$ and consider on the product manifold $M\times \R$
the $2$-vector $\tilde\Lambda$ given by
\[
\tilde{\Lambda} = e^{-t}(\Lambda+\frac{\partial}{\partial t} \wedge E),
\]
where $t$ is the usual coordinate on $\R.$ Then, $\tilde{\Lambda}$
defines a Poisson structure on $M\times \R$ (see \cite{L2}). The
manifold $M\times \R$ endowed with the structure $\tilde\Lambda$ is
called the {\it poissonization of the Jacobi manifold } $(M,\Lambda,E)$.
}
\end{remark}

\subsection{Examples of Jacobi manifolds}

In this section, we will present some examples of Jacobi manifolds.

\medskip

1. {\it Symplectic manifolds}.- A {\it symplectic manifold} is a pair
$(M, \Omega)$, where  $M$ is an even-dimensional manifold and
$\Omega$ is a closed non-degenerate $2$-form on $M$. We define a
Poisson $2$-vector $\Lambda$ on $M$ by
\begin{equation} \label{Posy}
\Lambda(\alpha,\beta) =
\Omega(\flat^{-1}(\alpha),\flat^{-1}(\beta)),
\end{equation}
for all $\alpha, \beta \in \Omega^{1}(M)$, where
$\flat : \frak X(M) \longrightarrow \Omega^1(M)$ is the
isomorphism of $C^{\infty}(M,\R)$-modules given by
$\flat(X) = i_{X}\Omega$ (see \cite{L1}).

Using the classical theorem of Darboux, around every point of $M$
there exist canonical coordinates $(q^1,\dots ,q^m,p_1,\dots ,p_m)$
on $M$ such that
\[
\Omega=\sum_idq^i\wedge dp_i,\makebox[1cm]{}
\Lambda=\sum_i\frac{\partial}{\partial q^i}\wedge \frac{\partial
}{\partial p_i}.
\]

\smallskip

2. {\it Lie-Poisson structures}.- Let $({\frak g},[\;,\;])$ be a real
Lie algebra of dimension $n$ with Lie bracket $[\;,\;]$ and denote by
${\frak g}^{*}$ the dual vector space of ${\frak g}$. Given two
functions $f, g \in C^{\infty}({\frak g}^{*}, \R)$, we define $\{f,
g\}$ as follows. For a point $x \in {\frak g}^{*}$, we linearize $f$
and $g$, namely, we take the tangent maps $df(x)$ and $dg(x)$ at $x$
and identify them with two elements $\hat{f}, \hat{g} \in {\frak g}$.
Thus, $[\hat{f}, \hat{g}] \in {\frak g}$, and we define
\[
\{f, g\}(x) = [\hat{f}, \hat{g}].
\]
$\{\;,\;\}$ is the so-called {\it Lie-Poisson bracket} on ${\frak g}^{*}$
(see \cite{Vai2,W}).

If $\bar{\Lambda}$ is the corresponding Poisson $2$-vector on ${\frak
g}^{*}$ and $(x_{i})$ are global coordinates for ${\frak g}^{*}$
obtained from a basis, we have
\begin{equation}\label{LieP}
\bar\Lambda =
\frac{1}{2}\sum_{i,j,k}c_{ij}^kx_k\frac{\partial}{\partial x_i}\wedge
\frac{\partial}{\partial x_j},
\end{equation}
$c_{ij}^k$ being the structure constants of ${\frak g}$.

From (\ref{LieP}), it follows that
\begin{equation}\label{LieP1}
{\cal L}_{A}\bar\Lambda = -\bar\Lambda,
\end{equation}
where $A$ is the radial vector field on ${\frak g}^{*}$. Note that
the expression of $A$ with respect to the coordinates $(x_{i})$ is
\begin{equation}\label{LieP2}
A = \displaystyle \sum_{i} x_{i} \frac{\partial}{\partial x_{i}}.
\end{equation}

\smallskip

3. {\it Contact manifolds}.-
Let $M$ be a $(2m+1)$-dimensional manifold and $\eta$ a $1$-form
on $M$. We say that $\eta$ is a contact $1$-form if $\eta
\wedge (d\eta)^{m} \neq 0$ at every point. In such a case
$(M,\eta)$ is termed a {\it contact manifold} (see, for example,
\cite{Bl,liber,L2}). If $(M, \eta)$ is a contact
manifold, we define a $2$-vector $\Lambda$ and a vector $E$ on $M$ as
follows 
\begin{equation} \label{Jacon}
\Lambda(\alpha,\beta) =
d\eta(\flat^{-1}(\alpha),\flat^{-1}(\beta)), \makebox[.4cm]{} 
E = \flat^{-1}(\eta)
\end{equation}
for all $\alpha, \beta \in \Omega^{1}(M)$, where
$\flat: \frak X(M) \longrightarrow \Omega^1(M)$ is the
isomorphism of $C^{\infty}(M,\R)$-modules given  by $\flat(X) =
i_{X}d\eta + \eta(X) \eta$. Then $(M, \Lambda, E)$ is a Jacobi manifold.
The vector field $E$ is just the {\it Reeb vector field } of $M$ and it is
characterized by the relations
\begin{equation}\label{Reeb}
i_{E} \eta = 1, \makebox[.4cm]{} i_{E}d\eta = 0.
\end{equation}
Using the generalized Darboux theorem, we deduce that around every
point of $M$ there exist canonical coordinates
$(t, q^{1}, \dots , q^{m}, p_{1}, \dots , p_{m})$ such that (see
\cite{liber,L2})
\begin{equation}\label{6''}
\eta = dt - \displaystyle \sum_{i} p_{i}dq^{i}, \makebox[1cm]{}
\Lambda = \displaystyle \sum_{i} (\frac{\partial}{\partial q^{i}} +
p_{i} \frac{\partial}{\partial t}) \wedge \frac{\partial}{\partial
p_{i}}, \makebox[1cm]{}E = \displaystyle \frac{\partial}{\partial t}.
\end{equation}

\begin{remark}\label{2.0'}{\rm
The poissonization of a contact structure is a symplectic structure
(see \cite{L2}).}
\end{remark}
\smallskip

4. {\it Locally conformal symplectic manifolds}.-
An {\it almost symplectic manifold} is a pair $(M,\Omega)$, where $M$
is an even dimensional manifold and $\Omega$ is a non-degenerate 2-form on
$M$. An almost symplectic manifold is
said to be {\it locally conformal symplectic (l.c.s.)}
if for each point $x \in M$ there is an open neighborhood $U$ such that
$d(e^{-f}\Omega) = 0$, for some function $f: U
\longrightarrow \R$ (see, for example, \cite{GL,Vai1}). So, $(U,=
 e^{-f}\Omega)$ is a symplectic manifold. If
$U=M$ then $M$ is said to be a {\it globally conformal
symplectic (g.c.s.)} manifold. An almost symplectic manifold
$(M,\Omega)$ is l.(g.)c.s.
if and only
if there exists a closed (exact) 1-form $\omega$ such that
\begin{equation}\label{conLcs}
d\Omega = \omega \wedge \Omega.
\end{equation}
The 1-form $\omega$ is called the {\it Lee 1-form} of $M$.
It is obvious that the l.c.s. manifolds with Lee 1-form
identically zero are just the symplectic manifolds.

In a similar way that for contact manifolds, we define a $2$-vector
$\Lambda$ and a vector field $E$ on $M$ which are given by
\begin{equation}\label{Lcs}
\Lambda(\alpha,\beta) =
\Omega(\flat^{-1}(\alpha),\flat^{-1}(\beta)), \makebox[1cm]{} E =
\flat^{-1}({\omega}) \; ,
\end{equation}
for all $\alpha , \beta \in \Omega^{1}(M)$, where $\flat :
\frak{X}(M) \longrightarrow \Omega^1(M)$ is the
isomorphism of $C^{\infty}(M,\R)$-modules defined by $\flat
(X) = i_{X}\Omega$. Then $(M,\Lambda, E)$ is a Jacobi
manifold (see \cite{GL}). Note that
\begin{equation}\label{Inv}
{\cal L}_{E}\Omega = 0.
\end{equation}
Using the
classical theorem of Darboux, around every point of $M$ there exist
canonical coordinates $(q^{1}, \dots , q^{m}, p_{1}, \dots ,p_{m})$
and a local differentiable function $f$ such that
\[
\begin{array}{lcl}
\Omega = e^{f} \displaystyle \sum_{i} dq^{i} \wedge dp_{i},&&
\omega = df = \displaystyle \sum_{i}(\frac{\partial
f}{\partial q^{i}} dq^{i} + \frac{\partial f}{\partial
p_{i}} dp_{i}),\\[5pt]
\Lambda = e^{-f} \displaystyle \sum_{i}
(\frac{\partial}{\partial q^{i}} \wedge \frac{\partial}{\partial
p_{i}}), &&
E = e^{-f} \displaystyle \sum_{i}(\frac{\partial
f}{\partial p_{i}} \frac{\partial}{\partial q^{i}} -
\frac{\partial f}{\partial q^{i}} \frac{\partial}{\partial
p_{i}}).
\end{array}
\]

\smallskip

5. {\it Unit sphere of a real Lie algebra}.- Let ${\frak g}$ be a
real Lie algebra of dimension $n$ with Lie bracket $[\;,\;]$ and let
$\bar{\Lambda}$ be the Poisson $2$-vector on the dual vector space
${\frak g}^{*}$ of ${\frak g}$.

Suppose that $<\;,\;>$ is a scalar product on ${\frak g}$ and that
$g$ is the corresponding Riemannian metric on ${\frak g}$.

Denote by $\flat_{<\;,\;>}:{\frak g}\rightarrow {\frak g}^*$ the
linear isomorphism between ${\frak g}$ and ${\frak g}^*$ given by
\begin{equation}\label{2.12'}
\flat_{<\;,\;>}(\xi)(\eta)=<\xi,\eta>,\makebox[1cm]{}\mbox{for all
}\xi,\eta\in {\frak g}.
\end{equation}
Using this isomorphism and the Lie-Poisson structure $\bar{\Lambda}$,
we can define a Poisson structure on ${\frak g}$ which we also denote
by $\bar{\Lambda}.$

Now, we consider the $2$-vector $\Lambda'$ and the vector field $E'$
on ${\frak g}$ given by
\begin{equation}\label{LinJ}
\Lambda'=\bar\Lambda-A\wedge i_\alpha\bar\Lambda,\makebox[1cm]{}
E'=i_\alpha\bar\Lambda,
\end{equation}
where $A$ is the radial vector field on ${\frak g}$ and $\alpha$ is
the $1$-form defined by $\alpha(X) = g(X, A)$, for $X \in {\frak
X}({\frak g})$. From (\ref{LieP2}), we obtain that
\begin{equation}\label{LinJ1}
\alpha = \displaystyle \frac{1}{2} d(\| \; , \; \|^{2}),
\makebox[.4cm]{} {\cal L}_{A}\alpha = 2\alpha,
\end{equation}
$\| \; , \; \|^{2}: {\frak g} \longrightarrow \R$ being the real
function on ${\frak g}$ given by
\[
\| \; , \; \|^{2}(\xi) = <\xi, \xi>,
\]
for all $\xi \in {\frak g}$. Using (\ref{LieP1}) and (\ref{LinJ1}), we
deduce that
\begin{equation}\label{EA}
[A, E'] = E', \makebox[.5cm]{} {\cal L}_{E'}\bar\Lambda = \displaystyle
\frac{1}{2} [[\bar\Lambda, \| \; , \; \|^{2}], \bar\Lambda]
= 0.
\end{equation}
Thus, the pair $(\Lambda', E')$ induces a Jacobi
structure on ${\frak g}$. Moreover, if $S^{n-1}({\frak g})$ is the
unit sphere in ${\frak g}$, it follows that the restrictions
${\Lambda}$ and $E$ to $S^{n-1}({\frak g})$ of $\Lambda'$ and $E'$,
respectively, are tangent to $S^{n-1}({\frak g})$. Therefore, the pair
$(\Lambda, E)$ defines a Jacobi structure on $S^{n-1}({\frak g})$
(for more details, see \cite{L3}).  

If for every $\xi \in {\frak g}$, we consider the function $<\xi, \;\;>:
S^{n-1}({\frak g}) \longrightarrow \R$ given by
\begin{equation}\label{Fusc}
<\xi, \;\; >(\eta) = <\xi, \eta>,
\end{equation}
then, from (\ref{LieP}), (\ref{LinJ}) and (\ref{Fusc}), we have that
\begin{equation}\label{CompL}
\{<\xi, \;\;>, <\eta, \;\;>\} = <[\xi, \eta], \;\;>
\end{equation}
for $\xi, \eta \in {\frak g}$, where $\{\;,\;\}$ is the Jacobi
bracket on $S^{n-1}({\frak g})$. 

Note that if $\xi\in S^{n-1}({\frak
g})$ it follows that $(d<\xi,\;>)(\xi)=0$ and consequently 
\[
E(\xi)=X_{<\xi, \;>}(\xi).
\]
This implies that the characteristic foliation of $S^{n-1}({\frak
g})$ is generated by the set of hamiltonian vector fields
$\{X_{<\xi,\;>}/ \xi\in {\frak g}\}.$

On the other hand, 
if $(x^{i})$ are global coordinates for ${\frak g}$
obtained from an orthonormal basis $\{\xi_{i}\}_{i=1,\dots,n}$ of
${\frak g}$ then
\[
\Lambda'=\sum_{i,j,k,h,r}(\frac{1}{2}c_{hj}^rx^r-c_{ij}^kx^kx^ix^h)
\frac{\partial}{\partial x^h}\wedge
\frac{\partial}{\partial x^j},\makebox[1cm]{}
E'=\sum_{i,j,k}c_{ij}^kx^kx^i\frac{\partial}{\partial
x^j},
\]
$c_{ij}^k$ being the structure constants for ${\frak g}$ with respect
to the basis $\{\xi_i\}_{i=1,\dots,n}$.

\begin{remark}\label{2.0''}{\rm 
Using the results of \cite{L3}, we
obtain that the poissonization of the Jacobi ma\-ni\-fold
$(S^{n-1}({\frak g}),\Lambda,E)$ is isomorphic to the Poisson
manifold $({\frak  g}-\{0\},\bar\Lambda_{|{\frak g}-\{0\}}).$ In fact,
an isomorphism between these Poisson manifolds is defined by 
\begin{equation}\label{2.17'}
F:{\frak g}-\{0\}\rightarrow S^{n-1}({\frak g})\times
\R,\makebox[1cm]{}\xi\mapsto
F(\xi)=(\frac{\xi}{\|\xi\|},\mbox{ln}\|\xi\|).
\end{equation}}
\end{remark}

\subsection{The characteristic foliation of a Jacobi manifold}

Let $(M,\Lambda,E)$ be a Jacobi manifold. Define a homomorphism of
$C^{\infty}(M, \R)$-modules $\#_{\Lambda}:\Omega^1(M) \longrightarrow
{\frak X}(M)$ by
\begin{equation}\label{e7}
(\#_{\Lambda}(\alpha))(\beta)=\Lambda(\alpha,\beta),
\end{equation}
for $\alpha,\beta\in \Omega^1(M)$.
This homomorphism can be extended to a homomorphism, which we also denote by
$\#_{\Lambda}$, from the space $\Omega^k(M)$ onto the space
${\cal V}^k(M)$ by putting:
\begin{equation}\label{e8}
\#_{\Lambda}(f)=f,\makebox[.5cm]{}
\#_{\Lambda}(\alpha)(\alpha_1,\dots,\alpha_k)
=(-1)^k\alpha(\#_{\Lambda}(\alpha_1),
\dots ,\#_{\Lambda}(\alpha_k)),
\end{equation}
for $f\in C^\infty(M,\R),$ $\alpha \in \Omega^k(M)$ and
$\alpha_1,\dots ,\alpha_k\in \Omega^1(M).$

\begin{remark}\label{r0}{\rm $i)$ If $M$ is a contact manifold with
Reeb vector field $E$, then $\#_{\Lambda}(\alpha) = -\flat
^{-1}(\alpha) + \alpha (E)E$, for all $\alpha \in \Omega^{1}(M)$.

$ii)$ If $M$ is a l.c.s. manifold then $\#_{\Lambda}(\alpha) =
-\flat^{-1}(\alpha)$, for all $\alpha \in \Omega^{1}(M)$. 
}
\end{remark}

\medskip

If $f$ is a $C^\infty$ real-valued function on a Jacobi manifold $M,$
the vector field $X_f$ defined by 
\begin{equation}\label{ch}
X_f=\#_{\Lambda}(df)+fE
\end{equation}
is called the {\it hamiltonian vector field} associated with $f$. It
should be noticed that the hamiltonian vector field associated with
the constant function $1$ is just $E$. A direct computation proves
that (see \cite{L2,M})
\begin{equation}\label{corh}
[X_f,X_g]=X_{\{f,g\}},
\end{equation}
which shows that the mapping 
\[
C^{\infty}(M,\R) \longrightarrow {\frak X}(M),\makebox[1cm]{}
f\mapsto X_f
\]
is a Lie algebra homomorphism.

Now, for every $x\in M,$ we consider the subspace
${\cal F}_{x}$ of $T_{x}M$
generated by all the hamiltonian vector
fields evaluated at the point $x$. In other words,
${\cal F}_{x} = (\#_{\Lambda})_{x}(T_{x}^*M) + \langle E_{x} \rangle$.
Since ${\cal F}$ is involutive, one easily follows that ${\cal F}$ defines a
generalized foliation in the sense of Sussmann \cite{Sus}, which is
called the {\it characteristic foliation} (see \cite{DLM,GL}).
Moreover, the Jacobi structure of $M$ induces  
a Jacobi structure on each leaf. In fact, if $L$ is the leaf over a
point $x$ of $M$ and $E_{x} \notin Im (\#_{\Lambda})_{x}$ (or
equivalently, the dimension of $L$ is odd) then $L$ is a contact
manifold with the induced Jacobi structure. If $E_{x} \in Im
(\#_{\Lambda})_{x}$ (or equivalently, the dimension of $L$ is even),
$L$ is a l.c.s. manifold (for a detailed study of the characteristic
foliation, we refer to \cite{DLM,GL}). If $M$ is a Poisson 
manifold then, from (\ref{e7}) and (\ref{ch}), we deduce that the
characteristic foliation of $M$ is just the {\it canonical symplectic
foliation } of $M$ (see \cite{Vai2,W}).

For a contact (respectively, l.c.s.) manifold $M$ there exists a
unique leaf of its characteristic foliation: the manifold $M$.

On the other hand, if ${\frak g}$ is a real Lie algebra of dimension
$n$ and $G$ is a connected Lie group with Lie algebra ${\frak g}$,
then the leaves of the symplectic foliation associated to the
Lie-Poisson structure on ${\frak g}^{*}$ are just the orbits of the
coadjoint action $Ad^{*}: G \times {\frak g}^{*} \longrightarrow {\frak
g}^{*}$ (see \cite{Vai2,W}).

Moreover, if $<\;,\;>$ is a scalar product on ${\frak g}$ then, under
the canonical identification between ${\frak g}$ and ${\frak g}^*$,
the coadjoint action induces an action of the Lie group $G$ on
${\frak g}$, which we will denote by $\widetilde{Ad^{*}}.$ Thus, we can
define an action of $G$ on the unit sphere $S^{n-1}({\frak g})$ as
follows: 
\begin{equation}\label{cobar}
\overline{Ad}^{*}: G\times S^{n-1}({\frak g})\longrightarrow
S^{n-1}({\frak g}), \makebox[1cm]{} (g,\xi)\mapsto
\overline{Ad}^{*}(g, \xi) =
\frac{\widetilde{Ad^{*}}(g, \xi)}{\|\widetilde{Ad^{*}}(g, \xi)\|}.
\end{equation}
Now, denote by $(\Lambda, E)$ the Jacobi structure on $S^{n-1}({\frak
g})$ defined in Section 2.2 and by $\xi_{S^{n-1}({\frak g})}$ the
infinitesimal generator, with respect to the action
$\overline{Ad}^{*}$, associated to $\xi \in {\frak g}$. Then, using the
results of \cite{L3}, we deduce that $\xi_{S^{n-1}({\frak g})}$ is
the hamiltonian vector field on $(S^{n-1}({\frak g}), \Lambda, E)$
associated to the function $<\xi, \; >: S^{n-1}({\frak g}) \longrightarrow
\R$ given by (\ref{Fusc}), that is,
\begin{equation}\label{Fund}
\xi_{S^{n-1}({\frak g})} = X_{<\xi, \; >}.
\end{equation}
This fact implies that the leaves of the characteristic foliation of
$S^{n-1}({\frak g})$ are just the orbits of the action
$\overline{Ad}^{*}$ (for more details, see \cite{L3}).

\subsection{Lie algebroid of a Jacobi manifold}

A {\it Lie algebroid structure }on a differentiable vector bundle
$\pi:K\longrightarrow M$ is a pair that consists of a Lie algebra
structure $\lcf \;,\; \rcf$ on the space $\Gamma(K)$ of the global cross
sections of $\pi:K \longrightarrow M$ and a homomorphism of
vector bundles $\varrho: K \longrightarrow TM$, the {\it anchor map},
such that if we also denote by $\varrho:\Gamma(K)\longrightarrow {\frak
X}(M)$ the homomorphism of $C^\infty(M,\R)$-modules
induced by the anchor map then:
\begin{enumerate}
\item
$\varrho:(\Gamma(K),\lcf\;,\;\rcf)\longrightarrow ({\frak X}(M),[\;,\;])$
is a Lie algebra homomorphism.
\item
For all $f\in C^\infty(M,\R)$ and for all $s_1,s_2\in \Gamma(K)$
one has
\[
\lcf s_1, fs_2 \rcf = f\lcf s_1, s_2\rcf + (\varrho(s_1)(f))s_2.
\]
\end{enumerate}
A triple $(K, \lcf\;,\;\rcf, \varrho)$ is called {\it a Lie algebroid
over $M$} (see \cite{Pra,Vai2}).

Let $(M,\Lambda,E)$ be a Jacobi manifold. In \cite{KS}, the authors
obtain a Lie algebroid structure on the vector bundle $J^1(M,\R)
\cong T^{*}M \times \R \longrightarrow M$ as follows.

Consider the homomorphism of $C^\infty(M, \R)$-modules
\[
(\#_{\Lambda},E):\Gamma(J^1(M,\R)) \cong \Omega^1(M)\times
C^\infty(M, \R)\rightarrow {\frak X}(M)
\]
defined by
\begin{equation}\label{46''}
(\#_{\Lambda},E)(\alpha,f)=\#_{\Lambda}(\alpha) + fE.
\end{equation}
It is clear that the vector field $(\#_{\Lambda},E)(\alpha,f)$ is tangent to
the characteristic foliation (note that $(\#_{\Lambda},E)(df,f)=X_f$).
Moreover, if $(\alpha,f),(\beta,g)\in \Omega^{1}(M) \times
C^{\infty}(M, \R)$ then (see \cite{KS})
\begin{equation}\label{4.6'}
[\#_{\Lambda}(\alpha)+fE,\#_{\Lambda}(\beta)+gE]=\#_{\Lambda}(\gamma)+hE,
\end{equation}
with $(\gamma,h)\in \Omega^1(M)\times C^\infty(M, \R)$ given by
\begin{equation}\label{4.7}
\begin{array}{lcl}
\gamma&=&{\cal L}_{\#_{\Lambda}(\alpha)}\beta-{\cal
L}_{\#_{\Lambda}(\beta)}\alpha-d(\Lambda(\alpha,\beta))+f{\cal L}_E\beta-g{\cal
L}_E\alpha -i_E(\alpha\wedge\beta),\\
h&=& \alpha(\#_{\Lambda}(\beta))+
\#_{\Lambda}(\alpha)(g)-\#_{\Lambda}(\beta)(f)+fE(g)-gE(f).
\end{array}
\end{equation}
This result suggests to introduce the mapping
$\lcf \; , \;\rcf_{(\Lambda, E)}: (\Omega^{1}(M) \times C^{\infty}(M,
\R))^{2} \longrightarrow \Omega^{1}(M) \times C^{\infty}(M, \R)$
defined by
\begin{equation}\label{4.7'}
\lcf (\alpha,f), (\beta,g)\rcf_{(\Lambda, E)} = (\gamma,h).
\end{equation}

This mapping gives a Lie algebra structure on $\Omega^{1}(M) \times
C^{\infty}(M, \R)$ in such a way that the triple
$(T^{*}M \times \R, \lcf \; , \;\rcf_{(\Lambda, E)},
(\#_{\Lambda},E))$ is a Lie algebroid over $M$ (see \cite{KS}).

\begin{remark}
{\rm
$i)$  If $\{\; , \; \}$ is the Jacobi bracket then the prolongation
mapping 
\begin{equation}\label{1j}
j^1:(C^\infty(M,\R), \{\; , \; \}) \longrightarrow (\Omega^1(M)\times
C^\infty(M,\R),\lcf \; , \;\rcf_{(\Lambda, E)}) \makebox[1cm]{}
f\mapsto j^1f=(df,f)
\end{equation}
is a Lie algebra homomorphism (see \cite{KS}).

$ii)$ In the particular case when $M$ is a Poisson manifold
we recover, by projection, the usual Lie algebroid structure on the
vector bundle $\pi:T^*M\longrightarrow M$ (see \cite{BV1,BV,CDW}).}
\end{remark}

\setcounter{equation}{0}

\section{Lichnerowicz-Jacobi cohomology of a Jacobi manifold}

\subsection{H-Chevalley-Eilenberg cohomology and Lichnerowicz-Jacobi
coho\-mo\-lo\-gy of a Jacobi manifold}

First of all, we recall the definition of the
cohomology of a Lie algebra ${\cal A}$ with
coefficients in an ${\cal A}$-module (we will follow
\cite{Vai2}).

Let $({\cal A},[\;,\;])$ be a real Lie algebra (not necessarily  finite
dimensional) and ${\cal M}$ a real vector space endowed with a $\R$-bilinear
multiplication
\[
{\cal A}\times {\cal M}\longrightarrow {\cal M}, \makebox[.4cm]{} (a,
m) \longrightarrow a.m
\]
such that
\begin{equation}\label{Comp}
[a_1,a_2].m=a_1.(a_2.m)-a_2.(a_1.m),
\end{equation}
for $a_{1}, a_{2} \in {\cal A}$ and $m \in {\cal M}$.
In other words, ${\cal A}$ acts on ${\cal M}$ on the left.
In such a case, a $k$-linear skew-symmetric  mapping
$c^k:{\cal A}^k\longrightarrow {\cal M}$ is called an {\it ${\cal
M}$-valued $k$-cochain }. These cochains form a real vector space  
$C^k({\cal A};{\cal M})$ and the linear operator $\partial^{k}:
C^k({\cal A};{\cal M}) \longrightarrow C^{k+1}({\cal A};{\cal M})$
given by
\begin{equation}\label{12''}
\begin{array}{rcl}
(\partial^k c^k)(a_0,\cdots,
a_k)&=&\displaystyle\sum_{i=0}^k(-1)^ia_i.c^k(a_0,\cdots,
\widehat{a_i},\cdots, a_k)+\\[5pt]
&&+ \displaystyle\sum_{i<
j}(-1)^{i+j}c^k([a_i,a_j],a_0,\cdots,
\widehat{a_i},\cdots,\widehat{a_j},\cdots ,a_k)
\end{array}
\end{equation}
defines a coboundary since $\partial^{k+1}\circ\partial^k=0$. Hence
we have the corresponding cohomology spaces
\[
H^k({\cal A};{\cal M}) = \frac{\ker \{\partial^k:C^k({\cal A};{\cal
M})\rightarrow C^{k+1}({\cal A};{\cal M})\}}{\mbox{\rm Im}
\{\partial^{k-1}:C^{k-1}({\cal A};{\cal
M})\rightarrow C^{k}({\cal A};{\cal M})\}}.
\]
This cohomology is called {\it the cohomology of the Lie algebra ${\cal
A}$ with coefficients in ${\cal M},$ or relative to the
given representation of ${\cal A}$ on ${\cal M}.$}

The {\it Chevalley-Eilenberg cohomology} of a Lie algebra $({\cal A},
[\;,\;])$ is just the cohomology of ${\cal A}$ relative to the
representation of ${\cal A}$ on itself given by
\[
a.m = [a, m].
\]

Now, let $(M,\Lambda,E)$ be a Jacobi manifold with Jacobi bracket
$\{\;,\;\}$. We consider the cohomology of the Lie algebra
$(C^\infty(M,\R),\{\;,\;\})$ relative to the representation defined
by the hamiltonian vector fields, that is,
\begin{equation}\label{RH}
C^\infty(M,\R)\times C^\infty(M,\R)\longrightarrow C^\infty(M,\R),
\makebox[1cm]{}(f,g)\longrightarrow X_f(g).
\end{equation}

This cohomology is denoted by $H^*_{HCE}(M)$ and it was called the
{\it H-Chevalley-Eilenberg cohomology } associated to $M$ (see
\cite{LME0,LME2,LME3,LME4}).
Explicitly, if $C_{HCE}^k(M)$ is the 
real vector space of the $k$-linear skew-symmetric mappings
$c^k:C^\infty(M,\R)\times \dots^{(k}\dots \times
C^\infty(M,\R)\longrightarrow C^\infty(M,\R)$ then
\[
H_{HCE}^k(M) = \frac{\ker \{\partial_H:C_{HCE}^k(
M)\rightarrow C_{HCE}^{k+1}(M)\}}{\mbox{\rm Im}
\{\partial_{H}:C_{HCE}^{k-1}(M)\rightarrow C_{HCE}^{k}(M)\}},
\]
where $\partial_H:C_{HCE}^r(M)\longrightarrow C_{HCE}^{r+1}(M)$ is the
linear differential operator defined by
\begin{equation}\label{ophce}
\begin{array}{rcl}
(\partial_H c^r)(f_0,\cdots,
f_r)&=&\displaystyle\sum_{i=0}^r(-1)^iX_{f_i}(c^r(f_0,\cdots,
\widehat{f_i},\cdots, f_r))+\\[5pt]
&&+ \displaystyle\sum_{i<
j}(-1)^{i+j}c^r(\{f_i,f_j\},f_0,\cdots,
\widehat{f_i},\cdots,\widehat{f_j},\cdots ,f_r)
\end{array}
\end{equation}
for $c^r\in C^r_{HCE}(M)$ and $f_0,\dots,f_r\in C^\infty(M,\R).$

\medskip

Note that for a Poisson manifold, $H_{HCE}^*(M)$ is the
{\it Chevalley-Eilenberg cohomology } of the Lie algebra
$(C^\infty(M,\R),\{\;,\;\})$ (see \cite{L1}). However, for arbitrary Jacobi
manifolds the Chevalley-Eilenberg cohomology (which is defined
with respect to the representation given by the Jacobi bracket
\cite{L2}) does not coincide in general with the  H-Chevalley-Eilenberg
cohomology defined above.

An interesting subcomplex of the H-Chevalley-Eilenberg complex is the
complex of the $1$-differentiable cochains.

A $k$-cochain $c^k\in C_{HCE}^k(M)$ is said to be {\it
$1$-differentiable }
if it is defined by a $k$-linear skew-symmetric differential operator
of order $1$. We can identify the space ${\cal
V}^k(M)\oplus {\cal V}^{k-1}(M)$ with the space of all
$1$-differentiable $k$-cochains $C_{HCE-1diff}^k(M)$ as follows
(see, for instance, \cite{L1}):
define $j^k:{\cal V}^k(M)\oplus {\cal V}^{k-1}(M)\longrightarrow
C_{HCE}^k(M)$ the monomorphism
given by
\begin{equation}\label{rp-hce}
j^k(P,Q)(f_1,\cdots ,f_k)=P(df_1,\dots,df_k)+
\displaystyle \sum_{q=1}^k (-1)^{q+1} f_q Q(df_1,\cdots
,\widehat{df_q},\cdots ,df_k).
\end{equation}
Then, $j^k({\cal V}^k(M)\oplus{\cal
V}^{k-1}(M))=C_{HCE1-diff}^k(M)$ which implies that the spaces ${\cal
V}^k(M)\oplus {\cal V}^{k-1}(M)$ and $C_{HCE1-diff}^k(M)$ are
isomorphic.

On the other hand, if $\tilde P\in C_{HCE1-diff}^k(M)$ then
$\partial_H \tilde P\in C_{HCE1-diff}^{k+1}(M).$ Thus, we
have the corresponding subcomplex
$(C_{HCE1-diff}^*(M),\partial_{H|C^*_{HCE1-diff}(M)})$ of the
H-Chevalley-Eilenberg complex whose cohomology $H^*_{HCE1-diff}(M)$
will be called the {\it $1$-differentiable H-Chevalley-Eilenberg
cohomology of $M$} (see \cite{LME3,LME4}). Moreover, using
(\ref{ophce}), (\ref{rp-hce}) and the properties of the
Schouten-Nijenhuis bracket, we can prove that
\begin{equation}\label{kj}
\partial_H(j^k(P,Q))=j^{k+1}(\sigma(P,Q)),
\end{equation}
where
\begin{equation}\label{opcoho}
\sigma(P,Q)=(-[\Lambda,P]+kE\wedge P+\Lambda\wedge
Q,[\Lambda, Q]-(k-1)E\wedge Q+[E,P]).
\end{equation}

The last equation defines a mapping
$\sigma:{\cal V}^k(M)\oplus {\cal V}^{k-1}(M)\longrightarrow {\cal
V}^{k+1}(M)\oplus {\cal V}^{k}(M)$
which is in fact a differential operator that verifies $\sigma^2=0.$
Thus, we have a complex $({\cal V}^*(M)\oplus {\cal
V}^{*-1}(M),\sigma)$ whose cohomology will be called
the {\it Lichnerowicz-Jacobi cohomology (LJ-cohomology) of $M$ } and
denoted by $H_{LJ}^{*}(M, \Lambda, E)$ or simply by $H_{LJ}^{*}(M)$
if there is not danger of confusion (see \cite{LME3,LME4}). This
cohomology is a generalization of the Lichnerowicz-Jacobi cohomology
introduced in \cite{LME0,LME1,LME2}.
In fact, the former one is the cohomology of the subcomplex
of the pairs $(P,0)$, where $P$ is invariant by $E$. For this
reason, we retain the name.

Notice that the mappings $j^k:{\cal
V}^k(M)\oplus {\cal V}^{k-1}(M)\longrightarrow C^k_{HCE}(M)$ given by
(\ref{rp-hce})
induce  an isomorphism between the complexes
\[
({\cal V}^*(M)\oplus {\cal V}^{*-1}(M),\sigma) \makebox[.4cm]{} \mbox{and}
\makebox[.4cm]{} (C_{HCE1-diff}^*(M),(\partial_H)_{|C_{HCE1-diff}^*(M)})
\]
and therefore the corresponding cohomologies are isomorphic.
\begin{remark}{\rm
If $\tilde\sigma$ denotes the cohomology operator in the
$1$-differentiable Chevalley-Ei\-len\-berg subcomplex then (see
\cite{L2})
\[
\tilde\sigma(P,Q)=(-[\Lambda,P]+(k-1)E\wedge P+\Lambda\wedge
Q,[\Lambda, Q]-(k-2)E\wedge Q+[E,P]),
\]
for $(P,Q)\in {\cal V}^k(M)\oplus {\cal V}^{k-1}(M).$ Thus, from
(\ref{opcoho}), we deduce that the $1$-differentiable
H-Chevalley-Eilenberg cohomology (that is, the LJ-cohomology) does not
coincide in general with the $1$-differentiable Chevalley-Eilenberg
cohomology. }
\end{remark}

To end this subsection, we will present another description of the
LJ-cohomology in terms of the Lie algebroid associated with the
Jacobi manifold (see \cite{LME3,LME4,Pva}).

Let $(K,\lcf\;, \;\rcf,\varrho)$ be a Lie algebroid over $M.$ For $r
\geq 0$, let $\Gamma(\Lambda^r K^*)$ be the space of the smooth
sections of $\Lambda^r K^*$, that is, $\Gamma(\Lambda^r K^*)$ is the
space of $C^{\infty}(M, \R)$-linear skew-symmetric mappings
$\xi^{r}: \Gamma(K) \times \dots^{(r} \dots \times \Gamma(K)
\longrightarrow C^{\infty}(M, \R).$
Define
\[
\tilde{\partial}^{r}: \Gamma(\Lambda^r K^*) \longrightarrow
\Gamma(\Lambda^{r+1} K^*)
\]
by
\begin{equation}\label{3.6'}
\begin{array}{rcl}
(\tilde{\partial}^r \xi^{r})(s_0,\cdots,
s_r)&=&\displaystyle\sum_{i=0}^r (-1)^i
\varrho(s_{i})(\xi^{r}(s_0,\cdots,
\widehat{s_i},\cdots, s_r))\\[5pt]
&&+ \displaystyle\sum_{i<
j}(-1)^{i+j}\xi^{r}(\lcf s_i,s_j\rcf ,s_0,\cdots,
\widehat{s_i},\cdots,\widehat{s_j},\cdots ,s_r).
\end{array}
\end{equation}
The operator $\tilde{\partial}^{r}$ satisfies $\tilde{\partial}^{r+1}
\circ \tilde{\partial}^{r} = 0$. The corresponding cohomology is
called the {\it Lie algebroid cohomology of $K$ with trivial coefficients}
(see \cite{Ma}).

\begin{remark}
{\rm Consider the representation of the Lie algebra
$(\Gamma(K),\lcf\;,\;\rcf)$ onto the space $C^\infty(M,\R)$ given by
$s.f=\varrho(s)(f),$ for all $s\in \Gamma(K)$ and $f\in
C^\infty(M,\R),$ and denote by
$(C^*(\Gamma(K);$ $C^\infty(M,\R)),$ $\partial)$ the corresponding
differential complex. Then, the Lie algebroid cohomology of $K$ with
trivial coefficients is the one of the subcomplex of
$(C^*(\Gamma(K);$ $C^\infty(M,\R)),$ $\partial)$ consisting of the cochains
which are $C^\infty(M,\R)$-linear. 
}
\end{remark}

Now, let $(M, \Lambda, E)$ be a Jacobi manifold of dimension $n$ and
$(J^{1}(M, \R), \lcf \; , \;\rcf_{(\Lambda, E)}, (\#_{\Lambda}, E))$
the Lie algebroid over $M$ (see Section 2.4). Then, the LJ-cohomology
of $M$ is just the Lie algebroid cohomology of $J^{1}(M, \R)$ with
trivial coefficients. In fact, the homomorphisms of $C^{\infty}(M,
\R)$-modules $\tilde{j^{k}}: {\cal V}^{k}(M) \oplus {\cal V}^{k-1}(M)
\longrightarrow \Gamma(\Lambda^k J^{1}(M, \R)^*)$ defined by
\[
\tilde{j^{k}}(P, Q)((\alpha_{1}, f_{1}), \cdots ,(\alpha_{k}, f_{k}))
= \displaystyle P(\alpha_{1}, \cdots , \alpha_{k}) + \sum_{q=1}^{k}
(-1)^{q+1} f_{q} Q(\alpha_{1}, \cdots , \widehat{\alpha_{q}}, \cdots , 
\alpha_{k})
\]
for $(P, Q) \in {\cal V}^{k}(M) \oplus {\cal V}^{k-1}(M)$ and
$(\alpha_{1}, f_{1}), \cdots ,(\alpha_{k}, f_{k}) \in \Omega^{1}(M)
\times C^{\infty}(M, \R)$, induce an isomorphism between the
complexes $({\cal V}^{*}(M) \oplus {\cal V}^{*-1}(M), \sigma)$ and
$\kern-2pt(\displaystyle \kern-2pt\bigoplus_{k=1}^{n+1}
\kern-1pt\Gamma (\Lambda^k \kern-2pt J^{1}(M, \R)^*),
\linebreak \tilde{\partial}^{*})$ (for more details, see
\cite{LME3,LME4,Pva}).

\subsection{Lichnerowicz-Jacobi cohomology and conformal changes of
Jacobi structures}

In this section, we will prove that the LJ-cohomology is invariant
under conformal changes.

First, we will recall some definitions and results related with the
theory of Lie algebroids which will be useful in the sequel (we will
follow \cite{Ma}).

Suppose that the pair $(\lcf \; , \rcf, \varrho)$ (respectively,
$(\lcf \; , \rcf', \varrho')$) is a Lie algebroid structure on the
vector bundle $\pi: K \longrightarrow M$ (respectively, $\pi': K'
\longrightarrow M$). An isomorphism between the Lie algebroids $(K,
\lcf \; , \rcf, \varrho)$ and $(K', \lcf \; , \rcf', \varrho')$ is an
isomorphism of vector bundles $\phi: K \longrightarrow K'$ such that
if we denote by $\phi_1: \Gamma(K) \longrightarrow \Gamma(K')$ the
isomorphism of $C^{\infty}(M, \R)$-modules induced by $\phi: K
\longrightarrow K'$ then:
\begin{enumerate}
\item
$\varrho' \circ \phi = \varrho.$
\item
For all $s_{1}, s_{2} \in \Gamma(K)$,
$\phi_1 \lcf s_{1}, s_{2} \rcf = \lcf \phi_1(s_{1}), \phi_1(s_{2}) \rcf',$
that is, $\phi_1: \Gamma(K) \longrightarrow \Gamma(K')$ is a Lie
algebra homomorphism.
\end{enumerate}
Assume that $\phi: K \longrightarrow K'$ is an isomorphism between
the Lie algebroids $(K, \lcf \; , \rcf, \varrho)$ and $(K', \lcf \; ,
\rcf', \varrho')$. Then, we can consider the isomorphism of
$C^{\infty}(M, \R)$-modules $\phi^{r}: \Gamma(\Lambda^r (K')^*)
\longrightarrow \Gamma(\Lambda^r K^*)$ given by
\begin{equation}\label{3.8'}
\phi^{r}(\xi^{r})(s_{1}, \dots , s_{r}) = \xi^{r}(\phi_1(s_{1}),
\dots , \phi_1(s_{r}))
\end{equation}
for $\xi^{r} \in \Gamma(\Lambda^r (K')^*)$ and $s_{1}, \dots , s_{r}
\in \Gamma(K)$. A direct computation, using (\ref{3.6'}), proves that
\begin{equation}\label{3.8''}
\tilde{\partial}^{k} \circ \phi^{k} = \phi^{k+1} \circ
(\tilde{\partial}')^{k},
\end{equation}
where $\tilde{\partial}^{*}$ (respectively,
$(\tilde{\partial}')^{*}$) is the cohomology operator induced by the
Lie algebroid structure $(\lcf \; , \; \rcf, \varrho)$ (respectively,
$(\lcf \; , \; \rcf', \varrho')$). Therefore, the Lie algebroids
cohomologies of $K$ and $K'$ with trivial coefficients are
isomorphic.

Now, let $(\Lambda, E)$ be a Jacobi structure on $M$. A {\it conformal
change } of $(\Lambda, E)$ is a new Jacobi structure $(\Lambda_{a},
E_{a})$ on $M$ defined by
\begin{equation}\label{3.6''}
\Lambda_{a} = a \Lambda, \makebox[.5cm]{} E_{a} = X_{a} =
\#_{\Lambda}(da) + a E,
\end{equation}
$a$ being a positive $C^{\infty}$ real-valued function on $M$ (see
\cite{DLM,GL}). We remark that $(\Lambda, E) =
((\Lambda_{a})_{\frac{1}{a}}, (E_{a})_{\frac{1}{a}})$. Moreover, we
have the following

\begin{theorem}\label{t3.2'}
Let $(M, \Lambda, E)$ be a Jacobi manifold and $(\Lambda_{a}, E_{a})$
a conformal change of the Jacobi structure $(\Lambda, E)$. Then,
\[
H^{k}_{LJ}(M, \Lambda, E) \cong H^{k}_{LJ}(M, \Lambda_{a}, E_{a}),
\]
for all $k$.
\end{theorem}
{\bf Proof:} We define the isomorphism of vector bundles $\phi:
T^{*}M \times \R \longrightarrow T^{*}M \times \R$ by
\begin{equation}\label{3.9'}
\phi(\alpha_{x}, \lambda) = (\displaystyle \frac{1}{a(x)} \alpha_{x}
+ \lambda d(\frac{1}{a})(x), \frac{\lambda}{a(x)})
\end{equation}
for $\alpha_{x} \in T_{x}^{*}M$ and $\lambda \in \R$. Note that the
isomorphism of $C^{\infty}(M, \R)$-modules $\phi_1: \Omega^{1}(M)
\times C^{\infty}(M, \R) \longrightarrow \Omega^{1}(M) \times
C^{\infty}(M, \R)$ induced by $\phi$ is given by
\begin{equation}\label{3.6'''}
\phi_1(\alpha, f) = (\displaystyle \frac{1}{a} \alpha + f
d(\frac{1}{a}), \frac{f}{a}) = (\frac{1}{a} \alpha - \frac{f}{a^{2}}
da, \frac{f}{a})
\end{equation}
for all $(\alpha, f) \in \Omega^{1}(M) \times C^{\infty}(M, \R)$.

A direct computation, using (\ref{46''}), (\ref{4.7}), (\ref{4.7'}),
(\ref{3.6''}) and (\ref{3.6'''}), proves that
\[
(\#_{\Lambda_{a}}, E_{a}) \circ \phi = (\#_{\Lambda}, E), \makebox[.3cm]{}
\phi_1 \lcf (\alpha, f), (\beta, g) \rcf_{(\Lambda, E)} = \lcf
\phi_1(\alpha, f), \phi_1(\beta, g) \rcf_{(\Lambda_{a}, E_{a})}
\]
for all $(\alpha, f), (\beta, g) \in \Omega^{1}(M) \times C^{\infty}(M,
\R)$. Thus, $\phi$ defines an isomorphism between the Lie algebroids
$(T^{*}M \times \R, \lcf \; , \;\rcf_{(\Lambda, E)}, (\#_{\Lambda},
E))$ and $(T^{*}M \times \R, \lcf \; , \;\rcf_{(\Lambda_{a}, E_{a})},
(\#_{\Lambda_{a}}, E_{a}))$ associated with the Jacobi structures
$(\Lambda, E)$ and $(\Lambda_{a}, E_{a})$, respectively. Therefore, from
the results in Section 3.1, it follows that
\[
H^{k}_{LJ}(M, \Lambda, E) \cong H^{k}_{LJ}(M, \Lambda_{a}, E_{a}),
\]
for all $k$.
\hfill$\Box$

Finally, using Theorem \ref{t3.2'}, we deduce the result announced at the
beginning of this section

\begin{corollary}\label{c3.2''}
The LJ-cohomology is invariant under conformal changes of the Jacobi
structure.
\end{corollary}

\subsection{Lichnerowicz-Jacobi cohomology of a Poisson manifold}

Let $(M,\Lambda)$ be a Poisson manifold and $\sigma$ the
LJ-cohomology operator. Using (\ref{opcoho}), we obtain that
\begin{equation}\label{opcohop}
\sigma(P,Q)=(-[\Lambda,P]+\Lambda\wedge Q,[\Lambda,Q]),
\end{equation}
for $(P,Q)\in {\cal V}^k(M)\oplus {\cal V}^{k-1}(M).$

Denote by $\bar\sigma$ the cohomology operator of the subcomplex of
the pairs $(P, 0)$. Under the canonical identification
${\cal V}^k(M)\oplus \{0\}\cong {\cal V}^k(M)$, we have that
\begin{equation}\label{oclp}
\bar \sigma(P)=-[\Lambda,P].
\end{equation}
The cohomology of the complex $({\cal V}^*(M),\bar\sigma)$ is called
the {\it Lichnerowicz-Poisson cohomology } (LP-cohomology) of $M$ and
denoted by $H^{*}_{LP}(M, \Lambda)$ or simply by $H^{*}_{LP}(M)$ if
there is not danger of confusion (see \cite{L1,Vai2}). Note
that $\bar\sigma(\Lambda) = 0$ and thus $\Lambda$ is a $2$-cocycle in
the LP-complex of $M$. Therefore, we can define the homomorphism
$L^{k} : H^{k}_{LP}(M) \longrightarrow H^{k+2}_{LP}(M)$ by
\[
L^{k} [P] = [P \wedge \Lambda]
\]
for all $[P] \in H^{k}_{LP}(M)$.

In \cite{L1} (see also \cite{LJMPA}), Lichnerowicz has exhibited the
relation between the LJ-cohomology (the $1$-differentiable
Chevalley-Eilenberg cohomology in his terminology) and the
LP-co\-ho\-mo\-lo\-gy of a Poisson manifold. In fact,
if $\hbox{dim} \, H^{k}_{LP}(M) < \infty$, for all $k$, we have that
\begin{equation}\label{nueva}
H^{k}_{LJ}(M) \cong \frac{H^{k}_{LP}(M)}{Im \, L^{k-2}} 
\oplus \ker \, L^{k-1} \; .
\end{equation}
Next, we will obtain
some consequences of the results of Lichnerowicz and of another authors
about the LJ-cohomology of symplectic and Lie-Poisson structures.

\subsubsection{Symplectic structures}

Let $(M,\Omega)$ be a symplectic manifold of dimension $2m$. Denote
by $\Lambda$ the Poisson $2$-vector and by $\#_{\Lambda}:
\Omega^{k}(M) \longrightarrow {\cal V}^{k}(M)$ the
homomorphism of $C^{\infty}(M, \R)$-modules given by (\ref{e7}) and
(\ref{e8}). Since, in this case, $\#_{\Lambda}$ is an isomorphism of
$C^{\infty}(M, \R)$-modules and
\begin{equation}\label{23'}
\#_{\Lambda}(d\alpha)=-\bar\sigma(\#_{\Lambda}(\alpha))
\end{equation}
for all $\alpha \in \Omega^k(M)$ (see \cite{L1,Vai2}), it
follows that $\#_{\Lambda}$ induces an isomorphism between the de Rham
cohomology of $M$, $H_{dR}^*(M)$, and the LP-cohomology. Thus,
$H_{dR}^k(M) \cong H_{LP}^k(M)$.

Under this identification and since $\#_{\Lambda}(\Omega) =
\Lambda$, the homomorphism $L^{k} : H^{k}_{LP}(M) \cong
H^{k}_{dR}(M) \longrightarrow H^{k+2}_{LP}(M) \cong H_{dR}^{k+2}(M)$
is given by
\begin{equation}\label{L}
L^k([\alpha])=[\alpha\wedge \Omega],\makebox[1cm]{}
\end{equation}
for all $[\alpha] \in H_{dR}^k(M)$ and $0 \leq k \leq 2m$.

From (\ref{nueva}), we have that
\begin{equation}\label{ljcohos}
H_{LJ}^k(M)\cong \frac{H_{dR}^k(M)}{{\rm Im}L^{k-2}}\oplus \ker L^{k-1}.
\end{equation}
Therefore, if $b_{r}(M)$ is the $r$-th Betti number of
$M$, we obtain
\begin{equation}\label{cotas}
\begin{array}{lcl}
\dim H_{LJ}^k(M)&\leq& b_k(M)+b_{k-1}(M),\\[8pt]
\dim H_{LJ}^k(M)&\geq&\max\{b_k(M)-b_{k-2}(M),\;
b_{k-1}(M)-b_{k+1}(M)\}.
\end{array}
\end{equation}
Next, we will discuss the behaviour of some examples of symplectic
manifolds with respect to the inequalities (\ref{cotas}).

\smallskip

{\it \bf Exact symplectic manifolds.-} If $\Omega$ is an exact $2$-form
then the homomorphisms $L^{k}$ are null and, from
(\ref{ljcohos}), it follows that $H_{LJ}^k(M)\cong H_{dR}^k(M)\oplus
H_{dR}^{k-1}(M)$ (see \cite{LJMPA}). Consequently,
\[
\dim H_{LJ}^k(M)=b_k(M)+b_{k-1}(M).
\]
In particular, the dimension of $H_{LJ}^k(M)$ is a topological
invariant of $M$, for all $k$.

\smallskip

{\it \bf Lefschetz symplectic manifolds.-} A symplectic manifold
$(M,\Omega)$ of dimension $2m$ is said to be a {\it Lefschetz symplectic
manifold} if it satisfies the strong Lefschetz theorem, that is, if
for every $k$, $0 \leq k \leq m$, the homomorphism
\[
\Delta^k=L^{k+2(m-k-1)}\circ \dots \circ L^{k+2}\circ
L^k:H_{dR}^k(M)\longrightarrow H_{dR}^{2m-k}(M),\makebox[.5cm]{}
[\alpha]\mapsto \Delta^k([\alpha])=[\alpha\wedge \Omega^{m-k}]
\]
is an isomorphism.

If $(M, \Omega)$ is a Lefschetz symplectic manifold then it is easy
to prove that:
\begin{enumerate}
\item
$L^k$ is a monomorphism, for $k \leq m-1.$
\item
$L^k$ is an epimorphism, for $k \geq m-1.$
\end{enumerate}
Thus, we deduce that
\[
\begin{array}{lcl}
b_k(M)-b_{k-2}(M)\leq 0,&\mbox{ for }& k\geq m+1\\[5pt]
b_{k-1}(M)-b_{k+1}(M)\leq 0,&\mbox{ for }& k\leq m.
\end{array}
\]
Moreover, using (\ref{ljcohos}), we have that
\[
\begin{array}{ll}
H_{LJ}^k(M)\cong \displaystyle\frac{H_{dR}^k(M)}{{\rm Im}L^{k-2}},&
\mbox{ for }k\leq m,\\
H_{LJ}^k(M)\cong \ker L^{k-1},&\mbox{ for }k\geq m+1,
\end{array}
\]
which implies that
\[
\begin{array}{ll}
\dim H_{LJ}^k(M)=b_k(M)-b_{k-2}(M), & \mbox{for }k\leq m,\\
\dim H_{LJ}^k(M)=b_{k-1}(M)-b_{k+1}(M), &\mbox{for }k\geq m+1.
\end{array}
\]
Therefore, the dimension of $H^{k}_{LJ}(M)$ is a topological
invariant of $M$, for all $k$.

\begin{remark}{\rm $i)$ A manifold $M$ endowed with a complex
structure $J$ is said to be {\it K\"ahler} if it admits a Riemannian
metric $g$ compatible with $J$ and such that the K\"ahler $2$-form
$\Omega$ given by
\[
\Omega(X,Y)=g(X,JY),
\]
is closed (see \cite{KN}). In such a case, $\Omega$ defines a
symplectic structure on $M$. Furthermore, if $M$ is compact then $(M,
\Omega)$ is a Lefschetz symplectic manifold (see \cite{We}).

$ii)$ There exist examples of compact Lefschetz symplectic manifolds
which do not admit K\"ahler structures (see \cite{Osa1,Osa2}).
}
\end{remark}

\smallskip

{\it \bf Compact symplectic nilmanifolds.-} Let $G$ be a simply connected
nilpotent Lie group of even dimension and let $\tilde\Omega$ be a
left-invariant symplectic $2$-form on $G$. Suppose that $\Gamma$ is a
discrete subgroup of $G$ such that the space of right cosets
$\Gamma \backslash G$ is a compact manifold. Then, the $2$-form
$\tilde\Omega$ induces a symplectic $2$-form $\Omega$ on $\Gamma
\backslash G$ and thus $\Gamma \backslash G$ is a compact symplectic
nilmanifold.

Now, denote by ${\frak g}$ the Lie algebra of $G$ and by $H^*({\frak
g})$ the cohomology of ${\frak g}$ relative to the trivial
representation of ${\frak g}$ on $\R$:
\[
{\frak g}\times \R\longrightarrow \R,\makebox[1cm]{}(a,t)\mapsto a.t=0.
\]
We define the homomorphism $(L_{\frak g})^k:H^k({\frak
g})\longrightarrow H^{k+2}({\frak g})$ by
\begin{equation}\label{Lg}
(L_{\frak g})^k[\alpha]=[\alpha\wedge \tilde\Omega_{\frak g}]
\end{equation}
for $[\alpha]\in H^k({\frak g}),$ where $\tilde\Omega_{\frak g}:{\frak
g}\times {\frak g}\longrightarrow \R$ is the symplectic $2$-form on
${\frak g}$ induced by $\tilde\Omega.$

Using Nomizu's Theorem \cite{N}, we have that the canonical
homomorphism $i^k:H^k({\frak g})\longrightarrow
H_{dR}^k(\Gamma\backslash G)$ is an isomorphism. Moreover, from 
(\ref{L}) and (\ref{Lg}), we deduce that
\[
i^{k+2}\circ (L_{\frak g})^k=L^k\circ i^k
\]
for all $k.$ Therefore (see (\ref{ljcohos}))
\begin{equation}\label{ljcohon}
H_{LJ}^k(\Gamma\backslash G)\cong \frac{H^k({\frak g})}{{\rm Im}
(L_{\frak g})^{k-2}}\oplus \ker(L_{\frak g})^{k-1}.
\end{equation}

\begin{remark}
{\rm $i)$ A compact Lefschetz symplectic nilmanifold is necessarily
a torus (see \cite{BG}).

$ii)$ If the Lie group $G$ is completely solvable then
(\ref{ljcohon}) also holds since, in such a case, the canonical
homomorphism $i^k:H^k({\frak g})\longrightarrow
H_{dR}^k(\Gamma\backslash G)$ is also an isomorphism, for all $k$
(see \cite{Ha}).}
\end{remark}

\begin{example}
{\rm Let $H$ be the {\it Heisenberg group} which consists of real
matrices of the form
\[\left(
\begin{array}{lcl}
1&x&z\\[5pt]
0&1&y\\[5pt]
0&0&1
\end{array}\right)
\]
with $x,y,z\in \R.$ $H$ is a simply connected nilpotent Lie group of
dimension $3$. Denote by $G$ the nilpotent Lie group of dimension $4$
defined by $G=H\times \R.$

If $t$ is the usual coordinate on $\R$, a basis for the
left-invariant $1$-forms on $G$ is given by
\[
\{\tilde\alpha=dx,\tilde\beta=dy,\tilde\eta=dz-xdy,\tilde\gamma=dt\}.
\]
We have that
\begin{equation}\label{3.16}
d\tilde\alpha=d\tilde\beta=0, \makebox[.5cm]{}
d\tilde{\eta}=-\tilde{\alpha}\wedge \tilde{\beta},
\makebox[.5cm]{}d\tilde\gamma=0.
\end{equation}
Thus,
\begin{equation}\label{fsh}
\tilde\Omega=\tilde\alpha\wedge \tilde\eta + \tilde\beta \wedge
\tilde\gamma 
\end{equation}
is a left-invariant symplectic $2$-form on $G$.

On the other hand, if $\bar\Gamma$ is the subgroup of $H$ consisting
of those matrices whose entries are integers then
$\Gamma=\bar\Gamma\times \Z$ is a discrete subgroup of $G$ and the
space of right cosets $\Gamma\backslash G=(\bar\Gamma\backslash
H)\times S^1$ is a compact nilmanifold. In fact, $\Gamma\backslash G$
is the {\it Kodaira-Thurston manifold } (see \cite{K,T}).

The $1$-forms $\tilde\alpha,\tilde\beta,\tilde\eta$ and
$\tilde\gamma$ all descend to $1$-forms $\alpha,\beta,\eta$ and
$\gamma$ on $\Gamma\backslash G$. Moreover, using (\ref{3.16})
and Nomizu's Theorem, it follows that
\begin{equation}\label{dR}
\begin{array}{lcl}
H^0({\frak g})\cong H_{dR}^0(\Gamma\backslash G)&=&<\{1\}>,\makebox[.5cm]{}
H^1({\frak g})\cong H_{dR}^1(\Gamma\backslash G)=<\{[\alpha], [\beta],
[\gamma]\}>,\\[8pt]
H^2({\frak g})\cong H_{dR}^2(\Gamma\backslash G)&=& <\{[\alpha\wedge
\eta], [\alpha\wedge \gamma],
                [\beta\wedge \eta],[\beta\wedge \gamma]\}>,\\[8pt]
H^3({\frak g})\cong H_{dR}^3(\Gamma\backslash G)&=& <\{[\alpha\wedge
\beta\wedge \eta],
[\alpha\wedge\eta\wedge \gamma], [\beta\wedge \eta\wedge \gamma]\}>,\\[8pt]
H^4({\frak g})\cong H_{dR}^4(\Gamma\backslash G)&=&<\{[\alpha\wedge
\beta\wedge \eta\wedge \gamma]\}>.
\end{array}
\end{equation}
Therefore,
\begin{equation}\label{3.18'}
b_0(\Gamma\backslash G)=b_4(\Gamma\backslash G)=1,\makebox[.5cm]{}
b_1(\Gamma\backslash G)=b_3(\Gamma\backslash G)=3,\makebox[.5cm]{}
b_2(\Gamma\backslash G)=4.
\end{equation}
Now, from (\ref{Lg}), (\ref{ljcohon}) and (\ref{dR}), we deduce that
\[
\begin{array}{ll}
\dim H^0_{LJ}(\Gamma\backslash G)=\dim H^5_{LJ}(\Gamma\backslash
G)=1,&\dim H^2_{LJ}(\Gamma\backslash G)=4\\
\dim H^1_{LJ}(\Gamma\backslash G)=\dim H^4_{LJ}(\Gamma\backslash
G)=3,& \dim H^3_{LJ}(\Gamma\backslash G)=5.
\end{array}
\]
Consequently (see (\ref{3.18'})), we have
\[
\max\{b_k(\Gamma\backslash G)-b_{k-2}(\Gamma\backslash G),
b_{k-1}(\Gamma\backslash G)-b_{k+1}(\Gamma\backslash G)\}< \dim
H_{LJ}^k(M)< b_k(\Gamma\backslash G)+b_{k-1}(\Gamma\backslash G)
\]
for $k=2,3.$
}
\end{example}

\subsubsection{Lie-Poisson structures}

Let $(M, \Lambda)$ be an {\it exact Poisson manifold}, that is, there
exists a vector field $X$ on $M$ such that
\[
\Lambda = \bar{\sigma} X = -{\cal L}_{X}\Lambda.
\]
In \cite{L1}, Lichnerowicz proved that, under this condition, we have 
\[
H^{k}_{LJ}(M) \cong H^{k}_{LP}(M) \oplus H^{k-1}_{LP}(M),
\]
for all $k$.

Now, suppose that ${\frak g}$ is a real Lie algebra of
dimension $n$ and consider the Lie-Poisson structure $\bar\Lambda$ on
${\frak g}^*$ (see Section 2.2). Using (\ref{LieP1}), it follows that
$({\frak g}^{*}, \bar\Lambda)$ is an exact Poisson manifold. Thus, 
\begin{equation}\label{colp}
H_{LJ}^k({\frak g}^*)\cong H_{LP}^k({\frak g}^*)\oplus
H_{LP}^{k-1}({\frak g}^*).
\end{equation}
On the other hand, if ${\frak g}$ is the Lie algebra of a compact Lie
group, in \cite{GW} the authors prove that
\begin{equation}\label{colp1}
H_{LP}^k({\frak g}^*)\cong H^k({\frak g})\otimes Inv,
\end{equation}
where $Inv$ is the algebra of all {\it Casimir functions} on ${\frak
g}^{*}$, that is,
\[
Inv = \{f \in C^{\infty}({\frak g}^{*}, \R) / X_{f}= 0\}.
\]
Therefore, from (\ref{colp}) and (\ref{colp1}), we conclude that for
the Lie algebra ${\frak g}$ of a compact Lie group
\begin{equation}\label{29'}
H_{LJ}^k({\frak g}^*)\cong (H^k({\frak g})\otimes Inv)\oplus
(H^{k-1}({\frak g})\otimes Inv).
\end{equation}

\subsubsection{A quadratic Poisson structure}

Let $\Lambda$ be the quadratic Poisson structure on $\R^{2}$ defined by
\[
\Lambda = xy \, \frac{\partial}{\partial x} \wedge
\frac{\partial}{\partial y} \; ,
\]
where $(x,y)$ stand for the usual coordinates on $\R^{2}$.

The restriction of $\Lambda$ to the open subset $\R^{2}-
(\{(x,0)/x\in \R\}\cup \{(0,y)/y\in \R\})$ is a symplectic structure
and the points $(x,0),(0,y)$ with $x,y\in \R,$ are  singular 
points. Moreover, $(\R^{2}, \Lambda)$ is not an exact Poisson manifold
and we have (see \cite{Naka})
\begin{equation}\label{3.29'}
H^{0}_{LP}(\R^{2}, \Lambda) \cong \R,\;
H^{1}_{LP}(\R^{2}, \Lambda) \cong \R^{2}, \; 
H^{2}_{LP}(\R^{2}, \Lambda) \cong \R^{2} \; .
\end{equation}
Using these facts and (\ref{nueva}) we deduce
\begin{equation}
\label{3.29''}
H^{0}_{LJ}(\R^{2}, \Lambda,0) \cong \R, \;
H^{1}_{LJ}(\R^{2}, \Lambda, 0) \cong \R^{2}, \;
H^{2}_{LJ}(\R^{2}, \Lambda, 0) \cong \R^{3},\;
H^{3}_{LJ}(\R^{2}, \Lambda, 0) \cong \R^{2}.
\end{equation}

\subsection{Lichnerowicz-Jacobi cohomology of a contact manifold}

In this section we will study the LJ-cohomology of a contact manifold. 

\medskip

First, we will obtain a general result for Jacobi manifolds which
relates the de Rham cohomology and the LJ-cohomology.

Let $(M,\Lambda,E)$ be a Jacobi manifold. Denote by
$\#_{\Lambda}:\Omega^k(M)\rightarrow {\cal V}^k(M)$ the homomorphism of
$C^\infty(M,\R)$-modules given by (\ref{e7}) and (\ref{e8}).
Then, we have (see \cite{LME1,LME2}):
\begin{equation}\label{s-dl}
{\cal L}_E(\#_{\Lambda}(\alpha))=\#_{\Lambda}({\cal
L}_E\alpha),\makebox[.6cm]{} 
-[\Lambda,\#_{\Lambda}(\alpha)]+kE\wedge
\#_{\Lambda}(\alpha)=-\#_{\Lambda}(d\alpha) +
\#_{\Lambda}(i_E\alpha)\wedge
\Lambda ,
\end{equation}
for all $\alpha\in \Omega^k(M)$. Using
(\ref{e1}), (\ref{opcoho}) and (\ref{s-dl}), we deduce the following

\begin{proposition}\label{5.1}
Let $(M,\Lambda,E)$ be a Jacobi manifold and
$\tilde{F}^k:\Omega^k(M)\oplus \Omega^{k-1}(M)\longrightarrow {\cal
V}^k(M)\oplus {\cal V}^{k-1}(M)$ the homomorphism of
$C^\infty(M,\R)$-modules defined by
\begin{equation}\label{st}
\tilde{F}^k(\alpha,\beta)=(\#_{\Lambda}(\alpha)+E\wedge
\#_{\Lambda}(\beta),-\#_{\Lambda}(i_E\alpha)+E\wedge \#_{\Lambda}(i_E\beta))
\end{equation}
for all $\alpha\in \Omega^k(M)$ and $\beta\in \Omega^{k-1}(M).$ Then,
the homomorphisms $\tilde{F}^k$ induce a homomorphism of complexes
\[
\tilde{F}:(\Omega^*(M),-d)\oplus (\Omega^{*-1}(M),d)\longrightarrow ({\cal
V}^*(M)\oplus {\cal V}^{*-1}(M),\sigma).
\]
Thus, if $H_{dR}^*(M)$ is the de Rham cohomology of $M$, we have the
corresponding homomorphism in cohomology
\[
\tilde{F} : H_{dR}^*(M)\oplus
H_{dR}^{*-1}(M) \longrightarrow H_{LJ}^{*}(M).
\]
\end{proposition}

Now, let $(M,\eta)$ be a contact manifold and $(\Lambda,E)$ its
associated Jacobi structure. Denote by $\flat:{\frak X}(M)
\rightarrow \Omega^1(M)$ the isomorphism of $C^\infty(M,\R)$-modules
given by $\flat(X)=i_X(d\eta)+\eta(X)\eta.$ The isomorphism
$\flat:{\frak X}(M) \rightarrow \Omega^1(M)$ can be extended to a
mapping, which we also denote by $\flat$, from the space ${\cal
V}^k(M)$ onto the space $\Omega^k(M)$ by putting:
\[
\flat(X_1\wedge \dots \wedge X_k)=\flat(X_1)\wedge \dots \wedge \flat(X_k)
\]
for all $X_1,\dots, X_k\in {\frak X}(M).$ This  extension is
also an isomorphism of $C^\infty(M,\R)$-modules. In fact, it follows that
\begin{equation}\label{sbc}
\#_{\Lambda}\alpha=(-1)^k\flat^{-1}(\alpha)+E\wedge
 \#_{\Lambda}(i_E\alpha),
\end{equation}
for $\alpha\in \Omega^k(M)$ (see \cite{LME2}). Moreover, we have
\begin{theorem}\label{5.2}
Let $(M,\eta)$ be a contact manifold of  dimension $2m+1.$
Then, the homomorphism
\[
\tilde{F}^k:H_{dR}^k(M)\oplus
H_{dR}^{k-1}(M)\rightarrow H_{LJ}^k(M)
\]
is an isomorphism for all $k.$ Thus, $H_{LJ}^k(M)\cong
H_{dR}^k(M)\oplus H^{k-1}_{dR}(M).$
\end{theorem}
{\bf Proof:}
Using (\ref{st}), (\ref{sbc}) and the fact that $i_E\circ\flat=\flat\circ
i_\eta,$ we deduce that the homomorphism of $C^\infty(M,\R)$-modules
$\tilde{G}^k:{\cal V}^k(M)\oplus {\cal V}^{k-1}(M)\longrightarrow
\Omega^k(M)\oplus \Omega^{k-1}(M)$ given by
\[
\tilde{G}^k(P,Q)=((-1)^k(\flat(P)+\eta\wedge \flat(Q)-\eta\wedge \flat
(i_\eta P)), (-1)^{k-1}(\flat(i_\eta P)-\eta\wedge \flat( i_\eta Q)))
\]
is just the inverse homomorphism of $\tilde{F}^k:\Omega^k(M)\oplus
\Omega^{k-1}(M)\rightarrow {\cal V}^k(M)\oplus {\cal V}^{k-1}(M).$
Consequently,
\[
\tilde{F}:(\Omega^*(M),-d)\oplus
(\Omega^{*-1}(M),d)\rightarrow ({\cal V}^*(M)\oplus {\cal
V}^{*-1}(M),\sigma)
\]
is an  isomorphism of complexes.
\hfill$\Box$

\begin{remark}{\rm In \cite{LJMPA}, Lichnerowicz proved that the
$1$-differentiable Chevalley-Eilenberg cohomology of a contact
manifold is trivial (compare this result with Theorem \ref{5.2}).}
\end{remark}

\subsection{Lichnerowicz-Jacobi cohomology of a locally conformal
symplectic manifold}

In this section, we will study the LJ-cohomology of a l.c.s.
manifold.

We will distinguish the two following cases:

$i)$ {\sc The particular case of a g.c.s. manifold:}
Let $(M, \Omega)$
be a g.c.s. manifold  with Lee $1$-form $\omega$ and let $(\Lambda,
E)$ be the associated Jacobi structure. Then, there exists a
$C^{\infty}$ real-valued function $f$ on $M$ such that $\omega = df$
and the $2$-form $\bar\Omega = e^{-f} \Omega$ is symplectic. Denote
by $\bar\Lambda$ the Poisson $2$-vector on $M$ associated with the
symplectic $2$-form $\bar\Omega$. Using (\ref{Posy}), (\ref{Lcs}) and
Remark \ref{r0}, we deduce that
\[
\Lambda = e^{-f}\bar\Lambda, \makebox[.5cm]{} E =
\#_{\bar\Lambda}(d(e^{-f})).
\]
Thus, the Jacobi structure $(\Lambda, E)$ is a conformal change of
the Poisson structure $\bar\Lambda$ (see (\ref{3.6''})). Therefore,
from (\ref{ljcohos}) and Theorem \ref{t3.2'}, we obtain
\begin{theorem}\label{t3.8'} Let $(M, \Omega)$ be a g.c.s. manifold
of finite type with Lee $1$-form $\omega = df$. Then,
\[
H_{LJ}^k(M)\cong \frac{H_{dR}^k(M)}{{\rm Im}\bar{L}^{k-2}}\oplus \ker
\bar{L}^{k-1},
\]
for all $k$, where $H^{*}_{dR}(M)$ is the de Rham cohomology of $M$
and $\bar{L}^{r}: H^{r}_{dR}(M) \longrightarrow H^{r+2}_{dR}(M)$ is
the homomorphism defined by
\[
\bar{L}^{r}[\alpha] = [e^{-f} \alpha \wedge \Omega],
\]
for all $[\alpha] \in H^{r}_{dR}(M)$.
\end{theorem}

\begin{remark}{\rm Relation (\ref{ljcohos}) follows
directly from Theorem \ref{t3.8'}.
}
\end{remark}

\begin{example}\label{3.8'''}
{\rm Let $(N,\eta)$ be a contact manifold of finite type.
Consider on the product manifold $M = N \times \R$ the $2$-form
$\Omega$ given by
\begin{equation}\label{Exgcs}
\Omega=(pr_1)^*(d\eta)-(pr_2)^*(dt)\wedge (pr_1)^*(\eta),
\end{equation}
where $t$ is the usual coordinate on $\R$ and $pr_i$ $(i=1,2)$ are the
canonical projections of $M$ onto the first and second factor,
respectively. Then, $(M, \Omega)$ is a g.c.s. manifold with
Lee $1$-form $\omega=(pr_2)^*(dt).$ Moreover, in this case, the
symplectic $2$-form $\bar\Omega = e^{-t} \Omega$ is exact which
implies that the homomorphism $\bar{L}^{r}$ is null, for all $r$.
Consequently, using Theorem \ref{t3.8'}, it follows that
\[
H_{LJ}^k(M) \cong H^k_{dR}(M) \oplus
H^{k-1}_{dR}(M) \cong H^k_{dR}(N) \oplus H_{dR}^{k-1}(N).
\]
}
\end{example}

$ii)$ {\sc The general case:} Now, we will study the LJ-cohomology of
an arbitrary l.c.s. manifold. First, we will obtain some results
about a certain cohomology, introduced by Guedira and Lichnerowicz
\cite{GL}, which is associated to an arbitrary differentiable
manifold endowed with a closed $1$-form.

Let $M$ be a differentiable manifoldd and $\omega$ a closed $1$-form
on $M$.

Define the differential operator $d_w$ by (see \cite{GL})
\begin{equation}\label{dw}
d_\omega=d+e(\omega),
\end{equation}
$d$ being the exterior differential and $e(\omega)$ the operator
given by
\begin{equation}\label{ew}
e(\omega)(\alpha)=\omega\wedge \alpha
\end{equation}
for $\alpha\in \Omega^*(M).$

Since $\omega$ is closed, it follows that $d_\omega^2=0$. This result
allows us to introduce the differential complex
\[
\cdots \longrightarrow \Omega^{k-1}(M) \stackrel{d_\omega}
\longrightarrow \Omega^k(M)
\stackrel{d_\omega}\longrightarrow \Omega^{k+1}(M)
\longrightarrow \cdots
\]
Denote by  $H_\omega^*(M)$ the cohomology of this complex.

\begin{proposition}\label{dw0}
Let $M$ be a differentiable manifold and $\omega$ a closed $1$-form
on $M$. Then:
\begin{enumerate}
\item[$i)$]
The differential complex $(\Omega^*(M),d_\omega)$ is elliptic.
Thus, if $M$ is compact the cohomology  groups $H_\omega^k(M)$ have
finite dimension.
\item[$ii)$]
If $\omega$ is exact and $f$ is a $C^{\infty}$ real-valued function
such that $\omega = df$ then the mapping
\[
H_{dR}^k(M)\longrightarrow H_\omega^k(M),\makebox[1cm]{}
[\alpha]\mapsto [e^{-f}\alpha],
\]
is an isomorphism. Therefore, $H_\omega^k(M)\cong H_{dR}^k(M).$
\end{enumerate}
\end{proposition}
{\bf Proof:} $i)$ It is easy to check that the differential
operators $d$ and $d_{\omega}$ have the same symbol which implies
that the complex $(\Omega^*(M),d_\omega)$ is elliptic.

$ii)$ A direct computation proves the result.
\hfill$\Box$

If the $1$-form $\omega$ is not exact then, in general,
\[
H^{*}_{\omega}(M) \ncong H^{*}_{dR}(M).
\]
In fact, we will show next that if $M$ is compact and $\omega$ is
non-null and parallel with respect to a Riemannian metric on $M$,
then the cohomology $H^{*}_{\omega}(M)$ is trivial. First, we will
recall some results proved by Guedira and Lichnerowicz \cite{GL}
which will be useful in the sequel.

Suppose that $M$ is a compact differentiable manifold of dimension
$n$, that $\omega$
is a closed $1$-form on $M$ and that $g$ is a Riemannian metric.

Consider the vector field $U$ on $M$ characterized by the condition
\begin{equation}\label{met}
\omega(X) = g(X,U),
\end{equation}
for all $X\in {\frak X}(M)$.

Denote by $\delta$ the codifferential operator given by
(see \cite{Gol})
\begin{equation}\label{del}
\delta\alpha=(-1)^{nk+n+1}(\star \circ d\circ \star)(\alpha),
\end{equation}
for all $\alpha\in \Omega^k(M),$ $\star$ being the Hodge star isomorphism.
Then, we define the operator $\delta_\omega:\Omega^{k}(M)\rightarrow
\Omega^{k-1}(M)$ by (see \cite{GL})
\begin{equation}\label{delw}
\delta_\omega=\delta+i_U,
\end{equation}
where $i_U$ denotes the contraction by the vector field $U$, that is
(see \cite{Gol}),
\begin{equation}\label{int}
i_U(\alpha)=(-1)^{nk+n}(\star\circ e(\omega)\circ \star)(\alpha),
\end{equation}
for $\alpha\in \Omega^k(M).$

Now, consider the standard scalar product $<\;,\;>$ on the space
$\Omega^k(M)$:
\begin{equation}\label{pe2}
<\;,\;> \; : \Omega^k(M)\times \Omega^k(M)\rightarrow \R,\makebox[1cm]{}
(\alpha,\beta)\mapsto <\alpha,\beta>=\int_M\alpha
\wedge\star\beta.
\end{equation}
Then, it is easy to prove that (see \cite{GL})
\begin{equation}\label{dual}
<d_\omega\alpha,
\beta>=<\alpha,\delta_\omega\beta>,
\end{equation}
for all $\alpha \in \Omega^{k-1}(M)$ and $\beta \in \Omega^{k}(M)$.

Thus, since $M$ is compact and the complex $(\Omega^*(M),d_\omega)$ is
elliptic, we obtain an orthogonal decomposition of the space
$\Omega^k(M)$ as follows
\begin{equation}\label{des}
\Omega^k(M)={\cal H}_\omega^k(M)\oplus
d_\omega(\Omega^{k-1}(M))\oplus \delta_\omega(\Omega^{k+1}(M)),
\end{equation}
where ${\cal H}_\omega^k(M)=\{\alpha\in \Omega^k(M)/
d_\omega(\alpha)=0, \;\; \delta_\omega(\alpha)=0\}$ (see \cite{GL}).

From (\ref{des}), it follows that
\begin{equation}\label{iso}
H_\omega^k(M)\cong {\cal H}_\omega^k(M).
\end{equation}

\medskip

Now, we will prove the announced result about the triviality of the
cohomology $H^{*}_{\omega}(M)$.

\begin{theorem}\label{dw1}
Let $M$ be a compact differentiable manifold and $\omega$ a closed
$1$-form on $M$, $\omega\not=0$. Suppose that $g$ is a Riemannian
metric on $M$ such that $\omega$ is parallel with respect to $g$.
Then, the cohomology $H_\omega^*(M)$ is trivial.
\end{theorem}

{\bf Proof:} Since $\omega$ is parallel and non-null it follows that $\|w\|=c$,
with $c$ constant, $c>0$. Assume, without the loss of generality, that $c
= 1$. Note that if $c\not=1$, we can consider the Riemannian metric
$g' = c^{2}g$ and it is clear that the module of $\omega$ with
respect to $g'$ is $1$ and that $\omega$ is also parallel with
respect to $g'$.

Under the hypothesis $c=1$, we have that
\begin{equation}\label{36}
\omega(U)=1.
\end{equation}
Using that $\omega$ is parallel and that $U$ is Killing, we obtain
that (see (\ref{del}), (\ref{int}) and \cite{Gol})
\begin{equation}\label{kill1}
{\cal L}_U = -\delta\circ e(\omega) - e(\omega)\circ \delta,
\end{equation}
\begin{equation}\label{kill}
\delta\circ {\cal L}_U={\cal L}_U\circ \delta.
\end{equation}
From (\ref{dw}), (\ref{ew}), (\ref{del}),
(\ref{delw}), (\ref{int}), (\ref{36}) and (\ref{kill}), we deduce the
following relations:
\begin{equation}\label{conm}
d_\omega\circ i_U=-i_U\circ d_\omega + {\cal L}_U+Id,
\makebox[1cm]{}\delta_\omega\circ i_U=-i_U\circ \delta_\omega,
\end{equation}
\begin{equation}\label{conm1}
d_\omega\circ {\cal L}_U={\cal L}_U\circ d_\omega,\makebox[1cm]{}
\delta_\omega\circ {\cal L}_U={\cal L}_U\circ \delta_\omega,
\end{equation}
where $Id$ denotes the identity transformation.

On the other hand, (\ref{kill1}) implies that
\[
\begin{array}{lcl}
<{\cal L}_U\alpha, \alpha>&=&-<\alpha, di_{U}\alpha + i_{U}d\alpha>\\
                          &=&-<\alpha, {\cal L}_U\alpha>
\end{array}
\]
for all $\alpha \in \Omega^{k}(M)$. Thus,
\begin{equation}\label{cond}
<{\cal L}_U\alpha,\alpha> = 0.
\end{equation}
Now, if $\alpha \in {\cal H}_\omega^k(M)$ then, using (\ref{conm}), we
have that
\[
{\cal L}_U\alpha=-\alpha+d_\omega(i_U\alpha).
\]
But, by (\ref{conm1}), we deduce that ${\cal L}_U\alpha\in {\cal
H}_\omega^k(M)$. Therefore (see (\ref{des})) we obtain that
\[
{\cal L}_{U}\alpha = -\alpha.
\]
Consequently, from (\ref{cond}), it follows that $\alpha = 0$.

This proves that ${\cal H}^{k}_{\omega}(M) = \{0\}$ which implies
that $H^{k}_{\omega}(M) = \{0\}$ (see (\ref{iso})).
\hfill$\Box$

Next, we will obtain some results which relate the LJ-cohomology of a
l.c.s. manifold $M$ with its de Rham cohomology and with the
cohomology $H_\omega^*(M),$ $\omega$ being the Lee $1$-form of $M$.

Let $(M,\Omega)$ be a l.c.s. manifold with Lee $1$-form $\omega$.
Suppose that $(\Lambda,E)$ is the associated Jacobi structure on $M$
and that $\flat:{\frak X}(M)\rightarrow \Omega^1(M)$ is the
isomorphism of $C^\infty(M,\R)$-modules defined by
$\flat(X)=i_X\Omega$. The isomorphism $\flat:{\frak X}(M)\rightarrow
\Omega^1(M)$ can be extended to a mapping, which we also denote by
$\flat$, from the space ${\cal V}^k(M)$ onto the space $\Omega^k(M)$
by putting:
\[
\flat(X_1\wedge \dots \wedge X_k)=\flat(X_1)\wedge\dots  \wedge \flat(X_k)
\]
for all $X_1,\dots , X_k\in {\frak X}(M).$ This extension is also an
isomorphism of $C^\infty(M,\R)$-modules. In fact, we have that
(see \cite{LME2})
\begin{equation}\label{sb}
\#_{\Lambda}\alpha=(-1)^k\flat^{-1}(\alpha)
\end{equation}
for all $\alpha\in \Omega^k(M),$ where
$\#_{\Lambda}:\Omega^k(M)\rightarrow {\cal V}^k(M)$ is the homomorphism
given by (\ref{e7}) and (\ref{e8}). 

Using (\ref{Lcs}), (\ref{e7}), (\ref{e8}), (\ref{sb}) and the fact that
$\#_{\Lambda}(\omega)=-E$, we obtain
\begin{equation}\label{sos}
\#_{\Lambda}\circ i_E=i_\omega\circ \#_{\Lambda},\makebox[1cm]{} i_E\circ
 \flat=-\flat
\circ i_\omega.
\end{equation}

Thus, from (\ref{s-dl}), (\ref{sb}) and (\ref{sos}), we deduce that
\begin{equation}\label{coho}
-\flat[\Lambda,P]+k\omega\wedge \flat(P)=d\flat(P)-i_E(\flat(P))\wedge
\Omega,\makebox[1cm]{} {\cal L}_E\flat (P)=\flat({\cal L}_EP)
\end{equation}
for all $P\in {\cal V}^k(M).$

Furthermore, we prove the following

\begin{theorem}\label{6.3}
Let $(M,\Omega)$ be a l.c.s. manifold with Lee $1$-form $\omega$.
Suppose that $(\Lambda,E)$ is the associated Jacobi structure on $M$
and that $\bar{F}^k:\Omega^k(M)\longrightarrow {\cal
V}^k(M)\oplus {\cal V}^{k-1}(M)$ and $\bar{G}^k:{\cal V}^k(M)\oplus
{\cal V}^{k-1}(M) \rightarrow \Omega^{k-1}(M)$ are the homomorphisms
of $C^\infty(M,\R)$-modules defined by
\[
\bar{F}^k(\alpha)=(\#_{\Lambda}\alpha, -\#_{\Lambda}(i_E\alpha))
\makebox[2cm]{and}
\bar{G}^k(P,Q)=(-1)^k(-\flat(Q)+i_E\flat(P))
\]
for all $\alpha\in \Omega^k(M)$ and $(P,Q)\in {\cal V}^k(M)\oplus {\cal
V}^{k-1}(M)$. Then:
\begin{enumerate}
\item
The mappings $\bar{F}^k$ and ${\bar G^k}$ induce an exact sequence of
complexes
\[
0\longrightarrow (\Omega^{*}(M),d)\stackrel{\bar{F}}
\longrightarrow ({\cal V}^{*}(M)\oplus {\cal V}^{*-1}(M)
,-\sigma)\stackrel{\bar{G}}\longrightarrow (\Omega^{*-1}(M),-d_\omega)
\longrightarrow 0,
\]
where $d$ is the exterior differential, $\sigma$ is the LJ-cohomology
operator and $d_\omega$ is the operator given by (\ref{dw}).
\item
This exact sequence induces a long exact cohomology sequence
\[
\cdots \longrightarrow H^{k}_{dR}(M) \stackrel{\bar{F}^k_*}
\longrightarrow H^k_{LJ}(M)
\stackrel{\bar G^k_*}\longrightarrow H_{\omega}^{k-1}(M)
\stackrel{L^{k-1}}\longrightarrow H_{dR}^{k+1}(M)
\longrightarrow \cdots  ,
\]
with connecting homomorphism $L^{k-1}$ defined by
\[
L^{k-1}([\alpha])=[\alpha\wedge \Omega]
\]
for all $[\alpha]\in H_\omega^{k-1}(M)$.
\end{enumerate}
\end{theorem}

{\bf Proof:} It follows from (\ref{Lcs}), (\ref{opcoho}),
(\ref{s-dl}), (\ref{dw}), (\ref{sb}), (\ref{sos}) and (\ref{coho}).
\hfill$\Box$

Using Theorem \ref{6.3}, we obtain
\begin{corollary}\label{6.4}
Let $(M,\Omega)$ be a l.c.s. manifold of finite type with Lee $1$-form
$\omega$. Suppose that the dimension of the $k$-th cohomology group
$H_\omega^k(M)$ is finite, for all $k.$ Then,
\[
H_{LJ}^k(M)\cong \frac{H_{dR}^k(M)}{{\rm Im}L^{k-2}}\oplus \ker L^{k-1},
\]
where $L^r:H_{\omega}^r(M)\longrightarrow H_{dR}^{r+2}(M)$ is the
homomorphism given by
\[
L^r([\alpha])=[\alpha \wedge \Omega]
\]
for $[\alpha]\in H_{\omega}^r(M).$ In particular, the dimension of
$H_{LJ}^k(M)$ is finite.
\end{corollary}

\begin{remark}
{\rm Theorem \ref{t3.8'} follows directly from Proposition \ref{dw0}
and Corollary \ref{6.4}.}
\end{remark}

Using Theorem \ref{dw1} and Corollary \ref{6.4}, we deduce the
following results
\begin{corollary}\label{6.5}
Let $(M,\Omega)$ be a l.c.s. manifold of finite type with Lee $1$-form
$\omega$ and such that the dimension of the $k$-th cohomology group
$H_\omega^k(M)$ is finite, for all $k.$ If the $2$-form $\Omega$ is
$d_{(-\omega)}$-exact, that is, there exists a $1$-form $\eta$ on $M$
such that
$\Omega=d\eta-\omega\wedge
\eta,$
then,
\[
H_{LJ}^k(M)\cong H_{dR}^k(M)\oplus
H_\omega^{k-1}(M),\makebox[1cm]{} \mbox{ for all $k$.} 
\]
\end{corollary}

\begin{corollary}\label{6.6}
Let $(M,\Omega)$ be a compact l.c.s. manifold with Lee $1$-form
$\omega$, $\omega\not=0$. Suppose that $g$ is a Riemannian metric on
$M$ such that $\omega$ is parallel with respect to $g$. Then,
\[
H_{LJ}^k(M)\cong
H_{dR}^k(M),\makebox[1cm]{}\mbox{  for all $k$.}
\]
\end{corollary}

\begin{example}\label{comcon}
{\rm Let $(N,\eta)$ be a compact contact manifold and consider on the
product manifold $M = N \times S^1$ the $2$-form $\Omega$ defined by
\begin{equation}\label{Exlcs}
\Omega=(pr_1)^*(d\eta)-(pr_2)^*(\theta)\wedge (pr_1)^*(\eta),
\end{equation}
$\theta$ being the length element of $S^1$. Then,
$(M,\Omega)$ is a l.c.s. manifold with Lee $1$-form $\omega =
(pr_2)^*(\theta).$ Furthermore, if $h$ is a Riemannian metric on $N$,
the $1$-form $\omega$ is parallel with respect to the Riemannian
metric $g$ on $M$ given by
\[
g = (pr_{1})^{*}(h) + \omega \otimes \omega.
\]
Therefore, using Corollary \ref{6.6}, we deduce
\[
H_{LJ}^k(M)\cong H_{dR}^k(M)\cong H_{dR}^k(N)\oplus H_{dR}^{k-1}(N).
\]
}
\end{example}

\subsection{Lichnerowicz-Jacobi cohomology of the unit sphere of a
real Lie algebra}

If ${\frak g}$ is a real Lie algebra of finite dimension endowed
with a scalar product then the unit sphere of ${\frak g}$ admits a
Jacobi structure (see Section 2.2). In this section, we will describe
the LJ-cohomology of the sphere for the case when ${\frak g}$ is the
Lie algebra of a compact Lie group.

First, we will prove some results which will be useful in the sequel.

Let $({\frak g},[\;,\;])$ be a real Lie algebra of dimension $n$ and
$<\;,\;>$ a scalar product on ${\frak g}$. Denote by $S^{n-1}({\frak
g})$ the unit sphere in ${\frak g}$.

\begin{lemma}\label{6.1}
If $\xi\in{\frak g}$ and $\tilde{\xi}:S^{n-1}({\frak g})\times \R\rightarrow
\R$ is the real $C^{\infty}$-function given by
\begin{equation}\label{tilde}
\tilde{\xi}(\eta,t)=e^t<\xi,\eta>
\end{equation}
for all $(\eta,t)\in S^{n-1}({\frak g})\times \R,$ then
\begin{equation}\label{partial}
\frac{\partial}{\partial t}(\tilde{\xi})=\tilde{\xi}.
\end{equation}

Moreover, if $\{\xi_i\}_{i=1,\dots n}$ is a basis of ${\frak g}$ we
have that the set $\{d\tilde{\xi}_i\}_{i=1,\dots ,n}$ is a global
basis of the space of $1$-forms on $S^{n-1}({\frak g})\times \R.$
\end{lemma}

{\bf Proof:} (\ref{partial}) follows directly from (\ref{tilde}).

On the other hand, let $F:{\frak g}-\{0\}\rightarrow S^{n-1}({\frak
g})\times \R$ be the diffeomorphism defined by (\ref{2.17'}).
Then, we deduce that
\[
\tilde{\xi}\circ F=<\xi,\;\;>
\]
for all $\xi\in {\frak g},$ where $<\xi,\;\;>:{\frak
g}-\{0\}\rightarrow \R$ is the real function given by
\[
<\xi,\;\;>(\eta)=<\xi,\eta>,
\]
for all $\eta \in {\frak g}-\{0\}$. This proves the second assertion of
Lemma \ref{6.1}.
\hfill$\Box$

Using Lemma \ref{6.1},  we obtain
\begin{lemma}
Let $P$ (respectively, $Q$) be a  $k$-vector (respectively, a
$(k-1)$-vector) on $S^{n-1}({\frak g})$. Denote by
$\widetilde{(P,Q)}$ the $k$-vector on  $S^{n-1}({\frak g})\times
\R$ given by
\begin{equation}\label{tildev}
\widetilde{(P,Q)}=e^{-kt}(P+\frac{\partial}{\partial t}\wedge Q).
\end{equation}

If ${\cal L}$ is the Lie derivative operator on $S^{n-1}({\frak
g})\times \R$ then,
\begin{equation}\label{inv}
{\cal L}_{\frac{\partial }{\partial
t}}\widetilde{(P,Q)}=-k\widetilde{(P,Q)},\makebox[1cm]{} 
\frac{\partial }{\partial t}(\widetilde{(P,Q)}(d\tilde \xi_1,\dots ,
d\tilde\xi_k))=0 
\end{equation}
for all $\xi_1,\dots \xi_k\in {\frak  g},$ where $\tilde{\xi_i}$
$(i=1,\dots ,k)$ is the real $C^{\infty}$-function on $S^{n-1}({\frak
g})\times \R$ given by (\ref{tilde}).
\end{lemma}

Now, we will describe the LJ-cohomology of $S^{n-1}({\frak g})$ for
the case when ${\frak g}$ is the Lie algebra of a compact Lie group.

\begin{theorem}\label{3.18}
Let ${\frak g}$ be the Lie algebra of a compact Lie group $G$
of dimension $n.$ Suppose that $<\;,\;>$ is a scalar product on
${\frak g}$ and consider on the unit sphere $S^{n-1}({\frak g})$ the
induced Jacobi structure. Then
\[
H_{LJ}^k(S^{n-1}({\frak g}))\cong H^k({\frak g})\otimes Inv
\]
for all $k$, where $H^*({\frak g})$ is the cohomology of
${\frak g}$ relative to the trivial representation of ${\frak g}$
on $\R$ and $Inv$ is the subalgebra of $C^\infty(S^{n-1}({\frak
g}),\R)$ defined by
\[
Inv=\{\varphi\in C^\infty(S^{n-1}({\frak g}),\R)/X_f(\varphi)=0, \;\;
\forall f\in C^\infty(S^{n-1}({\frak g}),\R)\}.
\]
\end{theorem}

{\bf Proof:} Denote by $\overline{Ad}^{*}$ the action of $G$ on
$S^{n-1}({\frak g})$ given by (\ref{cobar}).
This action induces a representation of ${\frak g}$ on the vector
space $C^\infty(S^{n-1}({\frak g}),\R)$ given by
\[
(\xi,\varphi)\in {\frak g}\times C^\infty(S^{n-1}({\frak
g}),\R)\rightarrow \xi_{S^{n-1}({\frak g})}(\varphi)\in
C^\infty(S^{n-1}({\frak g}),\R),
\]
$\xi_{S^{n-1}({\frak g})}$ being the infinitesimal generator, with
respect to the action $\overline{Ad}^{*}$, associated to $\xi
\in {\frak g}$.

The above representation allows us to consider the differential
complex
\[
(C^*({\frak g};C^\infty(S^{n-1}({\frak
g}),\R)),\partial)
\]
and its cohomology $H^*({\frak g};C^\infty(S^{n-1}({\frak g}),\R))$
(see Section 3.1).

We will show that
\[
H_{LJ}^k(S^{n-1}({\frak g}))\cong H^k({\frak
g};C^\infty(S^{n-1}({\frak g}),\R)),
\]
for all $k.$

\medskip

Let $C^k_{HCE}(S^{n-1}({\frak g}))$ be the space of $k$-cochains in
the H-Chevalley-Eilenberg complex of $S^{n-1}({\frak g}).$ We define
the homomorphism
\[
\mu^k:C^k_{HCE}(S^{n-1}({\frak g}))\longrightarrow
C^k({\frak g};C^\infty(S^{n-1}({\frak g}),\R))\]
by
\begin{equation}\label{musk}
(\mu^k(c^k))(\xi_1,\dots ,\xi_k)=c^k(<\xi_1,\;>,\dots, <\xi_k,\;>)
\end{equation}
for all $c^k\in C^k_{HCE}(S^{n-1}({\frak g}))$ and
$\xi_1,\dots ,\xi_k\in {\frak g},$  where $<\xi_j,\;\;>$ $(j=1,\dots
,k)$ is the real $C^{\infty}$-function on $S^{n-1}({\frak g})$ given
by (\ref{Fusc}).

Now, consider the homomorphism of $C^\infty(S^{n-1}({\frak g}),\R)$-modules
\[
\Phi^k:{\cal V}^k(S^{n-1}({\frak g}))\oplus {\cal
V}^{k-1}(S^{n-1}({\frak g}))\rightarrow C^k({\frak
g};C^\infty(S^{n-1}({\frak g}),\R))
\]
defined by
\begin{equation}\label{Desco}
\Phi^k=\mu^k\circ j^k,
\end{equation}
$j^k:{\cal V}^k(S^{n-1}({\frak g}))\oplus {\cal
V}^{k-1}(S^{n-1}({\frak g}))\rightarrow
C^k_{HCE}(S^{n-1}({\frak g}))$
being the mapping given by (\ref{rp-hce}).

A direct computation shows that 
\begin{equation}\label{3.47'}
\begin{array}{lcl}
(\Phi^k(P,Q))(\xi_1,\dots ,\xi_k)(\xi)&=&P(d<\xi_1,\;\;>,\dots,
d<\xi_k,\;\;>)(\xi)\\[8pt]
&&\kern-130pt\displaystyle + \sum_{i=1}^{k}(-1)^{i+1}
<\xi_i,\xi>Q(d<\xi_1,\;\;>,\dots ,\widehat{d<\xi_i,\;\;>},\dots,
d<\xi_k,\;\;>)(\xi)\\[10pt]
&=&(\widetilde{(P,Q)}(d\tilde{\xi}_1,\dots,
d\tilde{\xi}_k))(\xi,0),
\end{array}
\end{equation}
for all $(P,Q)\in {\cal V}^k(S^{n-1}({\frak g}))\oplus {\cal
V}^{k-1}(S^{n-1}({\frak g}))$, $\xi_1,\dots, \xi_k\in {\frak g}$ and
$\xi\in S^{n-1}({\frak g}),$ where $\widetilde{(P,Q)}$ is the
$k$-vector on $S^{n-1}({\frak g})\times \R$ defined by (\ref{tildev})
and $\widetilde{\xi}_i$ $(i=1,\dots ,k)$ is the function on
$S^{n-1}({\frak g})\times \R$ given by (\ref{tilde}).

\medskip

Using (\ref{CompL}), (\ref{Fund}), (\ref{12''}), (\ref{ophce}) and
(\ref{musk}), we have that the mappings $\mu^{k}$ induce a
homomorphism between the complexes $(C^*_{HCE}(S^{n-1}({\frak
g})),\partial_H)$ and $(C^{*}({\frak g};
C^{\infty}(S^{n-1}({\frak g}), \R)), \partial).$
Thus, the mappings
$\Phi^{k}$ induce a homomorphism between the complexes
$({\cal V}^{*}(S^{n-1}({\frak g}))\oplus {\cal
V}^{*-1}(S^{n-1}({\frak g})),\sigma)$ and $(C^{*}({\frak g};
C^{\infty}(S^{n-1}({\frak g}), \R)), \partial)$
(see (\ref{kj}) and (\ref{Desco})).

\medskip

On the other hand, if $\Phi^k(P,Q)=0$ then, from (\ref{inv}), 
(\ref{3.47'}) and Lemma \ref{6.1}, it follows that
\[
0=\widetilde{(P,Q)}=e^{-kt}(P+\frac{\partial }{\partial t}\wedge Q).
\]
Therefore,  $P=0$ and $Q=0.$

Consequently, $\Phi^k$ is a monomorphism.

\medskip

Next, we will see that $\Phi^k$ is an epimorphism.

Let $c^k:{\frak g}\times \dots^{(k}\dots \times {\frak g}\rightarrow
C^\infty(S^{n-1}({\frak g}),\R)$ be a $C^\infty(S^{n-1}({\frak
g}),\R)$-valued $k$-cochain.

We define a $k$-vector $R$ on $S^{n-1}({\frak g})\times \R$
characterized by the condition
\begin{equation}\label{R}
R(d\tilde{\xi_1},\dots d\tilde\xi_k)(\xi,t) = e^{kt}(c^k (\xi_1,\dots,
\xi_k)(\xi)) 
\end{equation}
for all $\xi_1,\dots ,\xi_k \in {\frak g}$ and $(\xi,t)\in
S^{n-1}({\frak g}) \times \R$.

From Lemma \ref{6.1}, we deduce that $R$ is well-defined and, using
(\ref{partial}) and (\ref{R}), we have that
\[
{\cal L}_{\frac{\partial}{\partial t}}R=0.
\]
This implies that
\begin{equation}\label{R1}
R=P+\frac{\partial }{\partial t}\wedge Q,
\end{equation}
with $(P,Q)\in {\cal V}^k(S^{n-1}({\frak g}))\oplus {\cal
V}^{k-1}(S^{n-1}({\frak g}))$.

Moreover, from (\ref{tildev}), (\ref{3.47'}), (\ref{R}) and
(\ref{R1}), it follows that
\[
\Phi^k(P,Q)=c^k.
\]
Thus, $\Phi^{k}$ is an epimorphism.

\medskip

Using the above facts, we conclude that
\[
H^{k}_{LJ}(S^{n-1}({\frak g})) \cong H^{k}({\frak g};
C^{\infty}(S^{n-1}({\frak g}), \R)),
\]
for all $k$.

Now, if we apply a general result of Ginzburg and Weinstein
(see Theorem 3.5 of \cite{GW}; see also
\cite{E}), we obtain that
\[
H^{k}({\frak g}; C^{\infty}(S^{n-1}({\frak g}), \R)) \cong
H^{k}({\frak g})\otimes \overline{Inv}
\]
where $\overline{Inv}$ is the algebra of $G$-invariant functions on
$S^{n-1}({\frak g})$ with respect to the action $\overline{Ad}^{*}$.

Finally, from (\ref{Fund}) and since the characteristic foliation of
$S^{n-1}({\frak g})$ is generated by the set of hamiltonian vector fields
\[
\{X_{<\xi, \;\; >} / \xi \in {\frak g} \},
\]
we deduce that $\overline{Inv} = Inv$.
\hfill$\Box$

It is well-known that if ${\frak g}$ is the Lie algebra of a compact
semisimple Lie group then $H^{2}({\frak g}) = \{0\}$. Therefore,
using Theorem \ref{3.18}, we have

\begin{corollary}\label{semisim}
Let ${\frak g}$ be the Lie algebra of a compact semisimple Lie group $G$
of dimension $n.$ Suppose that $<\;,\;>$ is a scalar product on
${\frak g}$ and consider on the unit sphere $S^{n-1}({\frak g})$ the
induced Jacobi structure. Then
\[
H_{LJ}^2(S^{n-1}({\frak g})) = \{0\}.
\]
\end{corollary}

\subsection{Table I}

The following table summarizes the main results obtained in Sections
3.3, 3.4, 3.5 and 3.6 about the LJ-cohomology of the different types
of Jacobi manifolds.

$\mbox{}$

\hspace{-21pt}
{\footnotesize
\begin{tabular}{|c|c|c|}
\hline\hline
&&\\[-3pt]
TYPE & LJ-COHOMOLOGY & REMARKS\\[5pt]
\hline\hline
&&\\[-8pt]
$(M,\Omega)$ symplectic &$H_{LJ}^k(M)\cong
\displaystyle\frac{H_{dR}^k(M)}{\mbox{Im}L^{k-2}}\oplus \ker
L^{k-1}$ & $L^r:H_{dR}^r(M)\longrightarrow H_{dR}^{r+2}(M)$\\
of finite type&&$[\alpha]\mapsto [\alpha\wedge \Omega]$\\[5pt]
\hline
&&\\[-6pt]
$M$ exact symplectic & $H_{LJ}^k(M)\cong H_{dR}^k(M)\oplus
H_{dR}^{k-1}(M) $ & $\dim
H_{LJ}^k(M)=b_k(M)+b_{k-1}(M)$\\
of finite type&&\\[5pt]
\hline
&&\\[-6pt]
$M^{2m}$ Lefschetz symplectic & $H_{LJ}^k(M)\cong
\displaystyle\frac{H_{dR}^k(M)}{\mbox{Im}L^{k-2}},\;\; k\leq m$&$\dim
H_{LJ}^k(M)=b_k(M)-b_{k-2}(M)$\\
of finite   type &&$k\leq m$\\
&$H_{LJ}^k(M)\cong \ker L^{k-1},\;\;k\geq m+1$&$\dim
H_{LJ}^k(M)=b_{k-1}(M)-b_{k+1}(M)$\\
&&$k\geq m+1$\\[5pt]
\hline
&&\\[-8pt]
& &${\frak g}$ Lie algebra of  $G$\\
$M=\Gamma\backslash G$ compact symplectic &
$H_{LJ}^k(M)\cong
\displaystyle\frac{H^k({\frak g})}{\mbox{Im}(L_{\frak
g})^{k-2}}\oplus \ker (L_{\frak g})^{k-1}$
& $\tilde{\Omega}_{\frak g}:{\frak g}\times {\frak g}\rightarrow
{\Bbb R}$
induced symplectic  \\[-5pt]
&& form \\[2pt]
nilmanifold && $(L_{\frak g})^r:H^r({\frak g})\rightarrow
H^{r+2}({\frak g})$\\
&&$[\alpha]\mapsto [\alpha\wedge \tilde{\Omega}_{\frak g}]$\\[2pt]
\hline
\hline
&&\\[-8pt]
Dual of a real Lie algebra &&\\
${\frak g}$ of finite dimension & $H_{LJ}^k({\frak g}^*)\cong
H_{LP}^k({\frak g}^*)\oplus H_{LP}^{k-1}({\frak g}^*)$&\\[4pt]
\hline
&&\\[-6pt]
Dual of the Lie algebra ${\frak g}$ &$H_{LJ}^k({\frak g}^*)\cong
(H^k({\frak g})\otimes Inv)\oplus$&$Inv\equiv$ subalgebra of the \\
of a compact Lie group &\hspace{1.5cm}$(H^{k-1}({\frak g})\otimes
Inv)$ &Casimir functions of ${\frak g}^*$\\[3pt]
\hline
\hline&&\\[-6pt]
$M$ contact & $H_{LJ}^k(M)\cong H_{dR}^k(M)\oplus H^{k-1}_{dR}(M)$ &\\[5pt]
\hline\hline&&\\[-8pt]
$(M,\Omega)$ g.c.s. of finite type&$H_{LJ}^k(M)\cong
\displaystyle\frac{H_{dR}^k(M)}{\mbox{Im}\bar{L}^{k-2}}\oplus \ker
\bar{L}^{k-1}$
&$\bar{L}^r:H_{dR}^r(M)\rightarrow H_{dR}^{r+2}(M)$ \\
with Lee $1$-form  $\omega=df$& &$[\alpha]\mapsto
 [e^{-f}\alpha\wedge\Omega]$\\[3pt]
\hline
$(M,\Omega)$ l.c.s. of finite type&&\\
with  Lee $1$-form $\omega$ &$H_{LJ}^k(M)\cong
\displaystyle\frac{H_{dR}^k(M)}{\mbox{Im}L^{k-2}}\oplus \ker
L^{k-1}$&$L^r:H_\omega^r(M)\rightarrow H_{dR}^{r+2}(M)$\\
and  $\dim H^*_\omega(M)<\infty$ &&$[\alpha]\mapsto [\alpha\wedge
\Omega]$\\[5pt]
\hline
&&\\[-8pt]
$M$ compact l.c.s.   &&\\
with Lee $1$-form $\omega$ &$H_{LJ}^k(M)\cong
H_{dR}^k(M)$&$\dim H_{LJ}^k(M)=b_k(M)$\\
$\omega$ parallel with respect to  &&\\
a Riemannian metric&&\\[3pt]
\hline
\hline&&\\[-8pt]
Unit sphere $S^{n-1}({\frak g})$ of the&&$Inv\equiv$ subalgebra of
the constant \\
Lie algebra ${\frak g}$ of  a compact &$H_{LJ}^{k}(S^{n-1}({\frak
g}))\cong H^k({\frak g})\otimes Inv$& functions on the leaves
\\
 Lie group $(\dim{\frak g}=n)$&& of the characteristic foliation\\[3pt]
\hline\end{tabular}
}

\vspace{-16pt}
\begin{center}
{\it Table I: LJ-cohomology}
\end{center}

\section{Lichnerowicz-Jacobi homology of a Jacobi manifold}
\setcounter{equation}{0}

\subsection{H-Chevalley-Eilenberg homology and Lichnerowicz-Jacobi
homology of a Jacobi manifold}

In a similar way that for the cohomology, firstly we recall the
definition of the homology of a Lie algebra ${\cal A}$ with
coefficients in an ${\cal A}$-module (see, for instance, \cite{CE}).

Let $({\cal A},[\;,\;])$ be a real Lie algebra (not necessarily
finite dimensional) and  ${\cal M}$ a  real vector space endowed with
a $\R$-bilinear multiplication
\[
{\cal A}\times {\cal M}\longrightarrow {\cal M},
\makebox[1cm]{}(a,m) \mapsto  a.m
\]
compatible with the bracket $[\;,\;],$ i.e., such that (\ref{Comp}) holds.

An {\it ${\cal M}$-valued  $k$-chain } is an element of the vector space
$C_k({\cal A};{\cal M})={\cal M}\otimes
\Lambda^k{\cal A},$ where 
$\Lambda^*{\cal A}$ is the exterior algebra of ${\cal A}.$  We can
consider the linear operator $\delta_k:C_k({\cal A};{\cal
M})\longrightarrow C_{k-1}({\cal A};{\cal M})$ characterized by
\begin{equation}\label{H-Ch-E}
\begin{array}{lcl}
\delta_k(m\otimes (a_1\wedge \dots \wedge a_k))& = &
\displaystyle\sum_{1\leq i \leq k}(-1)^ia_i.m \otimes (a_1\wedge\dots
\wedge \widehat{a_i} \wedge \dots \wedge a_k) \; + \\
&&\kern-70pt\displaystyle\sum_{1\leq i<j\leq k}(-1)^{i+j}m\otimes ([a_i,a_j]
\wedge a_1\wedge \dots \wedge \widehat{a_i} \wedge \dots \wedge
\widehat{a_j}\wedge \dots \wedge a_k),
\end{array}
\end{equation}
which satisfies  $\delta_{k-1}\circ \delta_k=0,$ for all  $k.$ Then,
we have the  corresponding homology spaces
\[
H_k({\cal A};{\cal M}) = \frac{\ker \{\delta_k:C_k({\cal A};{\cal
M})\rightarrow C_{k-1}({\cal A};{\cal M})\}}{\mbox{\rm Im}
\{\delta_{k+1}:C_{k+1}({\cal A};{\cal
M})\rightarrow C_{k}({\cal A};{\cal M})\}}.
\]
This  homology is said to be the {\it homology of the Lie algebra
${\cal A}$ with coefficients in ${\cal M}$ or relative to the given
representation of ${\cal A}$ on  ${\cal
M}.$}

Now, let $(M,\Lambda,E)$ be a Jacobi manifold and  $\{\;,\;\}$ the
associated Jacobi bracket. We consider the homology of the Lie
algebra $(C^\infty(M,\R),\{\;,\;\})$ relative to the
representation defined by the hamiltonian vector fields as in (\ref{RH}).
This homology is called  the {\it   
H-Chevalley-Eilenberg homology associated to $M.$}

We denote by $C_k^{HCE}(M)$ the space of the $k$-chains in the
H-Chevalley-Eilenberg complex, by  $\delta_H$ the homology operator
and  by  $H_k^{HCE}(M)$ the $k$-th homology group.
Then, if $f\otimes(f_1\wedge \dots \wedge f_k)\in
C_k^{HCE}(M)=C^\infty(M,\R)\otimes (\Lambda^k(C^\infty(M,\R))),$
\begin{equation}\label{61}
\begin{array}{lcl}
\delta_H(f\otimes (f_1\wedge \dots \wedge f_k)) & = &
\displaystyle\sum_{1\leq i\leq k}(-1)^i
X_{f_i}(f)\otimes(f_1\wedge\dots \wedge \widehat{f_i} \wedge
\dots \wedge f_k) \; + \\
&&\kern-70pt\displaystyle\sum_{1\leq i<j\leq k}(-1)^{i+j}f\otimes
(\{f_i,f_j\}
\wedge f_1\wedge \dots \wedge \widehat{f_i}\wedge \dots \wedge
\widehat{f_j}\wedge \dots \wedge f_k).
\end{array}
\end{equation}

On the other hand, the skew-symmetric $k$-multilinear mapping
$\tilde\pi_k:C^\infty(M,\R)\times \dots^{(k}\dots\times
C^\infty(M,\R)\rightarrow \Omega^k(M)\otimes \Omega^{k-1}(M)$
defined by
\[
\tilde\pi_k(f_1,\dots ,f_k) = (df_1\wedge \dots \wedge df_k,
\sum_{i=1}^k(-1)^{i+k}f_idf_1 \wedge \dots \wedge\widehat{df_i} \wedge
\dots \wedge df_k)
\]
induces a linear mapping $\pi_k:C^{HCE}_k(M) \rightarrow
\Omega^k(M)\oplus\Omega^{k-1}(M)$ characterized by
\begin{equation}\label{62}
\pi_k(f\otimes(f_1 \wedge \dots \wedge f_k)) = (fdf_1\wedge \dots
\wedge df_k, \sum_{i=1}^k(-1)^{i+k} f f_i df_1 \wedge \dots
\wedge\widehat{df_i} \wedge \dots \wedge df_k)
\end{equation}
for all $f\otimes(f_1\wedge \dots \wedge f_k)\in C_k^{HCE}(M).$

A direct computation, using  (\ref{cf}), (\ref{ch}), (\ref{61}) and
(\ref{62}), shows that
\begin{equation}\label{62'}
\delta\circ \pi_k=\pi_{k-1}\circ \delta_H,
\end{equation}
where $\delta:\Omega^r(M)\oplus \Omega^{r-1}(M)\longrightarrow
\Omega^{r-1}(M)\oplus \Omega^{r-2}(M)$ is the operator given by
\begin{equation}\label{63}
\begin{array}{rcl}
\delta(\alpha,\beta) &=&
(i(\Lambda)d\alpha-di(\Lambda)\alpha+ri_E\alpha + (-1)^r{\cal L}_E\beta,\\
&&i(\Lambda)d\beta-di(\Lambda)\beta+(r-1)i_E\beta+(-1)^ri(\Lambda)\alpha),
\end{array}
\end{equation}
$i(\Lambda)$ being the contraction by $\Lambda$.

Since the mappings $\pi_k$ are locally surjective, from
(\ref{62'}), it follows that
\begin{equation}\label{63'}
\delta^2=0.
\end{equation}
This fact allows us to consider the differential complex
\[
\cdots \longrightarrow \Omega^{k+1}(M)\oplus\Omega^k(M) \stackrel{\delta}
\longrightarrow \Omega^k(M)\oplus \Omega^{k-1}(M)
\stackrel{\delta}\longrightarrow \Omega^{k-1}(M)\oplus \Omega^{k-2}(M)
\longrightarrow \cdots
\]
whose homology is called the {\it  Lichnerowicz-Jacobi homology}
(LJ-homology) of $M$ and denoted by $H_*^{LJ}(M,\Lambda,E)$ or simply
by $H^{LJ}_*(M)$ if there is not danger  of confusion.

\begin{remark}\label{cano}
{\rm Let  $\Omega_B^k(M)$ be the space of the  basic $k$-forms
with respect to $E,$ that is, $\alpha\in \Omega_B^k(M)$ if and only if
\[
i_E\alpha=0,\makebox[1cm]{}{\cal L}_E\alpha=0.
\]
Denote by  $\bar\delta$ the homology operator of
the subcomplex of the  LJ-complex which consists  of the pairs
$(0,\alpha)$,  $\alpha$ being a basic form with respect to $E.$ Under
the canonical identification $\{0\}\oplus \Omega_B^k(M)\cong
\Omega_B^k(M)$ one has that
\[
\bar\delta\alpha = i(\Lambda) d\alpha - di(\Lambda)\alpha,
\makebox[1cm]{} \mbox{ for all } \; \alpha\in \Omega_B^k(M).
\]
The homology of the complex $(\Omega_B^*(M),\bar\delta)$ was
studied in  \cite{ChLM1} and \cite{ChLM2} and it was called the
{\it canonical homology of the Jacobi manifold } $M.$ This
name is justified by the fact that if $M$ is a
Poisson manifold $(E=0),$ then the homology of the complex
$(\Omega_B^*(M),\bar\delta)$ is  just the canonical homology
introduced by Brylinski \cite{Br} (see also
\cite{Ko}). Note that $d\bar\delta+\bar\delta d=0$ and thus one can
consider a double complex and the two spectral sequences associated
with it. The degeneration of these spectral sequences at the first
term and other related aspects were discussed in \cite{Br,FIL1,FIL2,I,40'}
(for the case of a Poisson manifold) and in \cite{ChLM1,ChLM2} (for
the case of a Jacobi manifold).}
\end{remark}

\subsection{Modular class of a Jacobi manifold and duality between
the Lichnerowicz-Jacobi cohomology and the Lichnerowicz-Jacobi homology}

In this section, we will show that the LJ-homology of a Jacobi
manifold $(M,\Lambda,E)$ of dimension $n$ is just the homology of the
Lie algebroid
$(T^{*}M \times \R, \lcf \; , \;\rcf_{(\Lambda,
E)},(\#_{\Lambda},E))$  with respect to  a certain flat
$(T^*M\times \R)$-connection on $\Lambda^{n+1}(T^*M\times \R).$ This
last homology was introduced  by
Vaisman in \cite{Pva}. Moreover, in this paper, Vaisman also
introduced the definition  of the
modular class of an orientable  Jacobi manifold and he proved that if
such a class is zero then there is a duality
between the LJ-homology and the LJ-cohomology.
 
Firstly, we will recall several results
of \cite{Xu} (see also \cite{ELW,K95,MX}).

Let $(K,\lcf \; , \;\rcf, \varrho)$ be a Lie algebroid over $M$.

Denote by $\lcf \; , \;\rcf_{\cal
A}:\Gamma(\Lambda^{r_1}K)\times
\Gamma(\Lambda^{r_2}K)\rightarrow \Gamma(\Lambda^{r_1+r_2-1}K) $ the
bracket characterized by the relations
\[
\lcf f, g \rcf_{\cal A}=0,\makebox[1cm]{} \lcf X, f \rcf_{\cal
A}=\varrho(X)(f),\makebox[1cm]{} \lcf X, Y \rcf_{\cal A}=\lcf X, Y
\rcf,
\]
\[
\lcf U, V\wedge W \rcf_{\cal A}=\lcf U, V \rcf_{\cal A}\wedge W +
(-1)^{(r_1+1)s_1}V\wedge \lcf U, W \rcf_{\cal A},
\]
for all $f,g\in C^\infty(M,\R),$ $X,Y\in \Gamma(K)$, $U\in
\Gamma(\Lambda^{r_1}K),$ $V\in \Gamma(\Lambda^{s_1}K)$ and $W\in
\Gamma(\Lambda^{s_2}K).$ Then, if $n$ is  the rank of $K$, it follows
that $({\cal A}=\bigoplus_{0\leq r\leq
n}\Gamma(\Lambda^rK),\lcf \; , \;\rcf_{\cal A})$ is a {\it Gerstenhaber
algebra } (see \cite{ELW,K95,MX,Xu}). 
  
On the other hand, a {\it $K$-connection on a vector bundle
$L\rightarrow M$} is a $\R$-bilinear mapping
\[
\nabla:\Gamma(K)\times
\Gamma(L)\rightarrow \Gamma(L),\makebox[1cm]{} (X,s)\mapsto
\nabla_Xs
\]
such that 
\[
\nabla_{fX}s=f\nabla_X s,\makebox[1cm]{} \nabla_{X}fs=f\nabla_Xs +
\varrho(X)(f)s,\makebox[1cm]{} \mbox{ for all } f\in C^\infty(M,\R).
\]
The curvature $R$ of a $K$-connection $\nabla$ may be defined as for
the usual connections. $\nabla$ is said to be {\it flat } if $R$ vanishes.

Any $K$-connection on $\Lambda^n K\rightarrow M$  defines a
differential operator
$D:\Gamma(\Lambda^rK)\rightarrow \Gamma(\Lambda^{r-1}K)$ locally
given by
\begin{equation}\label{D}
D(i(\omega)\Phi)=
(-1)^{n-k+1}(i(\tilde{\partial}^{n-r}\omega)\Phi +
\sum_{i=1}^n\alpha^h\wedge i(w)\nabla_{X_h}\Phi),
\end{equation}
where $\Phi\in \Gamma(\Lambda^nK)$, $\omega\in
\Gamma(\Lambda^{n-r}K^*)$, $\{X_h\}$ is a local basis of $\Gamma(K)$ and
$\{\alpha^h\}$ is the dual basis of $\Gamma(K^*).$
The operator $D$ generates the Gerstenhaber
algebra $({\cal A},\lcf \;, \; \rcf_{\cal A})$, that is,
for all $U_1\in\Gamma(\Lambda^{r_1} K)$ and $U_2\in \Gamma(\Lambda^{r_2}K)$
\[
\lcf U_1 , U_2\rcf_{{\cal A}}=(-1)^{r_1}(D(U_1\wedge U_2)-DU_1\wedge
U_2-(-1)^{r_1}U_1\wedge DU_2).
\]

Moreover, the connection $\nabla$ can be recovered from the operator
$D$. More precisely, we have that
\begin{equation}\label{rca}
\nabla_{X}\Phi=-X\wedge D\Phi,
\end{equation}
for all $X\in \Gamma(K)$ and $\Phi\in \Gamma(\Lambda^nK).$

In fact, (\ref{D}) and (\ref{rca}) define a  one-to-one correspondence
between $K$-connections on $\Lambda^nK$ and linear operators $D$
generating the Gerstenhaber algebra $({\cal A},\lcf \;, \; \rcf_{\cal
A})$. Under this correspondence,
a flat $K$-connection $\nabla$ corresponds to a operator $D$ of
square zero. Thus, a flat $K$-connection $\nabla$
induces a homology operator. The corresponding  
homology $H_*(K,\nabla)$ is the {\it
homology of the Lie algebroid $K$ with respect to the flat
$K$-connection $\nabla$.}
For two flat $K$-connections $\nabla$ and $\nabla'$  such that
their generating operators, $D$ and $D',$ satisfy
$D-D'=i(\tilde{\partial}f),$ with $f\in C^\infty(M,\R)$, one has
that $H_r(K,\nabla)\cong H_{r}(K,\nabla'),$ for all $r.$
Furthermore, if $\nu\in \Gamma(\Lambda^nK)$ is such that
$\nu(x)\not=0,$ for all $x\in M,$ and 
$\nabla\nu =0$
then it is possible to define a duality between the homology
$H_*(K,\nabla)$ and the
cohomology  of the Lie algebroid $K$ with trivial coefficients
$H^{*}(K).$  More precisely, the mapping
$\star:\Gamma(\Lambda^rK^*)\rightarrow \Gamma(\Lambda^{n-r}K)$ given
by
\[
\star \xi=i(\xi)\nu,
\]
induces an isomorphism between the cohomology group $H^r(K)$ and the
homology group $H_{n-r}(K,\nabla)$ (for more details, see \cite{Xu}).

Now, let $(M,\Lambda,E)$ be a Jacobi manifold of dimension $n$ and
$(T^*M\times \R,\lcf \;, \; \rcf_{(\Lambda,E)},$ $(\#_\Lambda,E))$ its
associated Lie algebroid (see Section 2.4).

The space $\Gamma(\Lambda^r(T^*M\times \R))$ can be identified with
$\Omega^r(M)\oplus \Omega^{r-1}(M)$ in such a way that the exterior
product of a section $(\alpha,\beta)$ of $\Lambda^r(T^*M\times
\R)\rightarrow M$ with a section $(\alpha',\beta')$ of
$\Lambda^{r'}(T^*M\times \R)\rightarrow M$ is given by
\begin{equation}\label{4.8'}
(\alpha,\beta)\wedge (\alpha',\beta')=(\alpha\wedge
\alpha',\alpha\wedge \beta'+ (-1)^{r'}\beta\wedge \alpha').
\end{equation}
On the other hand, under the identification of
$\Gamma(\Lambda^k(T^*M\times \R)^*)$ with ${\cal V}^k(M)\otimes {\cal
V}^{k-1}(M)$ (see Section 3.1)  the interior product of a section
$(\alpha,\beta)$ of $\Lambda^r(T^*M\times \R)\rightarrow M$ by a
section $(P,Q)$ of $\Lambda^k(T^*M\times \R)^*\rightarrow M$ is
given by
\[
\begin{array}{lcl}
\iota(P,Q)(\alpha,\beta)=(i(P)\alpha + (-1)^{r-1}i(Q)\beta,
i(P)\beta)&& \mbox{ if } k\leq r\\
\iota(P,Q)(\alpha,\beta)=0 &&\mbox{ if } k>r.
\end{array}
\]
In particular,
\begin{equation}\label{4.8''}
\iota(X,f)(\alpha,\beta)=(i(X)\alpha + (-1)^{r-1}f\beta, i(X)\beta),
\end{equation}
for all $(X,f)\in \Gamma(\Lambda^1(T^*M\times \R)^*)\cong {\frak
X}(M)\times C^\infty(M,\R).$

The Jacobi structure $(\Lambda,E)$ allows us to introduce a flat
$(T^*M\times \R)$-connection on $\Lambda^{n+1}(T^*M$ $\times \R)$
defined by (see \cite{Pva})
\begin{equation}\label{conexion}
\nabla_{(\alpha,f)}(0,\Phi)=
(0,fdi_E\Phi+\alpha\wedge (di(\Lambda)\Phi-ni_E\Phi)),
\end{equation}
for all $(\alpha,f)\in \Omega^1(M)\times C^\infty(M,\R)$ and $\Phi\in
\Omega^n(M).$ 
Then, if $\delta:\Omega^r(M)\oplus
\Omega^{r-1}(M)\rightarrow\Omega^{r-1}(M)\oplus \Omega^{r-2}(M)$ is
the LJ-homology operator (see (\ref{63})) and $D$  is the homology operator
associated with $\nabla,$ we have that $D=\delta$ (see relation
(2.10) in \cite{Pva}).
Therefore,

\begin{proposition}\label{a1}\cite{Pva}
Let $(M,\Lambda,E)$ be a Jacobi manifold of dimension $n$. Then the
LJ-homology $H^{LJ}_*(M)$ is the homology $H_*(T^*M\times \R,\nabla)$
of the Lie algebroid $(T^*M\times \R,\lcf \;, \; \rcf_{(\Lambda,E)},$ 
$(\#_\Lambda,E))$  with respect to the flat $(T^*M\times \R)$-connection
on $\Lambda^{n+1}(T^*M\times \R)$ defined by (\ref{conexion}).
\end{proposition}

Next, assume that $M$ is orientable and let $\nu$ be a volume form.
The volume form $\nu$ induces a flat $(T^*M\times
\R)$-connection $\nabla_0$ on $\Lambda^{n+1}(T^*M\times \R)$ by putting
\[
(\nabla_0)_{(\alpha,f)}(0,\nu)=(0,0),\makebox[1cm]{} \mbox{for
all }(\alpha,f) \in \Omega^1(M)\times C^\infty(M,\R).
\]
Then, using (\ref{63}), (\ref{D}), (\ref{rca}) and (\ref{4.8'}), we
have that for all $(\alpha,f)\in
\Omega^1(M)\times C^\infty(M,\R),$
\begin{equation}\label{modular1}
\begin{array}{rcl}
\nabla_{(\alpha,f)}(0,\nu)-(\nabla_0)_{(\alpha,f)}(0,\nu)&=&
-(\alpha,f)\wedge \delta(0,\nu)\\
&=&(0,(fdiv_\nu E-n\alpha(E))\nu+\alpha\wedge di(\Lambda)\nu),
\end{array}
\end{equation}
where $div_\nu E$ is the
divergence  of the vector field $E$ with respect to
$\nu$, that  is,
\[
{\cal L}_E\nu=(div_\nu E) \nu.
\]
Now, let ${\cal X}^\nu_{(\Lambda,E)}$ be the vector field characterized
by the relation
\begin{equation}\label{4.10'}
{\cal L}_{\#_{\Lambda}(df)}\nu={\cal
X}_{(\Lambda,E)}^\nu(f)\nu,\makebox[1cm]{} \mbox{ for all }f \in
 C^\infty(M,\R).
\end{equation}
Using (\ref{4.10'}) and the fact that
\begin{equation}\label{4.10''}
i_{\#_\Lambda(\alpha)}\nu=-\alpha\wedge i(\Lambda)\nu,
\end{equation}
for all $\alpha\in \Omega^1(M),$ it follows that
\begin{equation}\label{modular2}
\begin{array}{rcl}
\alpha({\cal X}_{(\Lambda,E)}^\nu)\nu&=&{\cal L}_{\#_\Lambda(\alpha)}\nu+
d\alpha\wedge i(\Lambda)\nu
\\&=&\alpha \wedge d i(\Lambda)\nu.
\end{array}
\end{equation}

Therefore,
\begin{equation}\label{n-n0}
\nabla_{(\alpha,f)}(0,\nu)-(\nabla_0)_{(\alpha,f)}(0,\nu)=
(0, (fdiv_\nu E+\alpha({\cal X}_{(\Lambda,E)}^\nu-nE))\nu).
\end{equation}
Denote by $D_0$ the corresponding homology operator associated with
$\nabla_0$. From (\ref{D}) and (\ref{n-n0}), we deduce that
\begin{equation}\label{d-d0}
D-D_0=\iota({\cal X}_{(\Lambda,E)}^\nu-nE,div_\nu E).
\end{equation}
The pair
\begin{equation}\label{4.13'}
{\cal M}_{(\Lambda,E)}^{\nu}=({\cal
X}_{(\Lambda,E)}^\nu-nE,div_\nu E)\in {\frak X}(M)\times C^\infty(M,\R)
\end{equation}
defines a $1$-cocycle in the LJ-complex of $M$, that is, $\sigma({\cal
M}_{(\Lambda,E)}^{\nu})=(0,0),$ where $\sigma$ is the LJ-cohomology
operator. Moreover, the corresponding cohomology
class ${\cal M}_{(\Lambda,E)}\in H_{LJ}^1(M)$ does not depend of the
volume form $\nu$ (see \cite{Pva}).

This cohomology class ${\cal M}_{(\Lambda,E)}$ is called the {\it
Jacobi modular class of $M$} (see \cite{Pva}). The manifold $M$ is said
to be {\it a unimodular Jacobi manifold} if the  Jacobi modular class
${\cal M}_{(\Lambda,E)}$ is zero. In such a case,
using (\ref{d-d0}), Proposition \ref{a1} and  the results of \cite{Xu}
described above, we conclude the following

\begin{theorem}\label{a2}
\cite{Pva} If $(M,\Lambda,E)$
is a unimodular Jacobi manifold of dimension $n$ then
\[
H_r^{LJ}(M)\cong H^{n+1-r}_{LJ}(M)
\]
for all $r\in \{0,\dots, n+1\}.$
\end{theorem}

\begin{remark}\label{r4.4}
{\rm Let $(M,\Lambda)$ be an orientable  Poisson manifold.
\medskip

$(i)$ Suppose that
$\nu$ is a volume form on $M.$ The {\it modular vector
field } of $M$ with respect to $\nu$ is the vector field ${\cal
X}_\Lambda^\nu$ characterized by
\begin{equation}\label{4.13''}
{\cal L}_{\#_\Lambda(df)}\nu={\cal X}_\Lambda^\nu(f)\nu,
\makebox[1cm]{} \mbox{ for all } f\in C^\infty(M,\R).
\end{equation}
If $\bar\sigma$ denotes the LP-cohomology operator (see (\ref{oclp}))
we have that $\bar{\sigma}({\cal X }_\Lambda^\nu)=0.$ Thus, ${\cal
X}_\Lambda^\nu$ defines a cohomology class ${\cal M}_\Lambda\in
H_{LP}^1(M).$
This class does not depend of the volume form $\nu$ and it is called
the {\it Poisson modular class } of $M.$ If ${\cal M}_\Lambda$ is
zero then $M$ is said to be a {\it unimodular Poisson manifold} (for
more details, we refer to \cite{W2,W3}).

\medskip
$(ii)$ Let $(Id^k,0):H_{LP}^k(M)\rightarrow H_{LJ}^k(M)$ be the
canonical homomorphism given by
\[
(Id^k,0)([P])=[(P,0)],\makebox[1cm]{}\mbox{for }[P]\in H_{LP}^k(M).
\]
Since the $0$-cochains in the LP-complex and in the LJ-complex are
the $C^\infty$ real-valued functions on $M,$ we deduce that
$(Id^1,0):H_{LP}^1(M)\rightarrow H^1_{LJ}(M)$ is a monomorphism.

On the other hand, from (\ref{4.10'}), (\ref{4.13'}) and
(\ref{4.13''}), it follows that
\[
(Id^1,0){\cal M}_\Lambda={\cal M}_{(\Lambda,0)},
\]
where ${\cal M}_{(\Lambda,0)}$ is the Jacobi modular class of $M.$
Therefore,  we conclude that
$(M,\Lambda)$ is a unimodular Poisson manifold if and only if
$(M,\Lambda,0)$
is a unimodular Jacobi manifold.}
\end{remark}
 
\begin{remark}\label{r4.5}{\rm Let $(M,\Lambda,E)$ be an orientable
Jacobi manifold and $(M\times \R,\tilde{\Lambda})$ the poissonization
of $M$. 

\medskip

$(i)$ Suppose that $\nu$ is a volume form on $M$ and consider in
$M\times \R$ the volume form
\[
\tilde{\nu}=e^{(n+1)t} \nu\wedge dt,
\]
where $n$ is the dimension of $M$ and $t$ is the usual  coordinate on
$\R.$ Using the results of Vaisman (see relations $(3.13)$ and $(3.14)$
in \cite{Pva}) we deduce that the
modular vector field ${\cal X}_{\tilde{\Lambda}}^{\tilde{\nu}}$ of
$(M\times \R,\tilde{\Lambda})$ with
respect to $\tilde{\nu}$ is
\[
{\cal X}_{\tilde{\Lambda}}^{\tilde\nu}=e^{-t}({\cal
X}_{(\Lambda,E)}^{\nu}-nE+(div_\nu E)\frac{\partial}{\partial t}).
\]
Thus, from (\ref{4.13'}), we conclude that ${\cal
X}_{\tilde{\Lambda}}^{\tilde{\nu}}$ is zero if and only if ${\cal
M}_{(\Lambda,E)}^\nu$ is zero.

\medskip

$(ii)$ Using again the results of Vaisman \cite{Pva}, we have that if
$(M,\Lambda,E)$ is a unimodular Jacobi manifold then $(M\times
\R,\tilde{\Lambda})$ is a unimodular Poisson manifold. However, in
general, the converse does not hold. In fact, the poissonization of a
contact manifold $M$ is unimodular (see Remark \ref{2.0'} and Section
4.4.1) and the LJ-cohomology and the
LJ-homology of $M$ are not dual one each other (see Theorems
\ref{5.2} and  \ref{trivial}).
}\end{remark}

\subsection{Lichnerowicz-Jacobi homology and conformal changes of
Jacobi structures}

In this section, we will show that the LJ-homology is also invariant under
conformal changes.

Suppose  that $(K,\lcf \; , \;\rcf, \varrho)$ (respectively, $(K',\lcf \; ,
\;\rcf', \varrho')$) is a Lie algebroid over $M$ of rank $n$ and that
$\phi:K\rightarrow K'$ is an isomorphism of Lie algebroids (see
Section 3.2). Denote by $\phi_1:\Gamma(K)\rightarrow \Gamma(K')$ the
isomorphism of $C^\infty(M,\R)$-modules induced by $\phi.$ This
isomorphism can be extended to an isomorphism
$\phi_r:\Gamma(\Lambda^rK)\rightarrow \Gamma(\Lambda^rK')$ by putting
\begin{equation}\label{4.13'''}
\phi_r(X_1\wedge \dots \wedge X_r)=\phi_1(X_1)\wedge \dots \wedge
\phi_1(X_r), 
\end{equation}
for all $X_1,\dots ,X_r\in \Gamma(K).$

Moreover, we have

\begin{proposition}\label{hcc1}
Let $\nabla$
(respectively, $\nabla'$) be a flat $K$-connection (respectively, a
$K'$-co\-nnec\-tion) on $\Lambda^nK\rightarrow M$ (respectively,
$\Lambda^nK'\rightarrow M$) such that $\nabla$ and $\nabla'$ are
$\phi$-related, that is,
\begin{equation}\label{13IV}
\phi_{n}(\nabla_X\Phi)=\nabla'_{\phi_1(X)}\phi_{n}(\Phi),
\end{equation}
for all $X\in \Gamma(K)$ and $\Phi\in \Gamma(\Lambda^nK).$ Then,
\begin{equation}\label{13V}
\phi_r\circ D=D'\circ \phi_{r+1},
\end{equation}
where $D$ (respectively, $D'$) is the homology operator associated
with $\nabla$ (respectively, $\nabla'$). Thus, the Lie algebroid
homologies $H_*(K,\nabla)$ and $H_*(K',\nabla')$
are isomorphic.
\end{proposition}
{\bf Proof:} Denote by $\phi^r:\Gamma(\Lambda^r(K')^*)\rightarrow
\Gamma(\Lambda^rK^*)$ the isomorphism of $C^\infty(M,\R)$-modules
given by (\ref{3.8'}). A direct computation proves that
\begin{equation}\label{13VI}
\phi_r(i(\phi^{n-r}\omega')\Phi)=i(\omega')(\phi_n\Phi),
\end{equation}
for all $\omega'\in \Gamma(\Lambda^{n-r}(K')^*)$ and $\Phi\in
\Gamma(\Lambda^nK)$. 

Therefore, using (\ref{3.8''}), (\ref{D}), (\ref{13IV}) and
(\ref{13VI}), we deduce (\ref{13V}).
\hfill$\Box$

Now, we have the following

\begin{theorem}\label{hcc2}
Let $(M,\Lambda,E)$ be a Jacobi manifold and $(\Lambda_a,E_a)$ a
conformal change of the Jacobi structure $(\Lambda,E)$. Then
\[
H_k^{LJ}(M,\Lambda,E)\cong H_k^{LJ}(M,\Lambda_a,E_a),
\]
for all $k.$
\end{theorem}
{\bf Proof:} We consider the isomorphism $\phi$
given by (\ref{3.9'}) between the Lie algebroids
$(T^{*}M \times \R, \lcf \; , \;\rcf_{(\Lambda, E)}, (\#_{\Lambda},
E))$ and $(T^{*}M \times \R, \lcf \; , \;\rcf_{(\Lambda_a,
E_a)}, (\#_{\Lambda_a}, E_a)).$

From (\ref{3.6'''}), (\ref{4.8'}) and (\ref{4.13'''}), we obtain that
\begin{equation}\label{13VII}
\phi_{n+1}(0,\Phi)=(0,\frac{1}{a^{n+1}}\Phi),
\end{equation}
for $\Phi\in \Omega^n(M),$ where $n$ is the dimension of $M.$

On the other hand, if $\nabla$ (respectively, $\nabla^a$) is the flat
$(T^*M\times \R)$-connection on $\Lambda^{n+1}(T^*M\times
\R)\rightarrow M$ defined by (\ref{conexion}) associated with the
Jacobi structure $(\Lambda,E)$ (respectively, $(\Lambda_a,E_a)$)
then, using (\ref{3.6'''}), (\ref{conexion}), (\ref{4.10''}) and
(\ref{13VII}), we prove that
\[
\phi_{n+1}(\nabla_{(\alpha,f)}(0,\Phi)) =
\nabla^a_{\phi_1(\alpha,f)}\phi_{n+1}(0,\Phi), 
\]
for all $(\alpha,f)\in \Omega^1(M)\times C^\infty(M,\R)$ and $\Phi\in
\Omega^n(M).$

Consequently, the result follows from 
Propositions \ref{a1} and \ref{hcc1}.
\hfill$\Box$

Finally, using Theorem \ref{hcc2}, we deduce the result announced at the
begining of this section

\begin{corollary}\label{hcc3}
The LJ-homology is invariant under conformal changes of the Jacobi
structure.
\end{corollary}

\subsection{Lichnerowicz-Jacobi homology of a
Poisson manifold}

Let $(M,\Lambda)$ be a Poisson manifold and $\delta$ (respectively,
$\bar\delta$) the operator of the LJ-homology
(respectively, of the  canonical homology) associated with $M$.
Then,
\begin{equation}\label{64'}
\begin{array}{rcl}
\bar\delta\alpha & = & i(\Lambda)d\alpha-di(\Lambda)\alpha\\[8pt]
\delta(\alpha,\beta) & = &
(\bar\delta\alpha,\bar\delta\beta+ (-1)^ki(\Lambda)\alpha)
\end{array}
\end{equation}
for all $\alpha\in \Omega^k(M)$ and  $\beta\in \Omega^{k-1}(M)$ (see
(\ref{63}) and Remark \ref{cano}).

Using  (\ref{64'}), we obtain the following relation between the
LJ-homology $H_*^{LJ}(M)$ and the  canonical homology $H_*^{can}(M).$

\begin{theorem}\label{4.3}
Let $(M,\Lambda)$ be a Poisson manifold. Suppose  that
$(0,Id_k):\Omega^{k-1}(M)\rightarrow \kern-4pt\Omega^k(M)\oplus
\Omega^{k-1}(M)$ and $(\pi_1)_k:\Omega^k(M)\oplus
\Omega^{k-1}(M)\rightarrow \Omega^k(M)$ are the  homomorphisms of
$C^\infty(M,\R)$-modules given by
\[
(0,Id_k)(\beta)=(0,\beta), \makebox[1cm]{}
(\pi_1)_k(\alpha,\beta)=\alpha
\]
for all $\alpha\in \Omega^k(M)$ and $\beta\in \Omega^{k-1}(M)$. Then:
\begin{enumerate}
\item
The mappings $(0,Id_k)$ and $(\pi_1)_k$ define an exact sequence
of  complexes
\[
0\longrightarrow (\Omega^{*-1}(M),\bar\delta)\stackrel{(0,Id)}
\longrightarrow (\Omega^{*}(M)\oplus \Omega^{*-1}(M), \delta)
\stackrel{\pi_1}\longrightarrow (\Omega^{*}(M),\bar\delta)
\longrightarrow 0
\]
\item The above exact sequence induces a long exact homology sequence
\[
\cdots \longrightarrow H_{k-1}^{can}(M) \stackrel{(0,Id_k)*}
\longrightarrow H_k^{LJ}(M)
\stackrel{((\pi_1)_k)_*}\longrightarrow H^{can}_{k}(M)
\stackrel{\Lambda_{k}}\longrightarrow H^{can}_{k-2}(M)
\longrightarrow \cdots
\]
where the connecting homomorphism  $\Lambda_k$ is defined by
\begin{equation}\label{64'''}
\Lambda_k[\alpha]=(-1)^k[i(\Lambda)(\alpha)]
\end{equation}
for all $[\alpha]\in H_k^{can}(M).$
\end{enumerate}
\end{theorem}

From Theorem \ref{4.3}, it follows that

\begin{corollary}\label{c4.5}
Let $M$ be a Poisson manifold such that its groups of canonical homology
have finite dimension. Then, the
LJ-homology groups have also finite dimension and
\[
H_k^{LJ}(M)\cong \frac{H_{k-1}^{can}(M)}{\mbox{\rm
Im}\Lambda_{k+1}}\oplus \ker \Lambda_k,
\]
where $\Lambda_r:H_r^{can}(M)\rightarrow H_{r-2}^{can}(M)$ is the
homomorphism given by (\ref{64'''}).
\end{corollary}

Next, we will obtain an explicit relation  between the
LJ-homology and the LJ-co\-ho\-mo\-lo\-gy for  the particular cases
of a  symplectic structure and of a  Lie-Poisson structure.

\subsubsection{Symplectic structures}

Let $(M,\Omega)$ be a symplectic manifold of dimension $2m$. If $f$
is a $C^\infty$ real-valued function on $M$, it follows that ${\cal
L}_{X_f}\Omega=0$ which implies that
\[
{\cal L}_{X_f}(\Omega\wedge \dots^{(m}\dots  \wedge \Omega)=0.
\]
Thus, $M$ is a unimodular Poisson manifold (see \cite{W2}). Using
this fact, (\ref{ljcohos}), Theorem \ref{a2} and Remark
\ref{r4.4}, we deduce the following result.

\begin{theorem}\label{4.4}
Let $(M,\Omega)$ be a symplectic manifold of  dimension $2m$.
Then
\[
H_k^{LJ}(M)\cong H_{LJ}^{2m-k+1}(M)
\]
for all $k.$ Moreover, if $M$ is of type finite, we have
\[
H^{LJ}_k(M)\cong \frac{H_{dR}^{2m-k+1}(M)}{{\rm
Im}L^{2m-k-1}}\oplus \ker L^{2m-k},
\]
where $H_{dR}^*(M)$ is the de Rham cohomology of $M$ and
$L^r:H_{dR}^r(M)\rightarrow H_{dR}^{r+2}(M)$ is the homomorphism
given by  (\ref{L}).
\end{theorem}

\subsubsection{Lie-Poisson structures}

Let $({\frak g},[\;,\;])$ be a real Lie algebra of dimension $n$
and consider on the dual space ${\frak g}^*$ the Lie-Poisson structure
$\bar \Lambda.$

Suppose that $\{\xi_i\}_{i=1,\dots ,n}$ is a basis of ${\frak g}$ and
that $(x_i)$ are the corresponding global coordinates for ${\frak
g}^*$. Denote by $\bar\nu$ the volume form on ${\frak g}^*$ given by
\[
\bar\nu=dx_1\wedge \dots \wedge dx_n
\]
and by $\mu_0$ the {\it modular character of ${\frak  g}$} , that is,
$\mu_0$ is the element of ${\frak g}^*$ defined by
\[
\mu_0(\xi)=\mbox{trace}(ad_\xi),\makebox[1cm]{} \mbox{ for all
}\xi \in {\frak g},
\]
where $ad_{\xi}:{\frak  g}\rightarrow {\frak g}$ is the endomorphism
given by
\[
ad_\xi(\eta)=[\xi,\eta],\makebox[1cm]{} \mbox{ for all } \eta\in
{\frak g}.
\]
$\mu_0$ induces a constant vector field on ${\frak g}^*$ which is
the modular vector field of the Poisson manifold $({\frak
g}^*,\bar\Lambda)$ with respect to the volume form $\bar\nu$ (see
\cite{Ko,W2}). Thus, if ${\frak g}$ is unimodular, i.e.,
if its modular character $\mu_0$ is zero then $({\frak
g}^*,\bar{\Lambda})$ is a unimodular Poisson manifold. Therefore,
from Theorem \ref{a2} and Remark \ref{r4.4}, we obtain

\begin{theorem}\label{ljauni}
Let $({\frak g},[\;,\;])$ be a unimodular real Lie algebra of dimension $n$
and  consider on  the dual space ${\frak g}^*$ the  Lie-Poisson
structure. Then, for all $k,$
\[
H^{LJ}_k({\frak g}^*)\cong H_{LJ}^{n-k+1}({\frak g}^*).
\]
\end{theorem}

Now, using (\ref{29'}), Theorem \ref{ljauni} and the fact that the
Lie  algebra of a compact Lie group is unimodular, we conclude
\begin{corollary}
Let ${\frak g}$ be the Lie algebra of a compact Lie
group of dimension $n$ and consider on the dual space ${\frak g}^*$
the Lie-Poisson structure. Then, for all $k,$
\[
H^{LJ}_k({\frak g}^*)\cong H_{LJ}^{n-k+1}({\frak g}^*)\cong
(H^{n-k+1}({\frak g})\otimes Inv)\oplus (H^{n-k}({\frak g})\otimes Inv),
\]
where $Inv$ is the algebra of Casimir functions on ${\frak g}^*$ and
$H^*({\frak g})$ is the cohomology of ${\frak g}$ relative to the
trivial representation of ${\frak  g}$ on $\R.$
\end{corollary}

\subsubsection{ A quadratic Poisson structure}
Let $\Lambda$ be the quadratic Poisson structure on $\R^2$ considered
in Section $3.3.3,$ that is, 
\[
\Lambda=xy\frac{\partial}{\partial x}\wedge \frac{\partial}{\partial y}.
\]
The modular vector field of $(\R^2,\Lambda)$ with respect to the
standar volumen $\nu=dx\wedge dy$ is 
\[
{\cal X}_\Lambda^{\nu}=x\frac{\partial}{\partial x}-y\frac{\partial
}{\partial y}.
\]
It is easy to confirm that $[{\cal  X}_\Lambda^{\nu}]\not=0$ in
$H_{LP}^1(\R^2,\Lambda).$ Thus, $(\R^2,\Lambda)$ is not a
unimodular Poisson manifold.

\medskip

{\bf  The canonical homology of
$(\R^2,\Lambda)$ }. First, we will compute
$H_2^{can}(\R^2,\Lambda).$ 

Suppose that $\beta=hdx\wedge dy$ is a $2$-form on $\R^2,$ with
$h\in C^\infty(\R^2,\R).$ We have that 
\begin{equation}\label{4.261}
\bar\delta\beta=-di(\Lambda)(\beta)=-d(xyh).
\end{equation}
Hence, we deduce that 
\[
\bar{\delta}\beta=0\Leftrightarrow xyh=cte\Leftrightarrow
h=0\Leftrightarrow \beta=0.
\]
Therefore,
\begin{equation}\label{4.262}
H_2^{can}(\R^2,\Lambda)=\{0\}.
\end{equation}
Now, let $\alpha$ be a $1$-form on $\R^2.$ It follows that
\[
\bar\delta\alpha=i(\Lambda)(d\alpha)=xyd\alpha(\frac{\partial}{\partial
x}, \frac{\partial}{\partial y}).
\]
Consequently, 
\[
\bar\delta\alpha=0\Leftrightarrow d\alpha=0\Leftrightarrow
\alpha=df,\makebox[.3cm]{} \mbox{ with } f\in C^\infty(\R^2,\R).
\]
This implies that (see (\ref{4.261})) 
\begin{equation}\label{4.263}
H^{can}_1(\R^2,\Lambda)=\frac{\{df/f\in
C^\infty(\R^2,\R)\}}{\bar\delta(\Omega^2(\R^2))}=\frac{\{df/f\in
C^\infty(\R^2,\R)\}}{\{d(xyh)/h\in C^\infty(\R^2,\R)\}}.
\end{equation}

In particular, we obtain that the dimension of $H^{can}_1(\R^2,\Lambda)$ is
not finite. In fact, if $n$ is an arbitrary integer, $n\geq 1$, and 
\[
\sum_{k=1}^n\lambda_k[dx^k]=0,\makebox[1cm]{} \mbox{ with }
\lambda_k\in \R,
\]
we have that 
\[
\sum_{k=0}^n\lambda_kx^k=xyh,
\]
where $\lambda_0\in \R$ and $h\in C^\infty(\R^2,\R).$ Thus, we
conclude that 
\[
\sum_{k=0}^n\lambda_kx^k=0,\makebox[1cm]{} \mbox{for all } x\in \R,
\]
and it follows that $\lambda_k=0,$ for all $k\in \{0,\dots ,n\}$
(note that $p(x)=\displaystyle\sum_{k=0}^n\lambda_kx^k$ is a polynomial of degree
$\leq n$).

Finally, we will compute $H_0^{can}(\R^2,\Lambda).$

If $\gamma=fdx+gdy$ is a $1$-form on $\R^2,$ we deduce that 
\[
\bar\delta\gamma=i(\Lambda)(d\gamma)=xy(\frac{\partial g}{\partial
x}-\frac{\partial f}{\partial y}),
\]
which implies that 
\[
H_{0}^{can}(\R^2,\Lambda)=\frac{C^\infty(\R^2)}{\{xy(\frac{\partial
g}{\partial x}-\frac{\partial f}{\partial y})/f,g\in C^\infty(\R^2,\R)\}}.
\]
On the other hand, using the fact that $H_{dR}^2(\R^2)=\{0\}$, we
obtain that 
\[
\{xy(\frac{\partial g}{\partial x}-\frac{\partial f}{\partial
y})/f,g\in C^\infty(\R^2,\R)\}=\{xyh/h\in C^\infty(\R^2,\R)\}
\]
and therefore,
\[
H_0^{can}(\R^2,\Lambda)=\frac{C^\infty(\R^2)}{\{xyh/h\in C^\infty(\R^2,\R)\}}.
\]
Now, we consider the $\R$-linear map 
\[
\psi:H_0^{can}(\R^2,\Lambda)\rightarrow
H^{can}_1(\R^2,\Lambda)\oplus \R
\]
defined by 
\[
\psi([f])=([df],f(0,0)),\makebox[1cm]{} \mbox{ for all } f\in C^\infty(\R^2,\R).
\]
An straightforward computation shows that $\psi$ is an isomorphism.
In fact, the $\R$-linear  map 
\[
\zeta:H^{can}_1(\R^2,\Lambda)\oplus \R\rightarrow H_0^{can}(\R^2,\Lambda)
\]
given by 
\[
\zeta([df],k)=[f-f(0,0)+k]
\]
is just the inverse of $\psi.$

Consequently,
\begin{equation}\label{4.264}
H_0^{can}(\R^2,\Lambda)\cong \frac{\{df/f\in
C^\infty(\R^2,\R)\}}{\{d(xyh)/h\in C^\infty(\R^2,\R)\}}\oplus \R.
\end{equation}

\begin{remark}
{\rm From (\ref{3.29'}), (\ref{4.262}), (\ref{4.263}) and (\ref{4.264}),
we conclude that 
\[
H^{can}_i(\R^2,\Lambda)\not\cong H_{LP}^{2-i}(\R^2,\Lambda),
\]
for $i\in \{0,1,2\}.$}
\end{remark}

\medskip

{\bf The LJ-homology of $(\R^2,\Lambda)$}. Using (\ref{4.262}),
(\ref{4.263}), (\ref{4.264}) and Theorem \ref{4.3}, we have that 
\begin{equation}\label{4.265}
\begin{array}{l}
H^{LJ}_3(\R^2,\Lambda,0)\cong
H_2^{can}(\R^2,\Lambda)=\{0\},\\
H_2^{LJ}(\R^2,\Lambda,0)\cong H_1^{can}(\R^2,\Lambda)=\displaystyle\frac{\{df/f\in
C^\infty(\R^2,\R)\}}{\{d(xyh)/h\in C^\infty(\R^2,\R)\}},\\
H_0^{LJ}(\R^2,\Lambda,0)\cong H_0^{can}(\R^2,\Lambda)\cong
\displaystyle\frac{\{df/f\in C^\infty(\R^2,\R)\}}{\{d(xyh)/h\in
C^\infty(\R^2,\R)\}}\oplus \R.
\end{array}
\end{equation}
Next, we will show that 
\[
H_1^{LJ}(\R^2,\Lambda,0)\cong H_1^{can}(\R^2,\Lambda)\oplus
H_0^{can}(\R^2,\Lambda).
\]
For this purpose, we consider the $\R$-linear map 
\[
\tilde\psi: H_1^{can}(\R^2,\Lambda)\oplus
H_0^{can}(\R^2,\Lambda)\rightarrow H_1^{LJ}(\R^2,\Lambda,0)
\]
given by 
\[
\tilde{\psi}([\alpha],[f])=[(\alpha,f)]
\]
for $\alpha\in \Omega^1(\R^2)$ and $f\in C^\infty(\R^2,\R),$ with
$\bar\delta\alpha=0.$ 

Note that if 
\[
\alpha'=\alpha+\bar\delta\beta,\makebox[1cm]{} f'=f+\bar\delta\gamma,
\]
with $\beta\in \Omega^2(\R^2)$ and $\gamma\in \Omega^1(\R^2)$ then,
since $H^2_{dR}(\R^2)=\{0\},$ there exists a $1$-form $\tilde{\gamma}$
on $\R^2$ satisfying 
\[
\beta=-d\tilde{\gamma}
\]
and 
\[
\delta(\beta,\gamma + \tilde\gamma)=(\bar\delta\beta,\bar\delta\gamma).
\]
Thus, 
\[
(\alpha',f')=(\alpha,f)+\delta(\beta,\gamma + \tilde\gamma).
\]
This proves that the map $\tilde\psi$ is well defined.

On the other hand, it is clear that $\tilde\psi$ is an epimorphism.
Moreover, using again that $H^2_{dR}(\R^2)=\{0\},$ it follows that
$\tilde\psi$ is a monomorphism. 

Therefore (see (\ref{4.263}) and (\ref{4.264})) 
\begin{equation}\label{4.266}
H_1^{LJ}(\R^2,\Lambda,0)\cong (\frac{\{df/f\in C^\infty(\R^2,\R)\}}{
\{d(xyh)/h\in C^\infty(\R^2,\R)\}})^2\oplus \R.
\end{equation}

\begin{remark}{\rm From (\ref{3.29''}), (\ref{4.265}) and
(\ref{4.266}), we conclude that 
\[
H_i^{LJ}(\R^2,\Lambda,0)\not\cong
H_{LJ}^{3-i}(\R^2,\Lambda,0),\makebox[1cm]{} \mbox{ for }i\in \{0,1,2,3\}.
\]}
\end{remark}

\subsection{Lichnerowicz-Jacobi homology of a
contact manifold}

Let $(M,\eta)$ be a contact manifold of dimension $2m+1$. Then
$\nu=\eta\wedge (d\eta)^m$ is a volumen form and (see \cite{Pva})
\[
{\cal M}_{(\Lambda,E)}^\nu=(-(m+1)E,0),
\]
where $(\Lambda,E)$ is the associated Jacobi structure on $M.$ 
Thus, $M$ is not a
unimodular Jacobi 
manifold and therefore it is not possible to apply Theorem
\ref{a2}. In fact, in this section, 
we will show that the LJ-homology of $M$ is
trivial. 

First, we prove the following result which will be useful in the
sequel.
\begin{lemma}\label{4.9}
Let $(M,\eta)$ be a contact manifold  of dimension $2m+1$ and
$(\Lambda,E)$ be the associated Jacobi structure on $M$. If
$e(\eta):\Omega^k(M)\rightarrow \Omega^{k+1}(M)$ and
$\tilde{L}^k:\Omega^k(M)\rightarrow \Omega^{k+2}(M)$ are the operators
defined by
\[
e(\eta)(\alpha)=\eta\wedge
\alpha,\makebox[1cm]{}\tilde{L}^k(\alpha)=\alpha\wedge d\eta,
\]
then,
\begin{equation}\label{con1}
i(\Lambda)\circ e(\eta)=e(\eta)\circ i(\Lambda),
\end{equation}
\begin{equation}\label{con2}
i(\Lambda)\circ {\tilde L}^k-{\tilde L}^{k-2} \circ
i(\Lambda)=(m-k)Id+e(\eta)\circ i_E,
\end{equation}
where $Id$ denotes the identity transformation.
\end{lemma}

{\bf Proof:} Let  $x$ be an arbitrary point of $M.$

Suppose that $(t,q^1,\dots ,q^m,p_1,\dots ,p_m)$ are canonical
coordinates in an open neighborhood $U$ of $x$ satisfying (\ref{6''}).

From (\ref{6''}),  it follows that
\[
(i(\Lambda)\circ e(\eta))(\alpha)=(e(\eta)\circ i(\Lambda))(\alpha),
\]
for all $\alpha\in \Omega^k(U).$
This proves (\ref{con1}).

On the other hand, if $\alpha_1,\dots ,\alpha_k$ are  $1$-forms on $U$
then a direct computation, using again  (\ref{6''}),
shows that
\[
\begin{array}{lcl}
(i(\Lambda)\circ \tilde{L}^k-\tilde{L}^{k-2}\circ i(\Lambda))(\alpha_1
\wedge \dots \wedge \alpha_k)&=&m(\alpha_1\wedge \dots \wedge \alpha_k)
- \displaystyle\sum_{i=1}^k\alpha_1\wedge \dots \wedge
\alpha_{i-1}\\[8pt]
&&\kern-76pt \wedge \;
[\displaystyle\sum_{j=1}^m(\alpha_i(\frac{\partial}{\partial q^j}
+ p_j\frac{\partial }{\partial t})dq^j+\alpha_i(\frac{\partial}{\partial
p_j}) dp_j)] \wedge \alpha_{i+1}\wedge \dots \wedge
\alpha_k \\[8pt]
&\kern-70pt=&\kern-40pt(m-k)(\alpha_1\wedge \dots \wedge \alpha_k)
+(e(\eta)\circ i_E)(\alpha_1 \wedge \dots \wedge\alpha_k).
\end{array}
\]
Hence, we have (\ref{con2}).
\hfill$\Box$

Now, we prove the following

\begin{theorem}\label{trivial}
The LJ-homology of a contact manifold  is trivial.
\end{theorem}

{\bf Proof:} Let  $(M,\eta)$ be a contact manifold of
dimension $2m+1$ and  $\delta$ the  LJ-homology operator.

If  $\alpha$ (respectively, $\beta$) is a $k$-form
(respectively, a $(k-1)$-form) on $M$ such that
\[
\delta(\alpha,\beta)=(0,0),
\]
then, using  (\ref{Reeb}), (\ref{63}) and Lemma \ref{4.9}, we deduce
that
\[
(\alpha,\beta)=\delta(\frac{\eta\wedge \alpha}{m+1},
\frac{\eta\wedge \beta}{m+1}).
\]
Therefore, $H_k^{LJ}(M)=\{0\}.$
\hfill$\Box$

\subsection{Lichnerowicz-Jacobi homology of a locally conformal
symplectic manifold} 

As in Section 3.5, we will distinguish the two following cases:

\medskip

$i)$ {\sc The particular case of a g.c.s. manifold:} Let $(M,\Omega)$
be a g.c.s. manifold with Lee $1$-form $\omega.$ Then, there exists
$f\in C^\infty(M,\R)$ such that $\omega=df$ and the Jacobi structure
of $M$ is a conformal change of the Poisson structure on $M$
associated with the symplectic form
$e^{-f}\Omega$ (see
Section 3.5). Thus, from Theorems \ref{t3.2'}, \ref{hcc2} and \ref{4.4}, we
conclude
\begin{theorem}\label{g4.4}
Let $(M,\Omega)$ be a g.c.s. manifold of  dimension $2m$.
Then,
\[
H_k^{LJ}(M)\cong H_{LJ}^{2m-k+1}(M)
\]
for all $k.$ Moreover, if $M$ is of type finite, we have 
\[
H^{LJ}_k(M)\cong \frac{H_{dR}^{2m-k+1}(M)}{{\rm
Im}\bar L^{2m-k-1}}\oplus \ker \bar L^{2m-k},
\]
where $H_{dR}^*(M)$ is the de Rham cohomology of $M$ and
$\bar L^r:H_{dR}^r(M)\rightarrow H_{dR}^{r+2}(M)$ is the homomorphism
given by
\[
\bar{L}^r[\alpha]=[e^{-f}\alpha\wedge \Omega],
\]
for all $[\alpha]\in H_{dR}^r(M).$
\end{theorem}

\begin{remark} {\rm In
\cite{Pva} Vaisman shows that a g.c.s. manifold is a
unimodular Jacobi manifold.
Using this fact, Theorem \ref{t3.8'} and Theorem \ref{a2}, we
also can prove  Theorem \ref{g4.4}.} 
\end{remark}

\begin{example}{\rm
Let $(N,\eta)$ be a contact manifold of type finite. Assume that the 
dimension of $N$ is $2m-1$
and consider on the product manifold $M=N\times \R$ the
g.c.s. structure $\Omega$ given by  (\ref{Exgcs}) (see Example
\ref{3.8'''}). Then,
 using  Theorem \ref{g4.4},  we have that
\[
H_k^{LJ}(M)\cong H_{dR}^{2m-k+1}(M)\oplus
H_{dR}^{2m-k}(M)\cong H_{dR}^{2m-k+1}(N)\oplus H_{dR}^{2m-k}(N).
\]
}
\end{example}

\medskip

$ii)$ {\sc The general case:} Now, we will study the LJ-homology of
an arbitrary l.c.s. manifold. 

Let $(M,\Omega)$ be a l.c.s. manifold of dimension $2m$ with Lee
$1$-form $\omega$. Then $\nu=\Omega^m$ is a volumen form and (see
\cite{Pva}) 
\[
{\cal M}_{(\Lambda,E)}^\nu=[(-(1+m)E,0)],
\]
where $(\Lambda,E)$ the associated
Jacobi structure on $M.$ Thus, in general, $M$ is not a unimodular
Jacobi manifold (see \cite{Pva}) and therefore it is not possible to
apply Theorem \ref{a2}. In fact, in this section, we will prove that
if $k\in \{0,\dots ,2m+1\}$ then, in general, the spaces
$H_{LJ}^k(M)$ and $H_{2m+1-k}^{LJ}(M),$ are not isomorphic.

First, we will introduce a certain cohomology in order to give an
explicit description of the LJ-homology of $M$.

We consider the  closed $1$-forms
$\omega_0$ and $\omega_1$ on $M$ defined by
\begin{equation}\label{omegas}
\omega_0=-m\omega,\makebox[1cm]{} \omega_1=-(m+1)\omega.
\end{equation}
Denote by $H_{\omega_0}^*(M)$ and  $H^*_{\omega_1}(M)$ the
cohomologies of the complexes $(\Omega^*(M),d_{\omega_0})$ and
$\kern-2pt (\Omega^*(\kern-1pt M),$ $d_{\omega_1}),$ where
$d_{\omega_0}$ and $d_{\omega_1}$ are the differential operators with
zero square given by 
(see (\ref{dw}) and (\ref{ew}))
\begin{equation}\label{domegas}
d_{\omega_0}=d+e(\omega_0),\makebox[1cm]{}d_{\omega_1}=d+e(\omega_1).
\end{equation}
Now, let $\tilde{d}:\Omega^k(M)\oplus \Omega^{k-1}(M)\rightarrow
\Omega^{k+1}(M)\oplus \Omega^k(M)$ be the differential operator defined
by
\begin{equation}\label{dtilde}
\tilde{d}(\alpha,\beta)=(d_{\omega_1}\alpha-\Omega\wedge
\beta,-d_{\omega_0}\beta).
\end{equation}
Using (\ref{conLcs}), it follows $\tilde{d}^2=0.$ Thus, we can
consider the complex
\[
\cdots \longrightarrow \Omega^{k-1}(M)\oplus \Omega^{k-2}(M)
\stackrel{\tilde{d}} \longrightarrow \Omega^k(M)\oplus \Omega^{k-1}(M)
\stackrel{\tilde{d}}\longrightarrow \Omega^{k+1}(M)\oplus \Omega^k(M)
\longrightarrow \cdots
\]
Denote by $\tilde{H}^*(M)$ the cohomology of this complex.

From (\ref{domegas}) and (\ref{dtilde}), we deduce the following
result which relates $\tilde{H}^*(M)$ with the cohomologies
$H_{\omega_0}^*(M)$ and $H_{\omega_1}^*(M).$

\begin{proposition}\label{4.11}
Let $(M,\Omega)$ be a l.c.s. manifold with Lee $1$-form $\omega$.
Suppose that  $(Id^k,0):\Omega^k(M)\rightarrow \Omega^k(M)\oplus
\Omega^{k-1}(M)$ and $(\pi_2)^k:\Omega^k(M)\oplus
\Omega^{k-1}(M)\rightarrow \Omega^{k-1}(M)$ are  the  homomorphism of
$C^\infty(M,\R)$-modules defined by
\[
(Id^k,0)(\alpha)=(\alpha,0),\makebox[1cm]{}(\pi_2)^k(\alpha,\beta) =
\beta ,
\]
for $\alpha\in \Omega^k(M)$ and  $\beta\in \Omega^{k-1}(M).$
Then:
\begin{enumerate}
\item
The mappings $(Id^k,0)$ and  $(\pi_2)^k$  induce an
exact sequence of complexes
\[
0\longrightarrow (\Omega^{*}(M),d_{\omega_1})\stackrel{(Id,0)}
\longrightarrow (\Omega^{*}(M)\oplus \Omega^{*-1}(M)
,\tilde{d})\stackrel{\pi_2}\longrightarrow (\Omega^{*-1}(M),-d_{\omega_0})
\longrightarrow 0.
\]
\item
This exact sequence induces  a long exact cohomology sequence
\[
\cdots \longrightarrow H^{k}_{\omega_1}(M) \stackrel{(Id^k,0)_*}
\longrightarrow \tilde{H}^k(M)
\stackrel{((\pi_2)^k)_*}\longrightarrow H_{\omega_0}^{k-1}(M)
\stackrel{-L^{k-1}}\longrightarrow H_{\omega_1}^{k+1}(M)
\longrightarrow \cdots  ,
\]
where the connector homomorphism $-L^{k-1}$ is defined by
\begin{equation}\label{Llcs}
(-L^{k-1})[\alpha]=[-\alpha\wedge \Omega],
\end{equation}
for all $[\alpha]\in H_{\omega_0}^{k-1}(M).$
\end{enumerate}
\end{proposition}

Now, from   Proposition \ref{4.11}, we obtain
\begin{corollary}\label{4.12}
Let $(M,\Omega)$ be a  l.c.s. manifold with Lee $1$-form  $\omega$ and
such that the cohomology groups  $H_{\omega_0}^k(M)$ and
$H_{\omega_1}^k(M)$ have  finite dimension, for all $k.$ Then,
the cohomology group  $\tilde{H}^k(M)$ has also finite dimension, for
all $k$, and
\[
\tilde{H}^k(M)\cong \frac{H_{\omega_1}^k(M)}{\mbox{\rm
Im}L^{k-2}}\oplus \ker L^{k-1},
\]
where $L^r:H_{\omega_0}^r(M)\rightarrow H^{r+2}_{\omega_1}(M)$ is the
homomorphism given by (\ref{Llcs}).
\end{corollary}

Using  (\ref{omegas}), (\ref{domegas}), Proposition  \ref{dw0}, Theorem
\ref{dw1} and  Corollary \ref{4.12}, we prove the following results:

\begin{corollary}\label{4.13}
Let $(M,\Omega)$ be a l.c.s. manifold with Lee $1$-form $\omega$ 
such that the dimensions of the cohomology groups
$H_{\omega_0}^k(M)$ and $H_{\omega_1}^k(M)$ are  finite, for all $k.$
Suppose that $\Omega$ is $d_{(-\omega)}$-exact, that is, there exists
a $1$-form $\eta$ on  $M$satisfying
\[
\Omega=d\eta-\omega\wedge \eta.
\]
Then, for all $k$, we have
\[
\tilde{H}^k(M)\cong H_{\omega_1}^k(M)\oplus H_{\omega_0}^{k-1}(M).
\]
\end{corollary}

\begin{corollary}\label{c4.14}
Let $(M,\Omega)$ be a compact l.c.s. manifold with  Lee $1$-form
$\omega\not=0.$ Suppose that $g$ is a Riemannian metric on $M$
such that  $\omega$ is parallel with respect to $g$.
Then, the  cohomology $\tilde{H}^*(M)$ is trivial.
\end{corollary}

Next, we will study the relation between the
LJ-homology of a l.c.s. manifold $M$ and the cohomology
$\tilde{H}^*(M).$

Let $(M,\Omega)$ be a l.c.s. manifold of  dimension $2m$ with
Lee $1$-form $\omega$. Denote by $(\Lambda,E)$ the associated Jacobi
structure on $M$ and by $\#_\Lambda:\Omega^k(M)\rightarrow {\cal
V}^k(M)$ the  isomorphism of $C^\infty(M,\R)$-mo\-du\-les given by
(\ref{e7}) and (\ref{e8}).

We define the star operator $\star:\Omega^k(M)\rightarrow
\Omega^{2m-k}(M)$ by
\begin{equation}\label{starlcs}
\star \alpha=(-1)^ki(\#_\Lambda(\alpha))\frac{\Omega^m}{m!}
\end{equation}
for all $\alpha\in \Omega^k(M).$

\begin{lemma}\label{4.15}
If $\alpha$ is a $k$-form on  $M$ then:
\begin{enumerate}
\item[$(i)$]
$\star(\star(\alpha))=\alpha.$
\item[$(ii)$]
$(i_E\circ \star )(\alpha)=(-1)^k(\star\circ e(\omega))(\alpha).$
\item[$(iii)$]
$({\cal L}_E\circ \star )(\alpha)=(\star\circ {\cal L}_E)(\alpha).$
\item[$(iv)$]
$(i(\Lambda)\circ \star )(\alpha)=\star (\alpha\wedge \Omega).$
\end{enumerate}
\end{lemma}
{\bf Proof:}
$(i)$ Let $(q^1,\dots, q^m,p_1,\dots ,p_m)$ be coordinates on  an
open subset $U$ of $M$ such that
\begin{equation}\label{local}
\omega=df, \makebox[1cm]{}\Omega=e^{f}\sum_{i}dq^i\wedge
dp_i, \makebox[1cm]{}\Lambda=e^{-f}\sum_i(\frac{\partial}{\partial
q^i}\wedge \frac{\partial}{\partial p_i}),
\end{equation}
where $f:U\rightarrow \R$ is a $C^\infty$ real-valued function on
$U$ (see Section 2.2).

Consider on $U$ the symplectic $2$-form
\begin{equation}\label{25}
\bar\Omega=e^{-f}\Omega=\sum_{i}dq^i\wedge dp_i.
\end{equation}
Denote by  $\bar\star:\Omega^k(U)\rightarrow \Omega^{2m-k}(U)$
the star operator induced by $\bar{\Omega},$ that is,
\[
\bar \star
\alpha=(-1)^ki(\#_{\bar\Lambda}(\alpha))\frac{\bar\Omega^m}{m!}, 
\]
where $\bar\Lambda$ is the Poisson structure associated with
$\bar\Omega$. 

In \cite{Br} (see also \cite{Li}), Brylinski shows that
\begin{equation}\label{69'}
\bar\star^2=Id.
\end{equation}
On the other hand, using
(\ref{e7}), (\ref{e8}), (\ref{starlcs}), (\ref{local}) and
(\ref{25}), we have that 
\begin{equation}\label{Rel}
\star(\alpha)=e^{(m-k)\sigma}\bar\star (\alpha),
\makebox[1cm]{}\mbox{ for all } \alpha \in \Omega^k(U).
\end{equation}
Thus, from (\ref{69'}) and (\ref{Rel}),
it follows that $\star^2(\alpha)=\alpha$, for all $\alpha\in
\Omega^k(U)$. This proves  $(i)$.

Using (\ref{Lcs}), (\ref{sb}) and (\ref{starlcs}), we deduce $(ii).$

From (\ref{Inv}), (\ref{starlcs}) and the first relation of
(\ref{s-dl}), we deduce that $(iii)$ holds.

Finally, $(iv)$ follows directly using (\ref{starlcs}) and the fact that
\[
\#_\Lambda(\Omega)=\Lambda.
\]
\hfill$\Box$

The star operator $\bar{\star}$ induced by a symplectic form
$\bar{\Omega}$ on a manifold $M$ satisfies the following relation
\[
\bar\star(\bar\delta\alpha)=(-1)^{k+1}d(\bar\star(\alpha)),
\]
for all $\alpha\in \Omega^k(M),$ where
$\bar\delta=i(\bar\Lambda)\circ d-d\circ i(\bar\Lambda)$ is the
canonical homology operator (see \cite{Br}). Therefore, proceeding as
in the proof of the first part of Lemma
\ref{4.15} and  using this lemma, we obtain

\begin{lemma}\label{4.16}
If $\alpha$ is a $k$-form on $M$ then
\[
\star(i(\Lambda) d\alpha - di(\Lambda) \alpha)
= (-1)^{k+1}(d (\star \alpha) - (m-k+1) e(\omega) (\star \alpha) -
\Omega \wedge i_E (\star \alpha)).
\]
\end{lemma}

Now, from  Lemmas \ref{4.15} and \ref{4.16}, we deduce the following

\begin{theorem}\label{4.17}
Let $(M,\Omega)$ be a l.c.s. manifold of dimension $2m$ and with
Lee $1$-form $\omega.$ Suppose that $(\Lambda,E)$ is the
associated Jacobi structure on $M$ and that
$\tilde{\phi}_k:\Omega^k(M)\oplus \Omega^{k-1}(M)\rightarrow
\Omega^{2m+1-k}(M)\oplus \Omega^{2m-k}(M)$ is the isomorphism of
$C^\infty(M,\R)$-modules defined by
\[
\tilde{\phi}_k(\alpha,\beta)=(\star \beta, i_E(\star \beta)-\star\alpha),
\]
where $\star$ is the star operator given by (\ref{starlcs}). If
$\delta$ is the   LJ-homology operator of $M$ and
$\tilde{d}$ is the differential  operator defined by (\ref{dtilde})
then,
\[
\tilde{\phi}_{k-1}(\delta(\alpha,\beta))=
(-1)^{k}\tilde{d}(\tilde{\phi}_k(\alpha,\beta))
\]
for all $(\alpha,\beta)\in \Omega^k(M)\oplus \Omega^{k-1}(M).$ Thus,
\[
H_k^{LJ}(M)\cong \tilde{H}^{2m+1-k}(M).
\]
\end{theorem}

Using  Theorem \ref{4.17} and Corollaries \ref{4.12},
\ref{4.13} and
\ref{c4.14}, we conclude
\begin{corollary}
Let $(M,\Omega)$ be  a l.c.s. manifold of dimension $2m,$ with
Lee $1$-form $\omega$ and  such that the cohomology groups
$H_{\omega_0}^k(M)$ and $H^k_{\omega_1}(M)$ have finite dimension,
for all $k$. Then:
\[
H^{LJ}_k(M)\cong \frac{H_{\omega_1}^{2m+1-k}(M)}{\mbox{\rm
Im}L^{2m-k-1}}\oplus \ker L^{2m-k},
\]
where $L^r:H_{\omega_0}^r(M)\rightarrow H_{\omega_1}^{r+2}(M)$ is the
homomorphism given by (\ref{Llcs}).
\end{corollary}

\begin{corollary}\label{4.19}
Let $(M,\Omega)$ be a l.c.s. manifold of  dimension $2m$, with
Lee $1$-form $\omega$ and such that the dimensions of the cohomology
groups $H_{\omega_0}^k(M)$ and $H_{\omega_1}^k(M)$ are
finite, for all  $k.$ Suppose that $\Omega$ is
$d_{(-\omega)}$-exact, that is, there exits a $1$-form $\eta$ on
$M$ which satisfies
\[
\Omega=d\eta-\omega\wedge \eta.
\]
Then,
\[
H^{LJ}_k(M)\cong H_{\omega_1}^{2m+1-k}(M)\oplus H_{\omega_0}^{2m-k}(M),
\]
for all $k.$
\end{corollary}

\begin{corollary}\label{4.20}
Let $(M,\Omega)$ be a compact l.c.s. manifold with Lee $1$-form
$\omega,$ 
$\omega\not=0.$ Suppose that $g$ is a Riemannian metric on $M$ such
that $\omega$ is parallel with respect to $g$. Then, the LJ-homology of
$M$ is trivial. 
\end{corollary}

\begin{remark}{\rm Note that under the same hypotheses as in
Corollary \ref{4.20}, we have that $H_{LJ}^k(M)\cong H_{dR}^k(M),$
for all $k$.}
\end{remark}
\begin{example}
{\rm Let $(N,\eta)$ be a  compact contact manifold of
dimension $2m-1.$
We consider on the product manifold
$M=N\times S^1$ the l.c.s. structure $\Omega$ given by (\ref{Exlcs})
(see  Example  \ref{comcon}). Then, from  Corollary \ref{4.20}, it
follows that the LJ-homology
of $M$ is trivial. 
}
\end{example}

\subsection{Lichnerowicz-Jacobi homology of the unit sphere of a
real Lie algebra}

In this section we will give an explicit description of the LJ-homology of
the unit sphere on the Lie algebra of a compact Lie group. For this
purpose, we will prove that the unit sphere of a unimodular real Lie
algebra is a unimodular Jacobi manifold.

\begin{theorem}\label{homesf}
Let ${\frak g}$ be a  unimodular real Lie algebra of dimension $n.$
Suppose that $<\;,\;>$ is a
scalar product on ${\frak g}$ and consider on the unit sphere
$S^{n-1}({\frak g})$ the induced Jacobi structure.
Then $S^{n-1}({\frak g})$ is a unimodular Jacobi manifold. Thus,
\begin{equation}\label{LJcanS}
H_k^{LJ}(S^{n-1}({\frak g}))\cong
H_{LJ}^{n-k}(S^{n-1}({\frak g})),
\end{equation}
for all $k$.
\end{theorem}
{\bf Proof:} Let $(x^{i})_{i=1,\dots,n}$ be the global coordinates
for ${\frak  g}$ obtained from an orthonormal basis
$\{\xi_{i}\}_{i=1,\dots,n}$ of ${\frak  g}$.

Consider the volume form $\bar\nu$ on ${\frak  g}$ defined by
\[
\bar \nu = dx^{1} \wedge \dots \wedge dx^{n}.
\]
Denote by $\flat_{<\;,\;>}: {\frak g} \longrightarrow {\frak g}^{*}$
the linear isomorphism between ${\frak g}$ and ${\frak g}^{*}$ given
by (\ref{2.12'}) and by $\bar{\Lambda}$ the Poisson structure on ${\frak
g}$ induced by the Lie-Poisson structure on ${\frak g}^{*}$ and by
the isomorphism $\flat_{<\;,\;>}$.

Suppose that $(x_{i})_{i=1,\dots,n}$ are the global coordinates for
${\frak g}^{*}$ obtained from the basis $\{\xi_{i}\}_{i=1,\dots,n}$.
Since $\flat_{<\;,\;>}^{*}(dx_{1} \wedge \dots \wedge dx_{n}) =
\bar\nu$, it follows that the modular vector field ${\cal
X}^{\bar \nu}_{\bar{\Lambda}}$ of $({\frak g}, \bar{\Lambda})$
with respect to $\bar\nu$ is zero (see Section 4.4.2).

Now, consider the $(n-1)$-form $\nu$ on $S^{n-1}({\frak g})$ defined by
\[
\nu = \displaystyle \sum_{i=1}^{n} (-1)^{n-i} <\xi_{i}, \;\;>
d<\xi_{1}, \;\;> \wedge \dots \wedge \widehat{d<\xi_{i}, \;\;>} \wedge \dots
\wedge d<\xi_{n}, \;\;>,
\]
where $<\xi_{j}, \;\;>: S^{n-1}({\frak g}) \longrightarrow \R$ is the
real function given by (\ref{Fusc}). Then, if $F: {\frak g} - \{0\}
\longrightarrow S^{n-1}({\frak g}) \times \R$ is the diffeomorphism
defined by (\ref{2.17'}), a direct computation proves that
\begin{equation}\label{volgs}
(F^{-1})^{*}({\bar\nu}_{|{\frak g} - \{0\}})=
e^{nt} \nu \wedge dt.
\end{equation}
Therefore, $\nu$ is a volume form on $S^{n-1}({\frak g})$ and, using
(\ref{volgs}), Remarks \ref{2.0''} and \ref{r4.5} and the fact that
the modular vector field of $({\frak g}-\{0\},\bar{\Lambda}_{|{\frak
g}-\{0\}})$ with respect to $\bar{\nu}_{|{\frak g}-\{0\}}$ is zero, we
deduce that
\[
{\cal M}_{(\Lambda,E)}^{\nu} = (0, 0),
\]
$(\Lambda, E)$ being the Jacobi structure of $S^{n-1}({\frak g})$.
\hfill$\Box$

From Theorems \ref{3.18} and \ref{homesf}, we conclude

\begin{theorem}
Let ${\frak g}$ be the Lie algebra of a compact Lie group
of dimension $n.$ Suppose that  $<\;,\;>$ is a
scalar product on ${\frak g}$ and consider on the unit sphere
$S^{n-1}({\frak g})$ the induced Jacobi structure.
Then,
\begin{equation}
H_k^{LJ}(S^{n-1}({\frak g}))\cong
H_{LJ}^{n-k}(S^{n-1}({\frak g}))\cong H^{n-k}({\frak g})\otimes Inv
\end{equation}
for all $k$, where  $H^*({\frak g})$ is the
cohomology of ${\frak g}$ relative to the  trivial representation of
${\frak g}$ on  $\R$ and  $Inv$ is the subalgebra of 
$C^\infty(S^{n-1}({\frak g}),\R)$ defined by
\[
Inv=\{\varphi\in C^\infty(S^{n-1}({\frak
g}),\R)/X_f(\varphi)=0,\;\;\forall f\in C^\infty(S^{n-1}({\frak g}),\R)\}.
\]
\end{theorem}

\subsection{Table II}

The following table summarizes the main results obtained in Sections
4.4, 4.5, 4.6 and  4.7 on the LJ-homology
and its relation with the LJ-cohomology  of the different types of
Jacobi manifolds. 

\vspace{20pt}

\hspace{-40pt}
{\footnotesize
\begin{tabular}{|c|c|c|c|}
\hline\hline
&&MODULAR&\\[-3pt]
TYPE&LJ-HOMOLOGY& CLASS & REMARKS\\[5pt]
\hline\hline &&&\\[-5pt]
$(M^{2m},\Omega)$ symplectic&$H_k^{LJ}(M)\cong
H_{LJ}^{2m+1-k}(M)$&$0$&$L^r:H_{dR}^r(M)\rightarrow H_{dR}^{r+2}(M)$
\\
of finite type& $\cong
\displaystyle\frac{H_{dR}^{2m+1-k}(M)}{\mbox{Im}L^{2m-k-1}}\oplus
\ker L^{2m-k}$&&$[\alpha]\mapsto [\alpha\wedge \Omega]$\\[5pt]
\hline\hline&&&\\[-5pt]
Dual of a unimodular &&&
\\
real Lie algebra ${\frak g}$&$H_k^{LJ}({\frak
g}^*)\cong 
H_{LJ}^{n-k+1}({\frak g}^*)$&$0$&\\[3pt]
of dimension $n$ &&&\\[3pt]
\hline&&&\\[-5pt]
Dual of the Lie algebra ${\frak g}$& $H_k^{LJ}({\frak g}^*)\cong
H_{LJ}^{n-k+1}({\frak g}^*)\cong$&&$Inv\equiv$ subalgebra of the \\
of a compact Lie group&$(H^{n-k+1}({\frak g})\otimes Inv)\oplus$&$0$&
Casimir functions of ${\frak g}^*$ \\
$(\dim {\frak g}=n)$&$\oplus (H^{n-k}({\frak g})\otimes Inv)$&&\\[3pt]
\hline\hline&&&\\[-5pt]
$M^{2m+1}$ contact& $H^{LJ}_k(M)=\{0\}$&$[-(m+1)E,0)]\not=0$&
$E\equiv$Reeb vector field\\[3pt]
\hline\hline&&&\\[-6pt]
$(M^{2m},\Omega)$ g.c.s. of finite type 
&$H_k^{LJ}(M)\cong$
&$0$&$\bar{L}^r:H_{dR}^r(M)\rightarrow H_{dR}^{r+2}(M)$\\ 
with Lee $1$-form $\omega=df$ 
&$\cong\displaystyle\frac{H_{dR}^{2m+1-k}(M)}{\mbox{Im}\bar{L}^{2m-k-1}}\oplus
\ker \bar{L}^{2m-k}$&&$[\alpha]\mapsto [e^{-f}\alpha\wedge 
\Omega]$\\[8pt]
\hline\hline&&&\\[-3pt]
$(M^{2m},\Omega)$ l.c.s. with Lee $1$-form 
&$H_k^{LJ}(M)\cong$&&$L^r:H_{\omega_0}^r(M)\rightarrow H_{\omega_1}^{r+2}(M)$\\
$\omega$ and $dim H_{\omega_i}^*(M) < \infty$ $(i=0,1)$
&$\cong
\displaystyle\frac{H_{\omega_1}^{2m+1-k}(M)}{\mbox{Im}L^{2m-k-1}}\oplus
\ker L^{2m-k}$& $[(-(1+m)E,0)]$&$[\alpha]\mapsto [\alpha\wedge 
\Omega]$\\
$\omega_0=-m\omega,\;\;\omega_1=-(m+1)\omega$&&&\\
\hline&&&\\[-7pt]
$M^{2m}$ compact l.c.s. & &&\\
with Lee $1$-form $\omega$ & $H_k^{LJ}(M)=\{0\}$&$[(-(1+m)E,0)]\not=0$&\\
$\omega$ parallel with respect to &&&\\
a Riemannian metric &&&\\[3pt]
\hline\hline&&&\\[-5pt]
Unit sphere $S^{n-1}({\frak g})$ &$H_k^{LJ}(S^{n-1}({\frak
g}))\cong$&&$Inv\equiv$
subalgebra of the \\
of the Lie algebra ${\frak g}$ of &$\cong H_{LJ}^{n-k}(S^{n-1}({\frak
g}))\cong$&$0$& constant functions on  \\ 
a compact Lie group & $H^{n-k}({\frak g})\otimes Inv$&&the leaves of 
\\
$(\dim{\frak g}=n)$&& &the characteristic foliation \\[3pt]
\hline
\end{tabular}
}\vspace{-10pt}
\begin{center}
{\it Table II: LJ-homology }
\end{center}

\bigskip

\section*{Acknowledgments}
This work has been partially supported through grants DGICYT
(Spain) Projects PB97-1257 and PB97-1487. ML
wishes to express his gratitude for the 
hospitality offered to him in the Departamento de Matem\'atica
Fundamental (University of La Laguna) where part of this work was
conceived.

\end{document}